\documentclass[a4paper,10pt,twoside]{article}
 
\newif\ifdraft

\drafttrue
\draftfalse
 
\usepackage{amscd}
\usepackage{amssymb}
\usepackage{amsmath}
 
\newtheorem{theorem}{\bf Theorem}[section]

\newtheorem{assumption}[theorem]{\bf Hypothesis}
\newtheorem{remark}[theorem]{\bf Remark}
\newtheorem{proposition}[theorem]{\bf Proposition}
\newenvironment{proof}{\noindent{\bf Proof.}}{\hbox{\, }~\kern-10pt~\hfill{\rule[-1pt]{5pt}{10pt}}}
\newtheorem{definition}[theorem]{\bf Definition}
\newtheorem{corollary}[theorem]{\bf Corollary}

\renewcommand\theequation{\thesection.\arabic{equation}}
\def\ssection#1{
\setcounter{theorem}{0}
\setcounter{subsection}{0}\setcounter{equation}{0}\addtocounter{section}{1}\pagebreak[2]\leavevmode\section*{\normalsize\bf\thesection. #1\hfil}\leavevmode\vskip-\parskip\nopagebreak}

\def\ssubsection#1{\addtocounter{subsection}{1}\pagebreak[2]\leavevmode\subsection*{\normalsize\bf\thesubsection. #1\hfil}\leavevmode\vskip-\parskip\nopagebreak}

\def\centreline{\centerline}

\newcounter{appendix}
\def\aappendix#1{\setcounter{theorem}{0}\let\ara=\Alph%
\setcounter{equation}{0}\addtocounter{appendix}{1}\setcounter{section}{\theappendix}
\pagebreak[2]\leavevmode\section*{\normalsize\bf Appendix \Alph{appendix}. #1\hfil}\leavevmode\vskip-\parskip\nopagebreak\renewcommand\theequation{\ara{appendix}.\arabic{equation}}
\renewcommand\thetheorem{\Alph{appendix}.\arabic{theorem}}}

\renewenvironment{abstract}{\small\bf Abstract\rm
        \quotation}

\global\evensidemargin=0in

\newtoks\ara
\let\ara=\arabic
\ifdraft
\gdef\BIBITEM[#1]#2{\bibitem[#1]{#2}{{\mbox{\small[{\bf\tt#2}]}}} }
\def\LABEL#1{\label{#1}{\ifmmode%
\gdef\theequation{\mbox{\small{%
\hbox{$\bf\tt#1$}\rlap{\mbox{${\ \ (\rm\ara{section}.\arabic{equation})}$}}%
}}\gdef\theequation{\ara{section}.\arabic{equation}}}%
\else\mbox{\small$({\bf\tt#1})$}\fi}%
}\else\def\LABEL#1{\label{#1}}\gdef\BIBITEM{\bibitem}\fi

\gdef\BIBITEM[#1]#2{\bibitem{#2}}

\long\def\new#1\endnew{{#1}}
\long\def\old#1\endold{}

\begin{document}

 \def\Eq.#1.{Eq.~(\ref{eq:#1})}
\def\Ineq.#1.{Ineq.~(\ref{eq:#1})}
\def\be{\begin{equation}}
\def\ee{\end{equation}}
\def\D{\b\nabla}
\def\dddot{\kern1pt\ddot{\ }\dot{\ }\kern-7pt}
\def\beq{\begin{eqnarray}}
\def\eq{\end{eqnarray}}
\def\frc{\displaystyle\frac}
\def\Int{\oint_{S_\infty}\kern-5pt}
\def\INT{\displaystyle\int}
\def\({\left(}
\def\){\right)}
\def\ni{\noindent}
\def\Var{\mathop{\rm Var}\nolimits}
\def\as#1->#2{\qquad\mbox{as $#1\to#2$}}
\def\Res{\mathop{\rm Res}\nolimits}
\def\d{{\boldmath\mbox{$d$}}}
\def\m{{\boldmath\mbox{$m$}}}
\def\ip<#1,#2>{\langle#1,#2\rangle}
\def\cosec{\mathop{\rm cosec}\nolimits}
\def\Sum{\displaystyle\sum}
\def\h{\mbox{$\frac12$}}
\def\bimath{\b{\imath}}
\def\bjmath{\b{\jmath}}
\def\adj{\mathop{\rm adj}\nolimits}
\def\and{\qquad\mbox{and}\qquad}
\def\tilde{\widetilde}
\def\<{\langle}
\def\>{\rangle}
\def\up{\ket|\kern-3pt\uparrow>}
\def\down{\ket|\kern-3pt\downarrow>}
\def\Dt{{\frc {\d}{\d t}}}
\def\frcn#1#2{\mbox{$\frac{#1}{#2}$}}
\def\twid#1{\smash{\mbox{${\mathrel{\sim}}^{\kern-7pt\strut#1}$}\,}}
\def\sup{\mathop{\rm sup}\nolimits}
\def\inf{\mathop{\rm inf}\nolimits}
\def\arg{\mathop{\rm arg}\nolimits}
\def\erf{\mathop{\rm erf}\nolimits}
\def\erfc{\mathop{\rm erfc}\nolimits}
\def\diag{\mathop{\rm diag}\nolimits}
\def\sech{\mathop{\rm sech}\nolimits}
\def\tanh{\mathop{\rm tanh}\nolimits}
\def\cosech{\mathop{\rm cosech}\nolimits}
\def\coth{\mathop{\rm coth}\nolimits}
\def\ket|#1>{|#1\>}
\def\bra<#1|{\<#1|}
\def\IP<#1|#2>{\<#1|#2\>}
\newcount\questionnumber
\questionnumber=1
\def\Q#1.{\bigskip\noindent{\bf #1. }}
\def\s#1,#2.{\left#1,#2\right}
\def\Stab{\mathop{\rm Stab}\nolimits}
\def\diam{\mathop{\rm diam}\nolimits}
\def\tr{\mathop{\rm tr}\nolimits}
\def\Ff{{\frak f}}
\def\BA{{\Bbb A}}
\def\BB{{\Bbb B}}
\def\BC{{\Bbb C}}
\def\BD{{\Bbb D}}
\def\BE{{\Bbb E}}
\def\BF{{\Bbb F}}
\def\BG{{\Bbb G}}
\def\BH{{\Bbb H}}
\def\BI{{\Bbb I}}
\def\BJ{{\Bbb J}}
\def\BK{{\Bbb K}}
\def\BL{{\Bbb L}}
\def\BM{{\Bbb M}}
\def\BN{{\Bbb N}}
\def\BO{{\Bbb O}}
\def\BP{{\Bbb P}}
\def\BQ{{\Bbb Q}}
\def\BR{{\Bbb R}}
\def\BS{{\Bbb S}}
\def\BT{{\Bbb T}}
\def\BU{{\Bbb U}}
\def\BV{{\Bbb V}}
\def\BW{{\Bbb W}}
\def\BX{{\Bbb X}}
\def\BY{{\Bbb Y}}
\def\BZ{{\Bbb Z}}
\def\CD{{\cal D}}
\def\ca{\cos\alpha}\def\sa{\sin\alpha}
\def\cg{\cos\gamma}\def\sg{\sin\gamma}
\def\Re{\mathop{\rm Re}\nolimits}
\def\mtx#1,#2,#3,#4.{\pmatrix{#1&#2\cr#3&#4}}\def\ge{\geqslant}
\def\le{\leqslant}
\def\hat{\widehat}
\def\emptyset{\varnothing}
\def\bar{\overline}
\def\mob#1,#2,#3,#4[#5]{\frc{#1#5+#2}{#3#5+#4}}
\def\Im{\mathop{\rm Im}\nolimits}
\def\Ker{\mathop{\rm Ker}\nolimits}
\def\n{\mathop{\rm n}\nolimits}
\def\r{\mathop{\rm r}\nolimits}
\def\p#1#2.{\frc{\partial#2}{\partial x_{#1}}}
\def\curl{\D\times}
\def\div{\D\cdot}
\def\dS{\cdot d\bS}
\def\pp#1,#2{\partial#1/\partial#2}
\def\ppp#1,#2{\frc{\partial#1}{\partial x_{#2}}}
\def\pd#1,#2{\frc{\partial#1}{\partial#2}}
\def\td#1,#2{\frc{\d#1}{\d#2}}
\def\pdc#1,#2,#3{\left.\pd#1,#2\right|_{#3}\kern-6pt}
\def\cross{\times}
\def\kmh{$\,\rm kmh^{-1}$}
\def\sign{\mathop{\rm sign}\nolimits}
\def\imp{\Rightarrow}
\def\impby{\Leftarrow}
\def\nimp{\not\Rightarrow}
\def\nimpby{\kern3pt\not\kern-3pt\Leftarrow}
\def\limsup{\mathop{\rm limsup}\limits}
\def\liminf{\mathop{\rm liminf}\limits}
\def\cupp{\mathop{\displaystyle\cup}\limits}
\def\capp{\mathop{\displaystyle\cap}\limits}
\def\openone{\leavevmode\hbox{\rm\small1\kern-3.8pt\normalsize1}}
\def\max{\mathop{\rm max} \limits}
\def\const{\mathop{\rm constant}\nolimits}
\def\Tr{\mathop{\rm Tr}\nolimits}

\def\cl#1,#2.{\left[\matrix{#1\cr#2}\right]}
\def\det{\mathop{\rm det}\nolimits}
\def\sup{\mathop{\rm sup}\limits}
\def\inf{\mathop{\rm inf}\limits}
\def\arg{\mathop{\rm arg}\nolimits}
\newcount\questionnumber
\questionnumber=1
\def\Q#1.{\bigskip\noindent{\bf #1. }}
\def\s#1,#2.{\left#1,#2\right}
\def\CA{{\cal A}}
\def\CB{{\cal B}}
\def\CC{{\cal C}}
\def\CD{{\cal D}}
\def\CE{{\cal E}}
\def\CF{{\cal F}}
\def\CG{{\cal G}}
\def\CH{{\cal H}}
\def\CI{{\cal I}}
\def\CJ{{\cal J}}
\def\CK{{\cal K}}
\def\CL{{\cal L}}
\def\CM{{\cal M}}
\def\CN{{\cal N}}
\def\CO{{\cal O}}
\def\CP{{\cal P}}
\def\CQ{{\cal Q}}
\def\CR{{\cal R}}
\def\CS{{\cal S}}
\def\CT{{\cal T}}
\def\CU{{\cal U}}
\def\CV{{\cal V}}
\def\CW{{\cal W}}
\def\CX{{\cal X}}
\def\CY{{\cal Y}}
\def\CZ{{\cal Z}}
\def\FA{{\frak A}}
\def\FB{{\frak B}}
\def\FC{{\frak C}}
\def\FD{{\frak D}}
\def\FE{{\frak E}}
\def\FF{{\frak F}}
\def\FG{{\frak G}}
\def\FH{{\frak H}}
\def\FI{{\frak I}}
\def\FJ{{\frak J}}
\def\FK{{\frak K}}
\def\FL{{\frak L}}
\def\FM{{\frak M}}
\def\FN{{\frak N}}
\def\FO{{\frak O}}
\def\FP{{\frak P}}
\def\FQ{{\frak Q}}
\def\FR{{\frak R}}
\def\FS{{\frak S}}
\def\FT{{\frak T}}
\def\FU{{\frak U}}
\def\FV{{\frak V}}
\def\FW{{\frak W}}
\def\FX{{\frak X}}
\def\FY{{\frak Y}}
\def\FZ{{\frak Z}}
\def\Fa{{\frak a}}
\def\Fb{{\frak b}}
\def\Fc{{\frak c}}
\def\Fd{{\frak d}}
\def\Fe{{\frak e}}
\def\Ff{{\frak f}}
\def\Fg{{\frak g}}
\def\Fh{{\frak h}}
\def\Fi{{\frak i}}
\def\Fj{{\frak j}}
\def\Fk{{\frak k}}
\def\Fl{{\frak l}}
\def\Fm{{\frak m}}
\def\Fn{{\frak n}}
\def\Fo{{\frak o}}
\def\Fp{{\frak p}}
\def\Fq{{\frak q}}
\def\Fr{{\frak r}}
\def\Fs{{\frak s}}
\def\Ft{{\frak t}}
\def\Fu{{\frak u}}
\def\Fv{{\frak v}}
\def\Fw{{\frak w}}
\def\Fx{{\frak x}}
\def\Fy{{\frak y}}
\def\Fz{{\frak z}}
\def\smilie{\leavevmode\raise1pt\hbox{$\bigcirc$\kern-8pt\raise1.5pt\hbox{$\cdot\kern0.5pt\cdot$}\kern-6.35pt\lower3.5pt\hbox{${}^{{}_\smile}$}\kern-6.2pt\raise3.3pt\hbox{${}_{_\circ}$} }}
\def\frownie{\leavevmode\raise1pt\hbox{$\bigcirc$\kern-8pt\raise1.5pt\hbox{$\cdot\kern0.5pt\cdot$}\kern-6.35pt\lower3.5pt\hbox{${}^{{}_\frown}$}\kern-6.2pt\raise3.7pt\hbox{${}_{_\circ}$} }}

\def\chem$#1${\mbox{$\rm\newdimen\fdfourteen\newdimen\fdsixteen\newdimen\fdseventeen\fdfourteen=\fontdimen14\textfont2\fdsixteen=\fontdimen16\textfont2
\fdseventeen=\fontdimen17\textfont2\fontdimen14\textfont2=2.7pt\fontdimen16\textfont2=2.7pt\fontdimen17\textfont2=2.7pt
#1$\fontdimen14\textfont2=\fdfourteen\fontdimen16\textfont2=\fdsixteen\fontdimen17\textfont2=\fdseventeen}}

\def\frcn#1#2{\mbox{$\frac{#1}{#2}$}}
\def\Proof{\noindent{\sc Proof}}
\def\QED{\hfill\mbox{$\square$}}
\def\D{\b\nabla}
\def\frc{\displaystyle\frac}
\def\Int{\oint_{S_\infty}\kern-5pt}
\def\INT{\displaystyle\int}
\def\({\left(}
\def\){\right)}
\def\as#1->#2{\qquad\mbox{as $#1\to#2$}}
\def\Res{\mathop{\rm Res}\nolimits}
\def\d{{\boldmath\mbox{$d$}}}
\def\d{{\rm d}}
\def\m{{\boldmath\mbox{$m$}}}
\def\b#1{\ifcat a#1\mathbf{#1}\else\bbb#1\fi}
\def\bbb#1{\mbox{\boldmath\mbox{$#1$}}}
\def\ip<#1,#2>{\langle#1,#2\rangle}
\def\cosec{\mathop{\rm cosec}\nolimits}
\def\h{\mbox{$\frac12$}}
\def\tilde{\widetilde}
\textheight=8.3truein
\def\<{\langle}
\def\>{\rangle}
\def\up{\ket|\kern-3pt\uparrow>}
\def\down{\ket|\kern-3pt\downarrow>}
\def\Dt{{\frc \d{\d t}}}
\def\twid#1{\smash{\mbox{${\mathrel{\sim}}^{\kern-7pt\strut#1}$}\,}}
\def\sup{\mathop{\rm sup}\nolimits}
\def\inf{\mathop{\rm inf}\nolimits}
\def\arg{\mathop{\rm arg}\nolimits}
\def\erf{\mathop{\rm erf}\nolimits}
\def\Aut{\mathop{\rm Aut}\nolimits}
\def\erfc{\mathop{\rm erfc}\nolimits}
\def\diag{\mathop{\rm diag}\nolimits}
\def\sech{\mathop{\rm sech}\nolimits}
\def\tanh{\mathop{\rm tanh}\nolimits}
\def\cosech{\mathop{\rm cosech}\nolimits}
\def\coth{\mathop{\rm coth}\nolimits}
\def\ket|#1>{|#1\>}
\def\bra<#1|{\<#1|}
\def\IP<#1|#2>{\<#1|#2\>}
\newcount\questionnumber
\questionnumber=1
\def\Q#1.{\bigskip\noindent{\bf #1. }}
\def\s#1,#2.{\left#1,#2\right}

\def\sb#1{\b{\scriptstyle#1}}
\def\sba{\b {\scriptstyle a}}
\def\sbb{\b {\scriptstyle b}}
\def\sbc{\b {\scriptstyle c}}
\def\sbd{\b {\scriptstyle d}}
\def\sbe{\b {\scriptstyle e}}
\def\sbf{\b {\scriptstyle f}}
\def\sbg{\b {\scriptstyle g}}
\def\sbh{\b {\scriptstyle h}}
\def\sbi{\b {\scriptstyle i}}
\def\sbj{\b {\scriptstyle j}}
\def\sbk{\b {\scriptstyle k}}
\def\sbl{\b {\scriptstyle l}}
\def\sbm{\b {\scriptstyle m}}
\def\sbn{\b {\scriptstyle n}}
\def\sbo{\b {\scriptstyle o}}
\def\sbp{\b {\scriptstyle p}}
\def\sbq{\b {\scriptstyle q}}
\def\sbr{\b {\scriptstyle r}}
\def\sbs{\b {\scriptstyle s}}
\def\sbt{\b {\scriptstyle t}}
\def\sbu{\b {\scriptstyle u}}
\def\sbv{\b {\scriptstyle v}}
\def\sbw{\b {\scriptstyle w}}
\def\sbx{\b {\scriptstyle x}}
\def\sby{\b {\scriptstyle y}}
\def\sbz{\b {\scriptstyle z}}
\def\sbA{\b {\scriptstyle A}}
\def\sbB{\b {\scriptstyle B}}
\def\sbC{\b {\scriptstyle C}}
\def\sbD{\b {\scriptstyle D}}
\def\sbE{\b {\scriptstyle E}}
\def\sbF{\b {\scriptstyle F}}
\def\sbG{\b {\scriptstyle G}}
\def\sbH{\b {\scriptstyle H}}
\def\sbI{\b {\scriptstyle I}}
\def\sbJ{\b {\scriptstyle J}}
\def\sbK{\b {\scriptstyle K}}
\def\sbL{\b {\scriptstyle L}}
\def\sbM{\b {\scriptstyle M}}
\def\sbN{\b {\scriptstyle N}}
\def\sbO{\b {\scriptstyle O}}
\def\sbP{\b {\scriptstyle P}}
\def\sbQ{\b {\scriptstyle Q}}
\def\sbR{\b {\scriptstyle R}}
\def\sbS{\b {\scriptstyle S}}
\def\sbT{\b {\scriptstyle T}}
\def\sbU{\b {\scriptstyle U}}
\def\sbV{\b {\scriptstyle V}}
\def\sbW{\b {\scriptstyle W}}
\def\sbX{\b {\scriptstyle X}}
\def\sbY{\b {\scriptstyle Y}}
\def\sbZ{\b {\scriptstyle Z}}
\def\ba{\b a}
\def\bb{\b b}
\def\bc{\b c}
\def\bd{\b d}
\def\bg{\b g}
\def\bh{\b h}
\def\bi{\b i}
\def\bj{\b j}
\def\bk{\b k}
\def\bl{\b l}
\def\bm{\b m}
\def\bn{\b n}
\def\bo{\b o}
\def\bp{\b p}
\def\bq{\b q}
\def\br{\b r}
\def\bs{\b s}
\def\bt{\b t}
\def\bu{\b u}
\def\bv{\b v}
\def\bw{\b w}
\def\bx{\b x}
\def\by{\b y}
\def\bz{\b z}
\def\bA{\b A}
\def\bB{\b B}
\def\bC{\b C}
\def\bD{\b D}
\def\bE{\b E}
\def\bF{\b F}
\def\bG{\b G}
\def\bH{\b H}
\def\bI{\b I}
\def\bJ{\b J}
\def\bK{\b K}
\def\bL{\b L}
\def\bM{\b M}
\def\bN{\b N}
\def\bO{\b O}
\def\bP{\b P}
\def\bQ{\b Q}
\def\bR{\b R}
\def\bS{\b S}
\def\bT{\b T}
\def\bU{\b U}
\def\bV{\b V}
\def\bW{\b W}
\def\bX{\b X}
\def\bY{\b Y}
\def\bZ{\b Z}
\def\BA{\mathbb{A}}
\def\BB{\mathbb{B}}
\def\BC{\mathbb{C}}
\def\BD{\mathbb{D}}
\def\BE{\mathbb{E}}
\def\BF{\mathbb{F}}
\def\BG{\mathbb{G}}
\def\BH{\mathbb{H}}
\def\BI{\mathbb{I}}
\def\BJ{\mathbb{J}}
\def\BK{\mathbb{K}}
\def\BL{\mathbb{L}}
\def\BM{\mathbb{M}}
\def\BN{\mathbb{N}}
\def\BO{\mathbb{O}}
\def\BP{\mathbb{P}}
\def\BQ{\mathbb{Q}}
\def\BR{\mathbb{R}}
\def\BS{\mathbb{S}}
\def\BT{\mathbb{T}}
\def\BU{\mathbb{U}}
\def\BV{\mathbb{V}}
\def\BW{\mathbb{W}}
\def\BX{\mathbb{X}}
\def\BY{\mathbb{Y}}
\def\BZ{\mathbb{Z}}
 \def\CD{{\cal D}}
 \def\Re{\mathop{\rm Re}\nolimits}
\def\mtx#1,#2,#3,#4.{\pmatrix{#1&#2\cr#3&#4}}
\def\ge{\geqslant}
\def\le{\leqslant}
\def\hat{\widehat}
\def\emptyset{\varnothing}
\def\bar{\overline}
\def\mob#1,#2,#3,#4[#5]{\frc{#1#5+#2}{#3#5+#4}}
\def\Im{\mathop{\rm Im}\nolimits}
\def\Ker{\mathop{\rm Ker}\nolimits}
\def\n{\mathop{\rm n}\nolimits}
\def\r{\mathop{\rm r}\nolimits}
\def\p#1#2.{\frc{\partial#2}{\partial x_{#1}}}
\def\curl{\D\times}
\def\div{\D\cdot}
\def\dS{\cdot\bd\bS}
\def\pp#1,#2{\partial#1/\partial#2}
\def\ppp#1,#2{\frc{\partial#1}{\partial x_{#2}}}
\def\pd#1,#2{\frc{\partial#1}{\partial#2}}
\def\pdc#1,#2,#3{\left.\pd#1,#2\right|_{#3}\kern-6pt}
\def\cross{\times}
\def\kmh{$\,\rm kmh^{-1}$}
\def\sign{\mathop{\rm sign}\nolimits}
\def\imp{\Rightarrow}
\def\impby{\Leftarrow}
\def\nimp{\not\Rightarrow}
\def\nimpby{\kern3pt\not\kern-3pt\Leftarrow}
\def\limsup{\mathop{\rm limsup}\limits}
\def\liminf{\mathop{\rm liminf}\limits}
\def\cupp{\mathop{\displaystyle\cup}\limits}
\def\capp{\mathop{\displaystyle\cap}\limits}
\def\openone{\leavevmode\hbox{\rm\small1\kern-3.0pt\normalsize1}}
\def\max{\mathop{\rm max} \limits}
\def\const{\mbox{constant}}
\def\Tr{\mathop{\rm Tr}\nolimits}
\def\esssup{\mathop{\rm ess\,sup}\limits}
\def\cl#1,#2.{\left[\matrix{#1\cr#2}\right]}
\def\det{\mathop{\rm det}\nolimits}
\def\sup{\mathop{\rm sup}\limits}
\def\inf{\mathop{\rm inf}\limits}
\def\arg{\mathop{\rm arg}\nolimits}
\newcount\questionnumber
\questionnumber=1
\def\Q#1.{\bigskip\noindent{\bf #1. }}
\def\s#1,#2.{\left#1,#2\right}
\def\CA{{\cal A}}
\def\CB{{\cal B}}
\def\CC{{\cal C}}
\def\CD{{\cal D}}
\def\CE{{\cal E}}
\def\CF{{\cal F}}
\def\CG{{\cal G}}
\def\CH{{\cal H}}
\def\CI{{\cal I}}
\def\CJ{{\cal J}}
\def\CK{{\cal K}}
\def\CL{{\cal L}}
\def\CM{{\cal M}}
\def\CN{{\cal N}}
\def\CO{{\cal O}}
\def\CP{{\cal P}}
\def\CQ{{\cal Q}}
\def\CR{{\cal R}}
\def\CS{{\cal S}}
\def\CT{{\cal T}}
\def\CU{{\cal U}}
\def\CV{{\cal V}}
\def\CW{{\cal W}}
\def\CX{{\cal X}}
\def\CY{{\cal Y}}
\def\CZ{{\cal Z}}
\def\FA{\mathfrak {A}}
\def\FB{\mathfrak {B}}
\def\FC{\mathfrak {C}}
\def\FD{\mathfrak {D}}
\def\FE{\mathfrak {E}}
\def\FF{\mathfrak {F}}
\def\FG{\mathfrak {G}}
\def\FH{\mathfrak {H}}
\def\FI{\mathfrak {I}}
\def\FJ{\mathfrak {J}}
\def\FK{\mathfrak {K}}
\def\FL{\mathfrak {L}}
\def\FM{\mathfrak {M}}
\def\FN{\mathfrak {N}}
\def\FO{\mathfrak {O}}
\def\FP{\mathfrak {P}}
\def\FQ{\mathfrak {Q}}
\def\FR{\mathfrak {R}}
\def\FS{\mathfrak {S}}
\def\FT{\mathfrak {T}}
\def\FU{\mathfrak {U}}
\def\FV{\mathfrak {V}}
\def\FW{\mathfrak {W}}
\def\FX{\mathfrak {X}}
\def\FY{\mathfrak {Y}}
\def\FZ{\mathfrak {Z}}
\def\Fa{\mathfrak {a}}
\def\Fb{\mathfrak {b}}
\def\Fc{\mathfrak {c}}
\def\Fd{\mathfrak {d}}
\def\Fe{\mathfrak {e}}
\def\Ff{\mathfrak {f}}
\def\Fg{\mathfrak {g}}
\def\Fh{\mathfrak {h}}
\def\Fi{\mathfrak {i}}
\def\Fj{\mathfrak {j}}
\def\Fk{\mathfrak {k}}
\def\Fl{\mathfrak {l}}
\def\Fm{\mathfrak {m}}
\def\Fn{\mathfrak {n}}
\def\Fo{\mathfrak {o}}
\def\Fp{\mathfrak {p}}
\def\Fq{\mathfrak {q}}
\def\Fr{\mathfrak {r}}
\def\Fs{\mathfrak {s}}
\def\Ft{\mathfrak {t}}
\def\Fu{\mathfrak {u}}
\def\Fv{\mathfrak {v}}
\def\Fw{\mathfrak {w}}
\def\Fx{\mathfrak {x}}
\def\Fy{\mathfrak {y}}
\def\Fz{\mathfrak {z}}
\def\R{{I\!\! R}}
\def\N{{I\!\! N}}
\def\H{{I\!\! H}}
\def\Z{{Z\!\!\! Z}}\def\ope{\; \lra_{\!\!\!\!\!\!\!\!\!\!\!\!\!_{\hbox{$_{n\to\infty}$}}} \;}\def\C{\hskip3.5pt\hbox{\vrule width 0.6pt height 7pt depth 0pt \hskip -3.5pt}C}
\def\P{{I\!\! P}}
\def\E{{I\!\! E}}
\def\A{{I\!\! A}}
\def\B{{I\!\! B}}
\def\DD{{I\!\! D}}
\def\F{{I\!\! F}}
\def\H{{I\!\! H}}
\def\I{{I\!\! I}}
\def\G{\hskip3.5pt\hbox{\vrule width 0.6pt height 7pt depth 0pt \hskip -3.5pt}G}\def\P{{I\!\! P}}
\def\J{{I\!\! J}}
\def\K{{I\!\! K}}
\def\L{{I\!\! L}}
\def\M{{I\!\! M}}
\def\S{\hskip4.5pt\hbox{\vrule width 0.6pt height 4.6pt depth 0pt \hskip -4.5pt}\hskip5pt\raise4.6pt\hbox{\vrule width 0.6pt height 3.6pt depth0pt\hskip-5pt}S}
\def\T{{T\!\! T}}
\def\U{{I\!\! U}}
\def\V{{V\!\!\! V}}
\def\W{{W\!\! W}}
\def\X{{X\!\! X}}
\def\Y{{Y\!\! Y}}
\def\O{\hskip3.5pt\hbox{\vrule width 0.6pt height 7pt depth 0pt \hskip -3.5pt}O}\def\Q{\hskip3.5pt\hbox{\vrule width 0.6pt height 7pt depth 0pt \hskip -3.5pt}Q}
\def\weakto{\stackrel{\rm  weak}\longrightarrow}
\def\dto{\stackrel d\to}
\def\bplim{\mathop{\mbox{\rm bp-lim}}\limits}
\def\openone{\mbox{\ \kern-4pt\kern2pt\vrule width 0.4pt height 5.9 pt depth 0 pt\kern-1.75pt\vrule width 4.3 pt height 0.2 pt depth 0 pt\kern-1.75pt\vrule width 0.4 pt height 7 pt depth 0pt\kern-4.2pt{\raise0.1pt\hbox{\normalsize\sf\'{}}}\kern-4.4pt{\raise2.4pt\hbox{\tiny\sf\'{}}}\kern-1.3pt\vrule width 1.0pt depth0pt height0.3pt\kern-0.6pt\vrule width 0.3pt depth0pt height0.4pt\kern1.4pt\vrule width 0.3pt depth0pt height0.4pt\vrule width 0.5pt depth0pt height0.3pt\kern1pt {}}}

\def\openone{\,1\kern-4.4pt\vrule width0.3pt depth0pt height 6pt\kern-1.8pt\kern0.4pt\raise5.7pt\hbox{\vrule height 0.3pt depth 0pt width 1.6pt}\kern-1.99pt\vrule height0.3pt depth0pt width1.8pt\kern-0.6pt\vrule height 0.4pt depth 0pt width0.3pt\kern0.0pt\vrule depth0pt height 0.3pt width 1pt\kern4pt{}}

\def\supp{\mathop{\rm supp}\nolimits}

\normalsize

\long\def\gobble#1\endgobble{
}

\def\Pd#1,#2,#3.{\(\pd#1,{#2}\)_{\!\!\!#3}}
\def\e{{\rm e}}
\def\D{{\nabla}}
\def\DD{\BD}
\def\div{\mathop{\rm div }}
\def\grad{\mathop{\rm grad }}
\def\grad{\D}
\def\curl{\mathop{\rm curl }}
\def\delsq{\triangle}
\def\:{|\kern-0.8pt\|}
\def\and{\mbox{\qquad and \qquad}}

\def\bplim{\mathop{\mbox{\rm bp-lim}}\limits}
\def\weakto{\rightharpoonup}
\def\Eq.#1.{equation~{\rm(\ref{eq:#1})}}
\def\Ineq.#1.{inequality~{\rm(\ref{eq:#1})}}
\def\FR{\FC}
\def\BS{\FM}
\def\ell{l}
\def\co{\mathop{\rm co}\nolimits}
\def\fext{f^\#}
\def\pim{\pi_{\rm m}}
\def\cN{c_N}
\def\cNn{c_{N_n}}
\def\aNn{m_{N_n}}
\def\aN{m_N}
\def\Id1{I\kern-3.5pt I_{d_1}}





\title{A stochastic approximation scheme and convergence theorem for particle interactions with perfectly reflecting boundaries}
\author{\sc \small C. G. Wells\\\ \\
 { \rm\small Department of Chemical
    Engineering \\ University of Cambridge \\ Pembroke Street \\  Cambridge CB2 3RA  \\
UK \\ E-Mail: cgw11@cam.ac.uk\\
}}







\def\bmu{{\mathchoice{{\b\mu}}{{\b\mu}}{{\mbox{\boldmath$\scriptstyle\mu$}}}{{\mbox{\boldmath$\scriptscriptstyle\mu$}}}}}
\parskip5pt
 
\thispagestyle{empty}
\markboth{\hfil C. G. Wells\hfil\hfil}{\hfil C. G. Wells}
\pagestyle{myheadings}

\gobble
\thispagestyle{empty}
\markboth{\hfil C. G. Wells\hfil\hfil}{\hfil C. G. Wells}
\pagestyle{myheadings}


\centreline{\huge A stochastic approximation scheme and convergence}
\centreline{\huge theorem for particle interactions with}
\centreline{\huge perfectly reflecting boundary conditions}

\vskip1cm
\centreline{\sc Clive G. Wells\def\thefootnote{}\footnote{{\ni\em Email address:\/} {\tt cgw11{\rm@}cam.ac.uk}, Tel: +44 (0)1223 762785, Fax: +44 (0)1223 334796}\footnote{Research supported by the Oppenheimer Foundation.}}
\vskip0.3cm

\centreline{{\it\small Department of Chemical Engineering, University of Cambridge, New Museum Site,}}
\centreline{{\it\small Pembroke Street, Cambridge CB2 3RA, United Kingdom.}}
 
\vskip1cm

\begin{abstract}
We prove the existence of a solution to an equation governing the number density within a compact domain of a discrete particle system for a prescribed class of particle interactions taking into account the effects of the diffusion and drift of the set of particles. Each particle carries a number of internal coordinates which may evolve continuously in time, determined by what we will refer to as the {\em internal drift\/}, or discretely via the interaction kernels. Perfectly reflecting boundary conditions  are imposed on the system and all the processes may be spatially and temporally inhomogeneous. We use a relative compactness argument to construct a sequence of measures that converge weakly to a solution of the governing equation. 
Since the proof of existence is a constructive one, it provides a stochastic approximation scheme that can be used for the numerical study of molecular dynamics. 
\vskip5mm

\end{abstract}
\vskip1cm
\vfill

\break

\endgobble


\centreline{\huge A stochastic approximation scheme and convergence}
\centreline{\huge theorem for particle interactions with}
\centreline{\huge perfectly reflecting boundary conditions}

\vskip1cm
\centreline{\sc Clive G. Wells\def\thefootnote{}\footnote{{\ni\em Email address:\/} {\tt cgw11{\rm@}cam.ac.uk}, Tel: +44 (0)1223 762785, Fax: +44 (0)1223 334796}\footnote{Research supported by the Oppenheimer Foundation.}}
\vskip0.3cm

\centreline{{\it\small Department of Chemical Engineering, University of Cambridge, New Museum Site,}}
\centreline{{\it\small Pembroke Street, Cambridge CB2 3RA, United Kingdom.}}
 
\vskip1cm
\hrule

\begin{abstract}
We prove the existence of a solution to an equation governing the number density within a compact domain of a discrete particle system for a prescribed class of particle interactions taking into account the effects of the diffusion and drift of the set of particles. Each particle carries a number of internal coordinates which may evolve continuously in time, determined by what we will refer to as the {\em internal drift\/}, or discretely via the interaction kernels. Perfectly reflecting boundary conditions  are imposed on the system and all the processes may be spatially and temporally inhomogeneous. We use a relative compactness argument to construct a sequence of measures that converge weakly to a solution of the governing equation. 
Since the proof of existence is a constructive one, it provides a stochastic approximation scheme that can be used for the numerical study of molecular dynamics. 
\vskip5mm

\end{abstract}
\vskip1cm
\hrule

\setcounter{section}{0}
\ssection{Introduction}

Under fairly general conditions, set out in section~2.2, we prove the existence of at least one measure-valued solution to the following equation:
\beq
\int_\FR f(z)P_t(\d z)&=&\int_\FR f(z)P_0(\d z)+\int_0^t\int_\FR\h a_s(z)\triangle f(z)\, P_s(\d z)\,\d s+\int_0^t\int_\FR b_s(z)\cdot\D f(z)\, P_s(\d z)\,\d s\cr
&&{}+\int_0^t\int_\FR H_s(z)\cdot\hat\D f(z)\, P_s(\d z)\,\d s+\int_0^t \int_\FR f(z)I_s(\d z)\,\d s\cr
&&\kern1cm{}+\sum_{i=1}^\CI\sum_{j=1}^\CJ\frc1{i!\,j!}\int_0^t\int_{\FR^{i+j}}K_{s,i,j}(z_1,\ldots,z_i,\d w_1,\ldots,\d w_j)\cr
&&\kern3cm\times\left\{\sum_{m=1}^jf(w_m)-\sum_{n=1}^if(z_n)\right\}\, P_s(\d z_1)\cdots P_s(\d z_i)
\,\d s\LABEL{eq:model}
\eq
for $t\ge0$ and $f$ a member of a suitable class of test functions obeying the boundary condition $\pp f,n=0$, where $n$ is the outward pointing unit normal. This condition corresponds to the imposition of perfectly reflecting boundaries. In \Eq.model., $\FR$ is the configuration space consisting of a physical domain $\Omega$ and a collection of internal degrees of freedom and $P_t$ is the measure representing the number density of particles in $\FR$. The functions $a_t$, $b_t$ and $H_t$ are rates associated with diffusion and the external and internal drifts respectively. The measure $I_t$ describes a source term for the system and $K_{t,i,j}$ are interaction kernels describing the interplay of particles within the domain.

This equation has many potential applications for the quantitative description of molecular systems. It generalizes Smoluchowski's coagulation equation~\cite{vS} and simple diffusion models, e.g.,~\cite{SV}. Existence theorems for coagulation and fragmentation processes have been developed for spatially and temporally inhomogeneous conditions in, e.g., \cite{Spouge,Stewart,BallCarr,Jeon,James,EW} (see also the references therein). The case of `cluster coagulation', 
extending this work, where the coagulation rate may have non-trivial dependence on some internal coordinates has been examined in \cite{James-cluster}. For a survey of some of the known exact solutions and a discussion of the mathematical and physical reasons for studying the coagulation equation see the review article~\cite{Aldous}.

Diffusion processes with boundary conditions have been studied using SDE techniques by, e.g., \cite{SV71,Skorohod,Watanabe}. 
However we do not follow this route since here it is our aim to produce a {\em constructive proof\/} for the existence of solutions to \Eq.model.. This is in order to provide a Monte Carlo approximation scheme capable of simple computer implementation that can simulate the underlying physical system accurately.

The existence problem for combined diffusion, coagulation and fragmentation processes has received some attention in~\cite{vanDongen,Collet,LM01,Deaconu} for the case of spatially and temporally homogeneous diffusion coefficients. Van Dongen discusses the case when the diffusion coefficients are non-zero and equal for masses no larger than some given mass and zero for larger masses. As a separate case he considers the situation when all diffusion coefficients are equal. 
Collet \& Poupaud found a local existence result in the former case for a restricted range of coagulation kernels and a global existence result in the latter case. Guia\c s~\cite{Guiasthesis,Gui01}, developing the work in~\cite{vanDongen}, has presented a numerical approach to this problem by dividing the domain into cells then solving  each cell as a homogeneous coagulation system with appropriate transfer of particles between neighbouring cells; thus combining a finite element and Monte Carlo approach. Guia\c s provides the necessary mathematical theorems to justify this strategy. In contrast to this combined finite element and Monte Carlo method, our resulting numerical approach is purely stochastic in nature, and the diffusion, drift and other rate coefficients may have spatial and temporal dependence, we also allow for more general forms of particle interaction to take place.

The paper is organized as follows.
In section~2 we give definitions of the physical domain and configuration space, and describe the hypotheses we will be making. Then we introduce the class of test functions and describe the boundary condition they obey. Before beginning the mathematical assault on the problem, we relate the boundary condition on the class of test functions to appropriate Neumann boundary conditions on the particle number density when such a density exists with respect to Lebesgue measure and certain extra differentiability conditions are assumed. Additionally in section~3 we characterize the boundary conditions associated with the internal degrees of freedom.

Section~4 is devoted primarily to the proofs of proposition~\ref{prop:boundedness} and theorem~\ref{thm:uniform-internal} that are needed in section~5. Proposition~\ref{prop:boundedness} provides a bound for an appropriate stochastic generator associated with the diffusion and drift terms in \Eq.model., whilst theorem~\ref{thm:uniform-internal} provides a result that allows us to deduce the uniform convergence of a discrete approximation to the derivative terms in the integrand of the RHS of \Eq.model.. To implement the diffusion process each particle's position is displaced by a multivariate normal random variable that is scaled proportional to the square root of the time step. The resulting displacement is modified to account for the drift terms and for any reflections off the boundaries. Since in a finite time the distance travelled by such a particle becomes unbounded as the time steps decrease, the proof relies on the cancellation of motions in opposite directions and it is this that requires delicate analysis, particularly in a neighbourhood of the boundary. We point out that the use of such a multivariate normal random variable makes direct contact with the mathematical definition of Brownian motion and with the paper of Roth~\cite{Roth} on generating a Feller semi-group associated with a simple diffusive process.

In section~5 we formulate a relative compactness argument for a sequence of stochastic processes taking values in an appropriate product space. The components of each term in the  sequence comprise a measure representing an approximation to the particle number density at each time and a non-negative real number approximating that time. This {\em fictitious time\/} is expedient for recasting \Eq.model. as a pair of equations with no explicit time-dependence. Tightness of the sequence is proved and Prohorov's theorem is used to deduce the existence of a convergent subsequence. We verify that the limit is a solution to \Eq.model. in the weak sense when the Linear Growth Condition holds and in the vague sense if the growth condition on the internal drift and self-interaction kernel are relaxed. Finally a new sequence of processes is presented with the same jump chain of measures but with modified holding times determined by the the fictitious time variable. This has the effect of removing the fictitious time quantity from the construction. 
This final result, theorem~\ref{thm:main}, provides a simple approximation scheme
that is in essence a simple extension of that of Marcus~\cite{Marcus} and Lushnikov~\cite{Lushnikov} except that the holding times are determined deterministically from the current state. We remark that proposition~\ref{prop:M_t} and theorem~\ref{thm:main} can be modified simply to remove all spatial dependence and in this way we recover the earlier results of \cite{James,EW}.

In section~6 we provide an explicit description of the approximation scheme for generating sample paths for the stochastic processes, while in section~7 we make some remarks about potential applications to physical situations. These include a brief discussion of the implications and limitations of the various hypotheses assumed in the proof. Appendices A, B and C contain important supporting results that do not fit naturally within the main text.

\ssection{Definitions and Hypotheses}

For any pair of Banach spaces $X$ and $Y$, let $C_b(X;Y)$ denote the bounded continuous  functions from $X$ to $Y$ and write $C_b(X)=C_b(X;\BR)$. Let $C^n(X;Y)$ denote the $n$-times continuously differentiable $Y$-valued functions on $X$.
We denote by $C_c(X)$ the compactly supported real-valued continuous functions on $X$ and define $C^n_c(X)=C^n(X)\cap C_c(X)$.
Additionally  $\CM_B(X)$ denotes the space of bounded non-negative Borel measures on $X$, $\FB(X)$ is the Borel $\sigma$-algebra of $X$ and $M(X)$ is the space of bounded Borel measurable functions on $X$. For any $P\in\CM_B(X)$  and  $f\in M(X)$  we will write $P(f)=\int_Xf(x)\,P(\d x)$.

\ssubsection{The Configuration Space $\FR$}

\begin{definition}\LABEL{def:Omega}
Let $\Omega=\omega^{-1}(-\infty,0)$ be a bounded open subset of the Euclidean inner product space $(\BR^{d_1},\,\cdot\,)$, with the induced norm denoted by $\|\cdot\|$, and where
$\omega\in C^2(\BR^{d_1})$
with $\|\D\omega(x)\|>1$ for all $x\in\partial\Omega=\omega^{-1}\{0\}$. Let $\bar\Omega$ be the closure of $\Omega$ and 
 write $\D=(\partial_1,\ldots,\partial_{d_1})$ for the gradient operator on $\BR^{d_1}$ and $\triangle$ for the Laplacian.
\end{definition}

\begin{definition}
Define the expressions
\be
\|\D\omega\|_\infty=\sup_{x\in\co(\bar\Omega)}\(\sum_{i=1}^{d_1}|\partial_i\omega(x)|^2\)^{\frac12}\and\|\D\D\omega\|_\infty=\sup_{x\in\co(\bar\Omega)}\(\sum_{i,j=1}^{d_1}|\partial_i\partial_j\omega(x)|^2\)^{\frac12}
\ee
where $\co(A)$ denotes the convex hull of the set $A$.
\end{definition}

\begin{definition}\LABEL{def:internal}
Let $a^{(m)}_\mu\le b^{(m)}_\mu$ for $m\in\BN$ and $\mu=1,\ldots,d_2$. Define the sequence of compact sets $\Gamma_m=\prod_{\mu=1}^{d_2}[a_\mu^{(m)},b_\mu^{(m)}]\subset\BR^{d_2}$.
We endow $\Gamma_m$ with the Euclidean inner product induced from $\BR^{d_2}$ and write a dot for this inner product and $\|\cdot\|$ for the Euclidean norm. In addition we let  $\hat\D=(\hat\partial_1,\ldots,\hat\partial_{d_2})$ denote the gradient operator on $\BR^{d_2}$.
\end{definition}

\begin{definition}\LABEL{def:Configuration}
Define the configuration space $\FR$, and its extension $\FR^\#$, both subsets of $\BN\times\BR^{d_1}\times\BR^{d_2}$ by
\be
\FR=\cupp_{m=1}^\infty(\{m\}\times\bar\Omega\times\Gamma_m)\and
\FR^\#=\cupp_{m=1}^\infty(\{m\}\times\co(\bar\Omega)\times\Gamma_m),
\ee
and endow each with the Euclidean metric induced from $\BR^{1+d_1+d_2}$. Define the three restrictions: $\pim:\FR\to\BN$, $\pi_{\bar\Omega}:\FR\to\bar\Omega$ and $\pi_\Gamma:\FR\to\BR^{d_2}$ by
\beq
\pim(z_0,z_1,\ldots,z_{d_1+d_2})&=&z_{0};\\
\pi_{\bar\Omega}(z_0,z_1,\ldots,z_{d_1+d_2})&=&(z_1,\ldots,z_{d_1});\\
\pi_{\Gamma}(z_0,z_1,\ldots,z_{d_1+d_2})&=&(z_{d_1+1},\ldots,z_{d_1+d_2}),
\eq
so that $\pi_\Gamma(z)\in\Gamma_{\pim(z)}$.
It is useful to define $e_0,e_1,\ldots, e_{d_1},\hat e_1,\ldots,\hat e_{d_2}$ to be the canonical basis for $\BR^{1+d_1+d_2}$ in the obvious way. For convenience we will refer to $\pim(z)$ as {\em the mass\/} of a particle with coordinates $z$.
\end{definition}

\ssubsection{Hypotheses}
In this subsection we present two choices of physical hypotheses that are sufficient to prove an appropriate convergence result.  

\begin{assumption}\LABEL{hyp:cont}
Suppose the source rate is given by the continuous map $I:[0,\infty)\to\CM_B(\FR)$, with $I:t\mapsto I(t)\equiv I_t$
satisfying the boundedness condition
\be
\Lambda^{(\CI\vee2)}\equiv
\sup\left\{I_t(\pim^{(\CI\vee2)}):0\le t<\infty\right\}<\infty
\ee
where we define $\pim^{(q)}(z)=(\pim(z))^q$ and where $\CI$ some the non-negative integer defined implicitly by \Eq.model..
\end{assumption}

\begin{assumption}\LABEL{hyp:diffusion}
Suppose the diffusion coefficient $a:[0,\infty)\times\FR\to\BR$ (with $a(u,z)\equiv a_u(z)\equiv\sigma_u(z)^2$) is continuous and there exists a bound $\sigma_\infty$ such that
\be
0<\sigma_u(z)\le\sigma_\infty
\ee
for all $u\in[0,\infty)$ and $z\in\FR$. 
Write
\be
\sigma^{\rm min}_{T,m}=1\wedge\inf\{\sigma_t(z):z\in\FR,\ \pim(z)\le m,\ 0\le t\le T\},
\ee
which is strictly greater than zero by continuity on the compact set $[0,T]\times\pim^{-1}\{1,\ldots,m\}$.

Furthermore assume that
the drift function $b:[0,\infty)\times\FR\to\BR^{d_1}$ (with $b(u,z)\equiv b_u(z)$) is continuous and obeys the boundedness condition $\|b_u(z)\|\le b_\infty$ for some $b_\infty\in\BR$
and additionally satisfies the boundary condition
\be
b_u(z)\cdot n(\pi_{\bar\Omega}(z))=0\qquad\mbox{where}\qquad n(x)=\frc{\D\omega(x)}{\|\D\omega(x)\|}
\ee
for all $u\in[0,\infty)$ and $z\in\FR$ such that $\pi_{\bar\Omega}(z)\in\partial\Omega$.
\end{assumption}

\begin{assumption}\LABEL{hyp:internal}
Suppose the internal drift rate $H:(t,z)\mapsto H_t(z)\equiv(H_{t,1}(z),\ldots,H_{t,d_2}(z))
$ is a continuous function from
$[0,\infty)\times\FR$ to $\BR^{d_2}$ where
\be
H^\infty_m\equiv1\vee\sup\{\|H_t(z)\|/\pim(z):z\in\FR,\ \pim(z)\le m,\ t\in[0,\infty)\}<\infty
\ee
for all $m\in\BN$. As a matter of notational convenience we set $H^\infty_0=1$.
We will also impose the boundary conditions: $H_{t,\mu}(z)\ge0$ whenever $(\pi_\Gamma(z))_\mu=a^{(\pim(z))}_\mu$ and $H_{t,\mu}(z)\le0$ whenever $(\pi_\Gamma(z))_\mu=b^{(\pim(z))}_\mu$.
\end{assumption}
The purpose of these boundary conditions is to ensure that the internal drift does not take particles outside of the domain of the internal coordinates.

\begin{assumption}\LABEL{hyp:interactionker}
Suppose there exist 
uniformly
continuous interaction kernels $K_{i,j}$ for $i=1,\ldots,\CI$, $j=1,\ldots,\CJ$ and $(i,j)\neq(1,1)$ with $K_{i,j}:[0,\infty)\times\FR^i\to\CM_B(\FR^j)\cong\bigotimes_{k=1}^j\CM_B(\FR)$ where we write
\be
K_{u,i,j}(z_1,\ldots,z_i,B_1\times\cdots\times B_j)=K_{i,j}(u,z_1,\ldots,z_i)(B_1\times\cdots\times B_j)
\ee
which we require to be totally symmetric in $B_1,\ldots,B_j\in\FB(\FR)$, totally symmetric in  $z_1,\ldots,z_i$ and to satisfy
\be
\sum_{\ell=1}^i\pim(z_\ell)=\sum_{n=1}^j\pim(w_n),\qquad K_{u,i,j}(z_1,\ldots,z_i,\d w_1,\ldots,\d w_j)\mbox{--a.e.}
\ee
for all $z_1,\ldots,z_i\in\FR$, $i=1,\ldots,\CI$, and $j=1,\ldots,\CJ$, i.e., the interaction is mass preserving. Moreover we assume that there exists a  constant $K_\infty$ such that
\be
\frc1{i!\,j!}{K_{u,i,j}(z_1,\ldots,z_i,\FR^j)}\le K_\infty\pim(z_1)\cdots\pim(z_i)\qquad\mbox{for $(i,j)\neq(1,1)$}.
\LABEL{eq:K_infty}
\ee
In addition we assume that there exists a continuous (though not necessarily uniformly continuous) {\em self-interaction kernel\/}, $K_{1,1}:[0,\infty)\times\FR\to\CM_B(\FR)$ where we write $K_{u,1,1}(z,B)=K_{1,1}(u,z)(B)$ satisfying 
$\pim(z)=\pim(w)$, $K_{u,1,1}(z,\d w)$--a.e.{}
and possessing
 the bound:
\be
K^\infty_m\equiv1\vee\sup\{K_{t,1,1}(z,\FR)/\pim(z):z\in\FR,\ \pim(z)\le m,\ t\in[0,\infty)\}<\infty.
\ee
Analogous to our definition of $H^\infty_0$ we set $K^\infty_0=1$.
\end{assumption}

\begin{assumption}\LABEL{hyp:P_0}
We assume that the initial measure-valued particle density has finite mass:
 $P_0(\pim)<\infty$.
\end{assumption}

\begin{remark}\LABEL{rem:lineargrowth}\em
If both $K^\infty_m$ and $H^\infty_m$ are bounded functions on $\BN$ then we say the system obeys the {\em Linear Growth Condition\/}. In which case we may chose $K_\infty$ and $H_\infty$ such that for all $m\in\BN_0$ we have $K^\infty_m\le K_\infty$, $H^\infty_m\le H_\infty$ and \Eq.K_infty. still holds.
When this additional condition holds we will be able to prove that \Eq.model. holds in the weak sense. When it does not we show that \Eq.model. is true in the vague sense.
\end{remark}

\break
\ssubsection{The Class of Test Functions}
\begin{definition}\LABEL{def:testfunctions}
Define the normed space of functions
\beq
\FF&=&\{f:f(m,\cdot)\in C^2(\bar\Omega\times\Gamma_m), \mbox{\ for all $m\in\BN$, $\:f\:<\infty$}
\cr&&\kern2cm\mbox{ and
$ \D\omega(\pi_{\bar\Omega}(z))\cdot\D f(z)=0$ for all $z\in\FR$ such that $\pi_{\bar\Omega}(z)\in\partial\Omega$}\}\LABEL{eq:frakF}\qquad
\eq
where $\FF$ and its subspaces have the norm $\:\cdot\:$ defined by
\be
\:f\:=\max\{\|f\|_\infty,\|\D f\|_\infty,\|\hat\D f\|_\infty,\|\D\D f\|_\infty\LABEL{eq:f-one}
\}
\ee
with
\be
\|\D f\|_\infty=\sup_{x\in\FR}\(\sum_{i=1}^{d_1}|\partial_if(z)|^2\)^{\frac12};\qquad\|\D\D f\|_\infty=\sup_{x\in\FR}\(\sum_{i,j=1}^{d_1}|\partial_i\partial_jf(z)|^2\)^{\frac12}
\ee
and
\be
\|\hat\D f\|_\infty=\sup_{x\in\FR}\(\sum_{\mu=1}^{d_2}|\hat\partial_\mu f(z)|^2\)^{\frac12}.\LABEL{eq:f-two}
\ee
Furthermore we define $\FF_c=\FF\cap C_c(\FR)$ and let the map $R:\FF\to\BN_0$ be given by
\be
R(f)=\begin{cases}\max\{\pim(\supp f)\}&\mbox{if $f\in\FF_c$;}\\0&\mbox{otherwise.}\end{cases}
\ee
\end{definition}

\begin{remark}\em
By proposition~\ref{prop:separable} $(\FF_c,\:\cdot\:)$ is separable, therefore there exists a dense sequence $f_n\in\FF_c$. We set $f_0=1$.
Proposition~\ref{prop:Whitney} 
implies that for each $f\in\FF_c\cupp\{f_0\}$ there exists a second-order continuously differentiable extension $\fext:\FR^\#\to\BR$ with the property that $\:\fext\:^\#$ is finite, where 
 we define the mapping $\:\cdot\:^\#$ on each extension  by equations~(\ref{eq:f-one}--\ref{eq:f-two})  with $\FR$ replaced by $\FR^\#$.
\end{remark}

The reason for introducing these extensions is to allow us to assume that the test function is defined on a convex domain for each integral mass. This allows us to make elementary estimates using Taylor's theorem applied to $\fext$ without regard for the possibility that the straight line segment between the two points of interest may not lie entirely within the domain of the unextended test function $f$. In this way we need not assume any convexity property of the physical domain $\bar\Omega$.

\begin{definition}\LABEL{def:testsequence}
Let the sequence $(f_n)_{n=1}^\infty\in\FF_c$ be a dense subset of $(\FF_c,\:\cdot\:)$ with the following properties
\be
\:\fext_n\:^\#\le 2^n;\qquad
H^\infty_{R(f_n)} K^\infty_{R(f_n)}\le 2^{\frac12n}
\ee
and $(f_n)_{n=1}^\infty$ is closed under the operation $f\mapsto f(\cdot)\openone_{\pim(\cdot)< k}$ for all $k\in\BN$.
\end{definition}

\begin{remark}\em
It follows from definition~\ref{def:testsequence} that for every $f\in\FF$ there exists a subsequence of $f_n$ such that
\be
\bplim_{k\to\infty} f_{n_k}=f;\qquad
\bplim_{k\to\infty}\partial_if_{n_k}=\partial_i f;\qquad
\bplim_{k\to\infty}\partial_i\partial_jf_{n_k}=\partial_i\partial_j f
\ee
for $i,j=1,\ldots,d_1$ and
\be
\bplim_{k\to\infty}\hat\partial_\mu f_{n_k}=\hat\partial_\mu f
\ee
for $\mu=1,\ldots,d_2$ (where $\bplim$ denotes bounded pointwise convergence, i.e., $\bplim_{n\to\infty} f_n=f$ if and only if $f_n(z)\to f(z)$ as $n\to\infty$ for all $z\in\FR$ and $\sup_n\sup_{z\in\FR} |f_n(z)|<\infty$). Moreover if $f\in\FF_c$ then we may assume that
\be
R(f_{n_k})\le R(f)
\ee
for all $k$.
\end{remark}
  
\begin{definition}
We endow $\CM_B(\FR)$
and its subspaces with the complete metric generating the weak topology given by
\be
d_{\rm weak}(P,Q)=\sum_{n=0}^\infty 4^{-n}(|P(f_n)-Q(f_n)|\wedge1) \LABEL{eq:metric0}
\ee
and where $f_0=1$. We write $P_n\weakto P$ to signify weak convergence.
\end{definition}

\begin{remark}\em
We note that since $\FR$ is separable $\CM_B(\FR^k)\cong\bigotimes_{j=1}^k\CM_B(\FR)$ and if $P_n\weakto P$ then
$P_n{}^{\otimes k}\weakto P^{\otimes k}$.
\end{remark}
Finally, to establish our conventions, let us make note of the following definition and the subsequent remark; they will play a frequent r\^ole in proving the convergence of various quantities in the rest of the paper.

\begin{definition}
If $(X,d_X)$ and $(Y,d_Y)$ are metric spaces and $f:X\to Y$ is a continuous function we define the modulus of continuity by
\be
w(f,\delta)=\sup_{d_X(x,y)\le \delta}d_Y(f(x),f(y))\LABEL{eq:modulus}
\ee
whenever the RHS exists.
\end{definition}

\begin{remark}\em
The function $f$ is uniformly continuous if and only if \be
\lim_{\delta\to0}w(f,\delta)=0.\LABEL{eq:ctsmod}
\ee
In particular if $f:X\to\BR^m$
is continuous and compactly supported then $f$ is uniformly continuous and \Eq.ctsmod.  holds.
\end{remark}

\ssection{Neumann Boundary Conditions on the Number Density}

In this section we relate the properties of weak/vague solutions to \Eq.model. to the boundary conditions on the number density with respect to Lebesgue measure when 
such a quantity exists. Let us write $z=(m,\CX)$ and $\CX=(x,X)$ with $m\in\BN$, $x\in\bar\Omega$, $X\in\Gamma_m$
and assume that the following definitions can be made with respect to Lebesgue measure: $I_t(\{m\}\times\d\CX)=I^m_t(\CX)\,\d\CX$,
\be K_{t,i,j}(z_1,\ldots,z_i,\{n_1\}\times\d\CY_1,\ldots,\{n_j\}\times\d\CY_j)=K_{t,i,j}^{m_1,\ldots,m_i,n_1,\ldots,n_j}(\CX_1,\ldots,\CX_i,\CY_1,\ldots,\CY_j)\,\d\CY_1\ldots\d\CY_j,
\ee
where $K_{t,i,j}^{m_1,\ldots,m_i,n_1,\ldots,n_j}$ and $I^m_t$ are continuous functions on their respective domains, $P_t(\{m\}\times\d\CX)=c_t^m(\CX)\,\d\CX$, with $a_t\in C_b^2(\FR)$, $b_t\in C^1_b(\FR;\BR^{d_1})$, $H_t\in C_b^1(\FR;\BR^{d_2})$.
Suppose too that $c_t^m(\CX)$ is continuously differentiable in time and twice continuously differentiable on $\bar\Omega\times\Gamma_m$, then we discover for all $\CX\in\Omega\times\(\Gamma_m\)^\circ$ (the interior of $\bar\Omega\times\Gamma_m$) that
\beq
\kern0cm&&\kern-1cm\pd c^m_t(\CX),t=\h\triangle(a_t(z)c^m_t(\CX))-\D\cdot(b_t(z)c_t^m(\CX))-\hat\D\cdot\(H_t(z)c^m_t(\CX)\)+I^m_t(\CX)\kern3cm{}\ \cr
&&{}+\sum_{i=1}^\CI\sum_{j=1}^\CJ\frc{1}{i!\,(j-1)!}\kern-10pt\sum_{\begin{array}{c}{\scriptstyle m_1,\ldots, m_i=1}\\[-5pt]{\scriptstyle n_1,\ldots,n_{j-1}=1}\end{array}}^\infty\kern-10pt \int_{\bar\Omega\times\Gamma_{n_1}}\kern-20pt\d\CY_1\cdots\int_{\bar\Omega\times\Gamma_{n_{j-1}}}\kern-20pt\d\CY_{j-1}
\int_{\bar\Omega\times\Gamma_{m_1}}\kern-20pt c_t^{m_1}(\CX_1)\,\d\CX_1\cdots\int_{\bar\Omega\times\Gamma_{m_{i}}}\kern-20pt c_t^{m_i}(\CX_i)\,\d\CX_i\cr
&&\kern4cm{}\times K_{t,i,j}^{m_1,\ldots,m_i,n_1,\ldots,n_{j-1},m}(\CX_1,\ldots,\CX_i,\CY_1,\ldots,\CY_{j-1},\CX)\cr
&&{}-\sum_{i=1}^\CI\sum_{j=1}^\CJ\frc{c^m_t(\CX)}{(i-1)!\,j!}
\kern-5pt\sum_{\begin{array}{c}{\scriptstyle n_1,\ldots,n_j=1} \\[-5pt]{\scriptstyle m_1,\ldots,m_{i-1}=1}\end{array}}^\infty\kern-5pt \int_{\bar\Omega\times\Gamma_{n_1}}\kern-20pt\d\CY_1\cdots\int_{\bar\Omega\times\Gamma_{n_j}}\kern-10pt\d\CY_j
\int_{\bar\Omega\times\Gamma_{m_1}}\kern-20pt c_t^{m_1}(\CX_1)\,\d\CX_1\cdots\int_{\bar\Omega\times\Gamma_{m_{i-1}}}\kern-20ptc_t^{m_{i-1}}(\CX_{i-1})\,\d\CX_{i-1}\cr
&&\kern4cm{}\times K_{t,i,j}^{m_1,\ldots,m_{i-1},m,n_1,\ldots,n_j}(\CX_1,\ldots,\CX_{i-1},\CX,\CY_1,\ldots,\CY_j).\LABEL{eq:c_t^m}
\eq
This follows by considering a suitable sequence of test functions concentrated at $z$ and compactly supported in $\Omega\times\(\Gamma_m\)^\circ$. Now multiply \Eq.c_t^m. by $f(z)\in\FF_c$, subtract the time derivative of \Eq.model. and apply the divergence theorem to leave the boundary terms:
\beq
0&=&\int_{\Gamma_m}\int_{\partial\Omega}[\h f(z)\D(a_t(z)c^m_t(\CX))-\h a_t(z)c^m_t(\CX)\D f(z)-f(z)c^m_t(\CX)b_t(z)]\cdot n(x)\d S\,\d X\cr
&&-\sum_{\mu=1}^{d_2}\int_{a^{(m)}_1}^{b^{(m)}_1}\kern-7pt\d X_1\cdots
\int_{a^{(m)}_{\mu-1}}^{b^{(m)}_{\mu-1}}\kern-7pt\d X_{\mu-1}
\int_{a^{(m)}_{\mu+1}}^{b^{(m)}_{\mu+1}}\kern-7pt\d X_{\mu+1}\cdots
\int_{a^{(m)}_{d_2}}^{b^{(m)}_{d_2}}\kern-7pt\d X_{d_2}\int_{\bar\Omega}\left[H_{t,\mu}(z)c^m_t(\CX)f(z)\right]_{X_\mu=a^{(m)}_\mu}^{X_\mu=b^{(m)}_\mu}\,\d x\qquad
\eq
where $n(x)$ is the outward pointing unit normal at $x$. The last two terms in the first integrand on the RHS vanish by the definition of the boundary condition on the test functions and by hypothesis~\ref{hyp:diffusion} respectively.
By taking a suitable sequence of test functions concentrated at $z$ we deduce that $\pp(a_t(z)c^m_t(\CX)),n=0$ for $x\in\partial\Omega$ and $c^m_t(\CX)=0$ whenever $X_\mu=a_\mu^{(m)}$ or $b_\mu^{(m)}$ and $H_{t,\mu}(z)\neq0$.

\ssection{Perfectly Reflecting Boundaries}

The diffusion, internal and external drifts of each individual particle are to be approximated by a sequence of discrete jumps in the particle's position in the configuration space $\FR$. To complicate matters these jumps may involve reflections off the boundaries of the physical domain. In this section we provide some detailed calculations giving estimates and convergence results for our diffusion and drift approximations which will be vital in section~5.

\begin{proposition}\LABEL{prop:normal}
(a)
If  $x,y\in\partial\Omega$ with $\|y-x\|=\theta\neq0$, and $n(x)=\frc{\D\omega(x)}{\|\D\omega(x)\|}$ then there exists a constant $B_0$ such that
\be
|n(x)\cdot \hat k|\le B_0\theta
\qquad\mbox{where $\hat k=\frc{y-x}{\|y-x\|}$}.
\ee
(b) If $x,y\in\bar\Omega$ with $\|\D\omega(x)\|\ge1$ and $\D\omega(y)\neq0$ then
there exists a constant $A_0$ such that
\be
\|n(y)-n(x)\|\le A_0\theta.
\ee
\end{proposition}
\begin{proof}
(a) Use Taylor's theorem with $y=x+\theta\hat k$ and the facts that $\omega(x)=\omega(y)=0$ and $\|\D\omega(x)\|>1$ to deduce that
\beq
|n(x)\cdot\hat k|&\le&|\theta^{-1}\omega(x+\theta\hat k)-\theta^{-1}\omega(x)-\D\omega(x)\cdot\hat k|\le B_0\theta
\eq
where $B_0=\h\|\D\D\omega\|_\infty$.

\ni(b) Compute
\beq
\|n(y)-n(x)\|&=&\left\|
\frc{\D\omega(y)}{\|\D\omega(y)\|}-\frc{\D\omega(x)}{\|\D\omega(x)\|}\right\|\\
&=&\left\|\frc{\|\D\omega(x)\|-\|\D\omega(y)\|}{\|\D\omega(x)\|\,\|\D\omega(y)\|}\D\omega(y)+\frc{\D\omega(y)-\D\omega(x)}{\|\D\omega(x)\|}\right\|\\
&\le&\frc{\bigl|\,\|\D\omega(x)\|-\|\D\omega(y)\|\,\bigr|+\|\D\omega(y)-\D\omega(x)\|}{\|\D\omega(x)\|}\\
&\le&2\|\D\omega(y)-\D\omega(x)\|\le A_0\theta
\eq
where $A_0=2\|\D\D\omega\|_\infty$.
\end{proof}

\begin{definition}\LABEL{def:xi}
Let $\xi:\bar\Omega\times\BR^{d_1}\to\bar\Omega$ be defined by
\be
\xi(x,k)=x^*,\qquad\mbox{for some $x^*$ such that $\|x+k-x^*\|=\min\{\|x+k-y\|:y\in\bar\Omega\}$}.
\ee
This exists by the compactness of $\bar\Omega$.
\end{definition}

\begin{proposition}\LABEL{prop:xi-estimate}
If $x\in\partial\Omega$, $k\cdot n(x)=0$ and $\|k\|<\h B_0^{-1}$ then
\be
\|\xi(x,k)-x-k\|\le2B_0\|k\|^2.
\ee
\end{proposition}
\begin{proof}
Consider $y=x+K$ where $K=k-2B_0\|\D\omega(x)\|^{-1}\|k\|^2n(x)$, then
\beq
\omega(y)
&\le&
 (k-2B_0\|\D\omega(x)\|^{-1}\|k\|^2n(x))\cdot\D\omega(x)
 +\h\|\D\D\omega\|_\infty\|K\|^2\\
&\le&-B_0\|k\|^2(1-4B_0^2\|k\|^2)\le0,
\eq
and thus $y\in\bar\Omega$, it follows that
\be
\|\xi(x,k)-x-k\|\le\|y-x-k\|=2B_0\|\D\omega(x)\|^{-1}\|k\|^2\le2B_0\|k\|^2
\ee
as required.
\end{proof}

\begin{definition}
Define the map $\gamma:\bar\Omega\times\BR^{d_1}\to\bar\Omega$ by the following procedure. If $x\in\bar\Omega$ and $k\in\BR^{d_1}\setminus\{0\}$ then let $x_0=x$ and $\hat k_0=k/\|k\|$, let
$\vartheta_n=\inf\{t\in[0,\|k\|-\sum_{i=1}^{n-1}\vartheta_i]:  x_n+t\hat k_n\notin\Omega\}$, $x_{n+1}=x_n+\vartheta_n\hat k_n$, $\hat k_{n+1}=\hat k_n-2(\hat k_n\cdot n(x_{n+1}))n(x_{n+1})$.
We set $\gamma'(x,k)=\lim_{n\to\infty}x_n$ (this exists since the $\{x_n\}$ are contained in the compact ball of radius $\|k\|$ centred at $x$,
so there exists an accumulation point,
which is unique since $\sum_{n=1}^\infty\|x_n-x_{n-1}\|=\sum_{n=0}^\infty\vartheta_n$ is finite).
We define $\gamma'(x,0)=x$. Finally we set $\gamma(x,k)=\gamma'(x,k)$ if 
 $\sum_{n=0}^\infty\vartheta_n=\|k\|$  otherwise let
$\gamma(x,k)=\xi(\gamma'(x,k),(\|k\|-\sum_{r=0}^\infty\vartheta_r)\hat k_\infty)$, where we remark that proposition~\ref{prop:normal} implies that for $n>m$, $\|\hat k_n-\hat k_m\|\le2\sum_{i=m+1}^n|\hat k_i\cdot n(x_i)|\le2B_0\sum_{i=m+1}^n\vartheta_i$
which
shows that the $\hat k_n$ are Cauchy and the completeness of the $(d_1-1)$-dimensional Euclidean sphere implies that $\hat k_n\to\hat k_\infty$ for some $\hat k_\infty$ with
$\|\hat k_\infty\|=1$ and $n(\gamma'(x,k))\cdot \hat k_\infty=0$.
\end{definition}

\begin{definition}
Define for $r>0$,
\be
\partial\Omega^r=\{y\in\bar\Omega:d(y,\partial\Omega)<r\}.
\ee
where $d(y,A)=\inf\{\|y-a\|:a\in A\}$.
\end{definition}

\begin{proposition}\LABEL{prop:delta}
There exists $0<\delta<\h B_0^{-1}$ such that for all $x,y\in\partial\Omega^{2\delta}$ with $\|x-y\|<2\delta$ we have $\|n(x)-n(y)\|\le A_0\|x-y\|$ and if additionally $x,y\in\partial\Omega$ then $|n(x)\cdot(x-y)|\le B_0\|x-y\|^2$.
\end{proposition}
\begin{proof}
Since $\|\D\omega\|>1$ whenever $x\in\partial\Omega$, $\partial\Omega$ is compact and $x\mapsto \|\D\omega(x)\|$ is continuous, it follows that there exists $\eta>0$ such that
$\|\D\omega(x)\|\ge1+\eta$ for all $x\in\partial\Omega$. The uniform continuity of $\D\omega$ on $\bar\Omega$ implies that there exists $\delta_0>0$ such that if $\|x-y\|<\delta_0$ then
$\|\D\omega(x)-\D\omega(y)\|<\eta$ and thus $\|\D\omega(y)\|\ge\|\D\omega(x)\|-\|\D\omega(x)-\D\omega(y)\|>1$. Let $2\delta=\delta_0\wedge B_0^{-1}
$ and the result follows by proposition~\ref{prop:normal}.
\end{proof}

\begin{proposition}\LABEL{prop:grad-estimate}
For all $z\in\FR$ with $\pi_{\bar\Omega}(z)=x$, $k$ such that $\Sum_{r=0}^\infty\vartheta_r=\|k\|$ and $\gamma(x,k)\neq x_2$ we have that
\be
|(\gamma(x,k)-x-k)\cdot\D f(z)|\le2(2B_0+A_0)\|\D\D f\|_\infty\|k\|^3
\ee
for each $f\in\FF$.
 \end{proposition}
\begin{proof} Suppose that $x_1,\ldots,x_N\in\partial\Omega$ with $\gamma(x,k)=x_{N+1}$ and $N>1$ 
then
\beq
|(\gamma(x,k)-x-k)\cdot\D f(z)|&=&\left|\sum_{r=1}^N \vartheta_r(\hat k_r-\hat k_0)\cdot\D f(z)\right|\\
&=&\left|\sum_{r=1}^N\sum_{p=1}^r \vartheta_r(\hat k_{p}-\hat k_{p-1})\cdot\D f(z)\right|\\
&=&\left|2\sum_{r=1}^N\sum_{p=1}^r \vartheta_r(\hat k_{p}\cdot n(x_p))n(x_p)\cdot[\D f(z)-\D f(z_p)]\right|\\ 
&\le&2B_0\|\D\D f\|_\infty\|k\|\sum_{r=1}^N\sum_{p=1}^r\vartheta_r\vartheta_p+
2\|\D\D f\|_\infty\|k\|\vartheta_N |\hat k_N\cdot n(x_N)|\quad\LABEL{eq:twozero}\\
&\le&2B_0\|\D\D f\|_\infty\|k\|^3+2\|\D\D f\|_\infty\|k\|\vartheta_N (|\hat k_{N-1}\cdot n(x_{N-1})|+A_0\vartheta_{N-1})\quad\cr
&&\LABEL{eq:twotwo}\\
&\le&2(2B_0+A_0)\|\D\D f\|_\infty\|k\|^3
\eq
where $z_p=(\pim(z),x_p,\pi_\Gamma(z))$. We have made use of proposition~\ref{prop:normal} to derive \Ineq.twotwo. and that
\be
\|\D f(z)-\D f(z_p)\|\le\sum_{r=1}^p\|\D f(z_{r-1})-\D f(z_r)\|\le\|\D\D f\|_\infty\sum_{r=0}^{p-1}\vartheta_r\le\|\D\D f\|_\infty\|k\|,
\ee
noting that the straight line segment $[z_{r-1},z_r]$ is contained within $\FR$ to derive \Ineq.twozero.. In the event that no such $N$ exists the above proof shows that $|(x_{N+1}-x-k)\cdot\D f(z)|\le2(2B_0+A_0)\|\D\D f\|_\infty\|k\|^3$ for all $N>1$ and the conclusion follows from $\gamma(x,k)=\gamma'(x,k)=\lim_{N\to\infty}x_N$.
\end{proof}

\begin{proposition}\LABEL{prop:onereflection}
(a) If $\delta$ is as in proposition~\ref{prop:delta}, $x\in\partial\Omega^\delta$ and $|k\cdot n(x)|>(A_0+B_0)\|k\|^2$ 
then $\gamma(x,k)=x_2$.

\ni(b) If $\gamma(x,k)=x_2$ and $x=\pi_{\bar\Omega}(z)$ for some $z\in\FR$ then $|(\gamma(x,k)-x-k)\cdot\D f(z)|\le2\|\D\D f\|_\infty\|k\|^2$ for each $f\in\FF$.
\end{proposition}
\begin{proof} (a)
If $x_1\notin\partial\Omega$ then $\gamma(x,k)=x+k=x_1=x_2$ and there is nothing to prove, so we assume that $x_1\in\partial\Omega$.
It follows by proposition~\ref{prop:normal} that $\|n(x)-n(x_1)\|\le A_0\|k\|$ and we have
\be
\hat k_1=\hat k_0-2(\hat k_0\cdot n(x_1))n(x_1),
\ee
so that
\be
|\hat k_1\cdot n(x_1)|=|\hat k_0\cdot n(x_1)|\ge|\hat k_0\cdot n(x)|-\|n(x)-n(x_1)\|> B_0\|k\|.
\ee
Proposition~\ref{prop:normal} also implies that $\gamma(x,k)=\gamma(x_1,(\|k\|-\vartheta_0)\hat k_1)=x_1+(\|k\|-\vartheta_0)\hat k_1=x_2\notin\partial\Omega$.

(b) We have
\beq
|\{\gamma(x,k)-x-k\}\cdot\D f(z)|&=&|\{(\|k\|-\vartheta_0)(\hat k_1-\hat k_0)\}\cdot\D f(z)|\\
&=&|2\vartheta_1(\hat k_0\cdot n(x_1))n(x_1)\cdot\{\D f(z)-\D f(\pim(z),x_1,\pi_\Gamma(z))\}|\qquad\ 
\\
&\le&2\|\D\D f\|_\infty\|k\|^2
\eq
where we have made use of the boundary condition $n(x_1)\cdot \D f(\pim(z),x_1,\pi_\Gamma(z))=0$ and that the line segment $[x,x_1]$ lies within $\bar\Omega$.
\end{proof}

\begin{definition}
Define the map $j_{\rm ID}^N:[0,\infty)\times\FR\to\BR^{d_2}$ associated with the internal drift to be
\be
\(j_{\rm ID}^N(u,z)\)_\mu=a_\mu^{(\pim(z))}\vee(\pi_\Gamma(z)_\mu+\cN ^{-2}H_{u,\mu}(z))\wedge b_\mu^{(\pim(z))}\qquad\mbox{
for $\mu=1,\ldots,d_2$.}\LABEL{eq:JID}
\ee
\end{definition}


\begin{proposition}\LABEL{prop:internallimit}
If $f\in\FF_c$, $u^N_s\to s$ and $\cN \to\infty$ as $N\to\infty$ then
\be
\cN ^2\(j_{\rm ID}^N(u^N_s,z)-\pi_\Gamma(z)-\frc{H_{u^N_s}(z)}{\cN ^2}\)\cdot\hat\D f(z)\to0
\ee
uniformly on $\FR$ as $N\to\infty$.
\end{proposition}
\begin{proof}
We may assume that $\pim(z)\le R(f)$ as otherwise the terms in the sequence are identically zero. By the boundary conditions on the internal drift (hypothesis~\ref{hyp:internal}) and \Eq.JID., if $\cN ^2(j_{\rm ID}^N(u^N_s,z)-\pi_\Gamma(z))_\mu\neq H_\mu(u^N_s,z)$ then either $H_\mu(u^N_s,z)\ge \cN ^2(j_{\rm ID}^N(u^N_s,z)-\pi_\Gamma(z))_\mu$ and $(j_{\rm ID}^N(u^N_s,z))_\mu=b_\mu^{(\pim(z))}$ letting us deduce that $H_\mu(u^N_s,z+(j_{\rm ID}^N(u^N_s,z)-\pi_\Gamma(z))_\mu \hat e_{\mu})\le0$, or $H_\mu(u^N_s,z)\le \cN ^2(j_{\rm ID}^N(u^N_s,z)-\pi_\Gamma(z))_\mu$ and $(j_{\rm ID}^N(u^N_s,z))_\mu= a_\mu^{(\pim(z))}$ leading to $H_\mu(u^N_s,z+(j_{\rm ID}^N(u^N_s,z)-\pi_\Gamma(z))_\mu\hat  e_{\mu})\ge0$. Let $g_\mu(\lambda)=H_\mu(u^N_s,z+\lambda (j_{\rm ID}^N(u^N_s,z)-\pi_\Gamma(z))_\mu\hat e_{\mu})-\cN ^2(j_{\rm ID}^N(u^N_s,z)-\pi_\Gamma(z))_\mu$ then $g_\mu:[0,1]\to\BR$ is continuous and additionally $g_\mu(0)$ and $g_\mu(1)$ are either zero or have opposite sign.
Therefore there exists $\lambda_\mu\in[0,1]$ such that $g_\mu(\lambda_\mu)=0$. Take $\lambda_\mu=1$ if $\cN ^2(j_{\rm ID}^N(u^N_s,z)-\pi_\Gamma(z))_\mu= H_\mu(u^N_s,z)$. Set $\zeta^{(\mu)}=z+\lambda_\mu (j_{\rm ID}^N(u^N_s,z)-\pi_\Gamma(z))_\mu\hat e_{\mu})$, then we have that $\|\zeta^{(\mu)}-z\|\le\|H(u^N_s,z)\openone_{\pim(z)\le R(f)}\|\cN ^{-2}\le H^\infty_{R(f)}R(f) \cN ^{-2}$ and
\beq
\left|\cN ^2\(j_{\rm ID}^N(u^N_s,z)-\pi_\Gamma(z)-\frc{H(u^N_s,z)}{\cN ^2}\)\cdot\hat \D f(z)\right|\kern-4.8cm&&\cr
&\le&
\|\hat\D f\|_\infty\sum_{\mu=1}^{d_2}|\cN ^2(j_{\rm ID}^N(u^N_s,z)-\pi_\Gamma(z))_\mu-H_\mu(u^N_s,z)\openone_{\pim(z)\le R(f)}|\\
&\le&\|\hat\D f\|_\infty\sum_{\mu=1}^{d_2}|H_\mu(u^N_s,\zeta^{(\mu)})\openone_{\pim(\zeta^{(\mu)})\le R(f)}-H_\mu(u^N_s,z)\openone_{\pim(z)\le R(f)}|
\qquad\\
&\le&\|\hat\D f\|_\infty\sum_{\mu=1}^{d_2}w(H_\mu(\cdot,\cdot)\openone_{\pim(\cdot)\le R(f)},|u^N_s-s|+H^\infty_{R(f)}R(f)\cN ^{-2})\qquad
\eq
which tends to zero uniformly on $\FR$ as $N\to\infty$ using the uniform continuity of $H_\mu(\cdot,\cdot)\openone_{\pim(\cdot)\le R(f)}$ on $[0,s+1]\times\FR$ and the definition of the modulus of continuity (\Eq.modulus.).
\end{proof}

\begin{proposition}\LABEL{prop:Zestimate}
If $Z\sim N(0,\Id1)$, $n\in\BR^{d_1}$ with $\|n\|=1$, $\ell\ge0$, $m>1$ and $\cN \to\infty$ as $N\to\infty$ then there exist constants $D_1^{\ell,m}$ and $D_2^{\ell,m}$ such that
\be
\BE^Z(\|Z\|^\ell\openone_{|Z\cdot n|\le \cN ^{1-m}(A+B\|Z\|^m)})\le
\frc{D_1^{\ell,m}A+D_2^{\ell,m}B}{\cN ^{m-1}}.
\ee
\end{proposition}
\begin{proof}
If $d_1>1$ then
write $\|Z\|=r$ and $Z\cdot n=r\sin\theta$, then compute:
\beq
\BE^Z(\|Z\|^\ell\openone_{|Z\cdot n|\le \cN ^{1-m}(A+B\|Z\|^m)})\kern-3cm&&\cr
&=&\frc{2^{2-\frac12d_1}}{\sqrt\pi\Gamma(\h d_1-\h)}\int_0^\infty\int_0^{\sin^{-1}((\cN ^{1-m}(Ar^{-1}+Br^{m-1}))\wedge1)}\kern-2.5cm\cos^{d_1-2}\theta\,\d\theta\  r^{\ell+d_1-1}\e^{-\frac12r^2}\,\d r\\
&\le&\frc{\sqrt\pi}{2^{\frac12d_1-1}\Gamma(\h d_1-\h)\cN ^{m-1}}\int_0^\infty(Ar^{-1}+Br^{m-1})r^{\ell+d_1-1}\e^{-\frac12r^2}\,\d r\\
&\le&\frc{2^{\frac \ell2-\frac12}\sqrt\pi(A\Gamma(\h(\ell+d_1-1))+2^{\frac m2}B\Gamma(\h(\ell+m+d_1-1)))}{\Gamma(\h d_1-\h)\cN ^{m-1}},
\eq
whereas if $d_1=1$ and $\ell>0$ then
\beq
\BE^Z(\|Z\|^\ell\openone_{|Z\cdot n|\le \cN ^{1-m}(A+B\|Z\|^m)})&=&\sqrt{\frc2\pi}\int_0^\infty r^\ell\e^{-\frac12r^2}\openone_{r\le A+Br^m}\,\d r\\
&\le&\frc1{\cN ^{m-1}}\sqrt{\frc2\pi}\int_0^\infty r^\ell\e^{-\frac12r^2}(Ar^{-1}+Br^{m-1})\,\d r\\
&\le&\frc{2^{\frac \ell2-\frac12}(\Gamma(\h \ell)A+2^{\frac m2}\Gamma(\h(\ell+m))B)}{\sqrt\pi\,\cN ^{m-1}}.
\eq
The situation is more delicate when $d_1=1$ and $\ell=0$. First consider $B\neq0$, in this case we use that
\be
\openone_{r\le \cN ^{1-m}(A+Br^m)}\le\openone_{r\in[0,2Ac^{1-m}]\cup[\frac12 \cN B^{-\frac1{m-1}},\infty)}
\ee
which follows by noting that the function $f(r)=\cN ^{1-m}(A+Br^m)-r$ has exactly two positive roots with  $f(2A\cN ^{1-m}),f(\frac12 B^{-\frac1{m-1}}\cN )\le0$ provided that
\be
\cN >
2^{\frcn1{m-1}}A^{\frcn1{m}}B^{\frcn1{m(m-1)}}.
\LABEL{eq:bounds}
\ee
To prove this, the inequality $2(1-2^{-X})^X\ge1$ for $X>0$ is useful. It follows that
\beq
\BE^Z(\|Z\|^\ell\openone_{|Z\cdot n|\le \cN ^{1-m}(A+B\|Z\|^m)})&\le&\sqrt{\frc2\pi}\int_0^{2A\cN ^{1-m}} \kern-0.0cm\e^{-\frac12r^2}\,\d r\cr
&&\kern1cm{}+\sqrt{\frc2\pi}\int_{\frac12\cN B^{-\frac1{m-1}}}^\infty\kern-0.0cm \e^{-\frac12r^2}(2B^{\frac1{m-1}}\cN ^{-1}r)^{m-1}\,\d r\qquad\\
&\le&\frc{2^{\frac32}A+2^{\frac{3(m-1)}2}\Gamma(\h m) B}{\sqrt\pi\,\cN ^{m-1}}.
\eq
When $B=0$ the computation proceeds as above except that we note that there is only one zero of the function $f$. The finite number of times when \Ineq.bounds. may be violated can be dealt with by taking either $D_1^{\ell,m}$ or $D_2^{\ell,m}$ suitably large.
\end{proof}

\begin{corollary}\LABEL{cor:K2-bounds}
If $\cN \to\infty$ as $N\to\infty$ then there exist a function $F_\ell(A)$ such that
\be
\BE^Z\(\|k\|^\ell\openone_{|k\cdot n(x)|\le A\|k\|^2}\)\le\frc{F_\ell(A)}{\sigma^{\rm min}_{u,\pim(z)}\cN ^{\ell+1}}
\ee
for all $z\in\FR$, where $\ell\in\BN_0$, $x=\pi_{\bar\Omega}(z)$,  $Z\sim N(0,\Id1)$ and $
k=\frc{\sigma_u(z)Z}{\cN }+\frc{b_u(z)}{\cN ^2}$.
\end{corollary}
\begin{proof}
We make use of a special case of H\"older's inequality: $|X+Y|^\ell\le 2^{\ell-1}(|X|^\ell+|Y|^\ell)$ which follows for $\ell\ge1$ by considering the vectors $(X,Y)$ and $(1,1)$
together with the pair of conjugate indices $\ell$ and $\frac \ell{\ell-1}$. By inspection this inequality also holds for $\ell=0$.
If $|k\cdot n(x)|<A\|k\|^2$ then
\be
|Z\cdot n(x)|\le\frc{\alpha+\beta\|Z\|^2}{\cN }
\ee
where
\be
\alpha=\frc{b_\infty+2Ab^2_\infty\sup_N\cN ^{-2}}{\sigma_{u,\pim(z)}^{\rm min}}\and
\beta=\frc{2A\sigma_\infty^2}{\sigma_{u,\pim(z)}^{\rm min}}.
\ee
The result is completed by the use of proposition~\ref{prop:Zestimate}.
\end{proof}

\begin{proposition}\LABEL{prop:Taylor}
The following is a (Taylor's theorem) identity
\beq
F(x+h,y+k)-F(x,y)&=&\int_0^1hD_1 F(x,y+tk)\,\d t+\int_0^1kD_2F(x+th,y+tk)\,\d t\cr
&&+\int_0^1\int_0^s h^2D_1^2F(x+th,y+sk)\,\d t\,\d s
\eq
where  $F:\BR^2\to\BR$ is continuously differentiable, has continuous second order partial derivative in the first argument and the domain of $f$ contains the closed cuboid consisting of the points with coordinates $(x+th,y+sk)$ for $s,t\in[0,1]$.
\end{proposition}
\begin{proof}
Write $t=sr$ in the last integral on the RHS and perform the integrals.
\end{proof}

\begin{remark}\em\LABEL{remark:GN}
Proposition~\ref{prop:Taylor} implies that
if $k=\frc{\sigma_u(z)Z}{\cN }+\frc{b_u(z)}{\cN ^2}$, $x=\pi_{\bar\Omega}(z)$ and
\beq
G^N(u,z,Z)&=&\int_0^1[\gamma(x,k)-x-k]\cdot\D f(\pim(z),x,(1-t)\pi_\Gamma(z)+tj_{\rm ID}^N(u,z))\,\d t\cr
&&\kern-0.5cm{}+\int_0^1k\cdot\(\D f(\pim(z),x,(1-t)\pi_\Gamma(z)+tj_{\rm ID}^N(u,z))-\D f(z)\)\,\d t\cr
&&\kern-0.5cm{}+\sum_{i,j=1}^{d_1}\int_0^1\int_0^s(\gamma(x,k)-x)_i(\gamma(x,k)-x)_j\cr
&&\kern0cm\times
\(\partial_i\partial_j\fext(\pim(z),(1-t)x+t\gamma(x,k),(1-s)\pi_\Gamma(z)+sj_{\rm ID}^N(u,z))
-\partial_i\partial_jf(z)\)\,\d t\,\d s
\cr
&&\kern-0.5cm{}+\sum_{i,j=1}^{d_1}\h(\gamma(x,k)-x-k)_i(\gamma(x,k)-x+k)_j\partial_{i}\partial_jf(z)\cr
&&\kern-0.5cm{}+\int_0^1[j_{\rm ID}^N(u,z)-\pi_\Gamma(z)]\cdot\(\hat\D \fext(\pim(z),(1-t)x+t\gamma(x,k),(1-t)\pi_\Gamma(z)+tj_{\rm ID}^N(u,z))\right.\cr
&&\left.\kern10cm{}-\hat\D f(z)\)\,\d t\cr
&&\kern-0.5cm{}+\(j_{\rm ID}^N(u,z)-\pi_\Gamma(z)-\frc{H_u(z)}{\cN ^2}\)\cdot\hat\D f(z)+\sum_{i,j=1}^{d_1}\frc{b_ib_j}{2\cN ^4}\partial_i\partial_jf(z)
\eq
then if $Z\sim N(0,\Id1)$,
\beq
\BE^Z\( \cN ^2G^N(u,z,Z)\)
&=& \BE^Z\(\cN ^2
f\(\pim(z),\gamma(x,k),j_{\rm ID}^N(u,z)\)-f(z)\)
\cr
&&\kern2cm{}-\h a_u(z)\triangle f(z)-b_u(z)\cdot\D f(z)-H_u(z)\cdot\hat\D f(z).\qquad
\eq
\end{remark}

\begin{proposition}\LABEL{prop:boundedness}
If $f\in\FF_c\cup\{f_0\}$, $\cN \to\infty$ as $N\to\infty$ then there exists some constant $L$ such that for all $z\in\FR$
\be
\left|\cN ^2\BE^Z\(f(\pim(z),\gamma(\pi_{\bar\Omega}(z),k),j_{\rm ID}^N(u,z))-f(z)\)\right|\le L\pim(z)H^\infty_{R(f)}\:\fext\:^\#\LABEL{eq:estimate1}
\ee
for all $N$ where $Z\sim N(0,\Id1)$ and
$
k=\frc{\sigma_u(z)Z}{\cN }+\frc{b_u(z)}{\cN ^2}.
$
\end{proposition}
\begin{proof}
If $f=f_0$ there is nothing to prove, so we assume that $f\in\FF_c$.
If $x=\pi_{\bar\Omega}(z)\in\partial\Omega^\delta$ then propositions~\ref{prop:grad-estimate} and  \ref{prop:onereflection}(b) imply
 that
\beq
\left|\cN ^2\BE^Z\(\{\gamma(x,k)-x-k\}\cdot \D f(z)\)\right|\kern-4cm&&\cr
&\le&2\cN ^2\|\D f\|_\infty\BE^Z(\openone_{\|k\|>\delta}\|k\|)
+\BE^Z(\cN ^2\{\gamma(x,k)-x-k\}\cdot \D f(z)\openone_{\gamma(x,k)=x_2})\cr
&&\kern3cm{}+\BE^Z(\cN ^2\{\gamma(x,k)-x-k\}\cdot \D f(z)\openone_{\gamma(x,k)\neq x_2})\qquad\\
&\le&
2(\delta^{-1}\|\D f\|_\infty+\|\D\D f\|_\infty+(2B_0+A_0)\delta^{-1}\|\D\D f\|_\infty)\BE^Z(\cN ^2\|k\|^2)\\
&\le&4(\delta^{-1}\|\D f\|_\infty+\|\D\D f\|_\infty+(2B_0+A_0)\delta^{-1}\|\D\D f\|_\infty)(\sigma_\infty^2d_1+b_\infty^2\sup_N{\cN ^{-2}})\qquad
\eq
whereas if $x\not\in\partial\Omega^\delta$ then
\beq
\left|\cN ^2\BE^Z\(\{\gamma(x,k)-x-k\}\cdot \D f(z)\)\right|
&\le&2\cN ^2\|\D f\|_\infty\BE^Z(\openone_{\|k\|>\delta}\|k\|)\\
&\le&2\|\D f\|_\infty\delta^{-1}\BE^Z(\cN ^2\|k\|^2)\\
&\le&4\|\D f\|_\infty\delta^{-1}(\sigma_\infty^2d_1+b_\infty^2\sup_N\cN ^{-2}).
\eq
Note that
\beq
\left|\cN ^2\BE^Z\(k\cdot\D f(z)\)\right|&=&\left|\cN \BE^Z\(\sigma_u(z)Z\cdot\D f(z)\)+b_u(z)\cdot\D  f(z)\right|\\
&=&|b_u(z)\cdot\D  f(z)|\\
&\le&b_\infty\|\D f\|_\infty
\eq
so that if we set
\be
L_0=4(1+\delta^{-1}+(2B_0+A_0)\delta^{-1})(\sigma_\infty^2d_1+b_\infty^2\sup_N{\cN ^
{-2}})+b_\infty
\ee
then Taylor's theorem in the form of proposition~\ref{prop:Taylor} and hypothesis~\ref{hyp:internal} give the upper bound
\beq
\left|
\BE^Z\(\cN ^2[
f\(\pim(z),\gamma\(x,k\),j_{\rm ID}^N(u,z)\)-f(z)]\)
\right|&&\cr
&&\kern-6cm{}\le \left|
\BE^Z\(\cN ^2
\{\gamma(x,k)-x\}\cdot\D f(z)\)
\right|+\h \|\D\D \fext\|_\infty^\#\BE^Z(\cN ^2\|k\|^2)\cr
&&\kern-4cm{}+\cN ^2\|\hat\D \fext\|^\#_\infty
\BE^Z\(
\|j_{\rm ID}^N(u,z)-\pi_\Gamma(z)\|\openone_{[1,R(f)]}(\pim(z))\)
\\
&&\kern-6cm{}\le L_0\:\fext\:^\#+\h \|\D\D \fext\|^\#_\infty\(\sigma_\infty^2d_1+b_\infty^2\sup_N\cN ^{-2}\)
+H^\infty_{R(f)}\|\hat\D \fext\|^\#_\infty \pim(z)\qquad\\
&&\kern-6cm{}\le L\pim(z)H^\infty_{R(f)}\:\fext\:^\#\LABEL{eq:bded}
\eq
where
\be
L=L_0+1+\h(\sigma_\infty^2d_1+b_\infty^2\sup_N\cN ^{-2}).
\LABEL{eq:L_1}
\ee
\end{proof}

\begin{proposition}\LABEL{prop:uniform-result}
Suppose $f\in\FF_c$, $\cN \to\infty$ as $N\to\infty$, $Z\sim N(0,\Id1)$ and $z^t=(\pim(z),x,(1-t)\pi_\Gamma(z)+tj^N_{\rm ID}(u^N_s,z))$ for $t\in[0,1]$ and  $k=\frc{\sigma_{u^N_s}(z)Z}{\cN }+\frc{b_{u^N_s}(z)}{\cN ^2}$ then for each $s\in[0,\infty)$,
\be
I^N(z,t)=\sum_{i=1}^{d_1}\cN ^2\BE^Z\([\gamma(x,k)-x-k]_i[\partial_i f(z^t)+\sum_{j=1}^{d_1}[\h(\gamma(x,k)-x+k)]_j\partial_i\partial_jf(z)]\)\to0
\ee
uniformly on $\FR\times[0,1]$ as $N\to\infty$.
\end{proposition}
\begin{proof}
We may assume that $N>N_0$ where for all $N$ with $N>N_0$ we have that $v^N_s\le s+1$. 
Set $\FR^\delta=\pi_{\bar\Omega}^{-1}(\partial\Omega^\delta)$, if $z\notin\FR^\delta$ then
\beq
|I^N(z,t)|&\le&2\delta^{-2}\BE^Z\(\cN ^2\|k\|^3\|\D f\|_\infty+\cN ^2\|k\|^3\delta\|\D\D f\|_\infty\)\\
&\le&\frc{8(\|\D f\|_\infty+\delta\|\D\D f\|_\infty)}{\delta^2\cN }\(\sigma_\infty^3\BE^Z(\|Z\|^3)+\frc{b_\infty^3}{\cN ^3}\)\to0\LABEL{eq:outer2}
\eq
uniformly on $(\FR\setminus\FR^\delta)\times[0,1]$. Now consider the more tricky case when $z\in\FR^\delta$. Let $z_1^t=(\pim(z),x_1(k),(1-t)\pi_\Gamma(z)+tj^N_{\rm ID}(u^N_s,z))$. We drop the $(z,t)$ dependence for notational convenience and write $I^N=I_1^N+I_2^N+I_3^N$ where
\beq
I^N_1&=&\sum_{i=1}^{d_1}\cN ^2\BE^Z\([\gamma(x,k)-x-k]_i[\partial_i f(z^t)+\sum_{j=1}^{d_1}[\h(\gamma(x,k)-x+k)]_j\partial_i\partial_jf(z)]\openone_{\|k\|\ge\delta}\);\\
I^N_2&=&\cN ^2\BE^Z\([\gamma(x,k)-x-k]\cdot\D f(z_1^t)\openone_{\|k\|<\delta}\);\\
I^N_3&=&\sum_{i=1}^{d_1}\cN ^2\BE^Z\([\gamma(x,k)-x-k]_i[\partial_i f(z^t)-\partial_if(z^t_1)+\sum_{j=1}^{d_1}[\h(\gamma(x,k)-x+k)]_j\partial_i\partial_jf(z)]\openone_{\|k\|<\delta}\),\cr&&
\eq
then $I^N_1$ obeys the same bound as \Ineq.outer2. and therefore vanishes uniformly on $\FR^\delta\times[0,1]$ in the limit $N\to\infty$. If $\gamma(x,k)= x_2$ then
$\gamma(x,k)-x-k=-2\vartheta_1\hat k_0\cdot n(x_1)n(x_1)$ and therefore the boundary condition on the class of test functions: $n(x_1)\cdot\D f(z_1^t)=0$ (definition~\ref{def:testfunctions}) implies that
\be
I^N_2=\cN ^2\BE^Z\([\gamma(x,k)-x-k]\cdot\D f(z_1^t)\openone_{\gamma(x,k)\neq x_2}\openone_{\|k\|<\delta}\).\\
\ee
Proposition~\ref{prop:onereflection} implies that $\gamma(x,k)\neq x_2$ only if $|n(x)\cdot k|\le(A_0+B_0)\|k\|^2$. If $\gamma(x,k)\neq x_2$ and $\sum_{n=0}^\infty\vartheta_n=\|k\|$ then
 proposition~\ref{prop:grad-estimate} gives an estimate for $|(\gamma(x,k)-x-k)\cdot\D f(z)|$. If $\gamma(x,k)\neq x_2$ and $\sum_{n=0}^\infty\vartheta_n\neq\|k\|$ then we may still use proposition~\ref{prop:grad-estimate} and add an extra term given by
proposition~\ref{prop:xi-estimate}. In consequence
\beq
|(\gamma(x,k)-x-k)\cdot\D f(z_1^t)|&\le&|(\gamma(x,k)-x-k)\cdot\D f(z^t)|+|(\gamma(x,k)-x-k)\cdot[\D f(z_1^t)-\D f(z^t)]|\cr
&&\\
&\le&2(2B_0+A_0)\|\D\D f\|_\infty\(\sum_{n=0}^\infty\vartheta_n\)^3+2B_0\(\|k\|-\sum_{n=0}^\infty\vartheta_n\)^2\|\D f\|_\infty\cr&&\kern6cm{}+2\|\D\D f\|_\infty\|k\|^2\\
&\le&[2\{(2B_0+A_0)\delta+1\}\|\D\D f\|_\infty+2B_0\|\D f\|_\infty]\|k\|^2.
\eq
Accordingly
\beq
|I^N_2|&\le& [2(2B_0+A_0)\delta\|\D\D f\|_\infty+2B_0\|\D f\|_\infty]\BE^Z(\cN ^2\|k\|^2\openone_{|n(x)\cdot k|\le(A_0+B_0)\|k\|^2})\\
&\le&\frc{[2(2B_0+A_0)\delta\|\D\D f\|_\infty+2B_0\|\D f\|_\infty] F_2(A_0+B_0)}{\cN \sigma^{\rm min}_{s+1,R(f)}}\to0
\eq
uniformly on $\FR^\delta\times[0,1]$ using corollary~\ref{cor:K2-bounds}. Split $I^N_3=I^N_4+I^N_5$ where
\beq
I^N_4&=&\sum_{i=1}^{d_1}\cN ^2\BE^Z\([\gamma(x,k)-x-k]_i[\partial_i f(z^t)-\partial_if(z^t_1)+\sum_{j=1}^{d_1}[\h(\gamma(x,k)-x+k)]_j\partial_i\partial_jf(z)]\right.\cr
&&\kern8.5cm{}\left.\phantom{\sum_{n=0}^\infty}{}\times\openone_{\gamma(x,k)\neq x_2}\openone_{\|k\|<\delta}\);\\
I^N_5&=&\sum_{i=1}^{d_1}\cN ^2\BE^Z\([\gamma(x,k)-x-k]_i[\partial_i f(z^t)-\partial_if(z^t_1)+\sum_{j=1}^{d_1}[\h(\gamma(x,k)-x+k)]_j\partial_i\partial_jf(z)]\right.\qquad\qquad\cr
&&\kern8.5cm{}\left.\phantom{\sum_{n=0}^\infty}{}\times\openone_{\gamma(x,k)=x_2}\openone_{\|k\|<\delta}\)
\eq
and notice that
\beq
|I^N_4|&\le& 4\|\D\D f\|_\infty\BE^Z\(\cN ^2\|k\|^2\openone_{|n(x)\cdot k|\le(A_0+B_0)\|k\|^2}\)\\
&\le&\frc{4\|\D\D f\|_\infty F_2(A_0+B_0)}{\cN \sigma^{\rm min}_{s+1,R(f)}}\to0\LABEL{eq:bound4}
\eq
uniformly on $\FR^\delta\times[0,1]$ using corollary~\ref{cor:K2-bounds}. Now decompose $I^N_5$ as $I^N_5=I^N_6+I^N_7+I^N_8+I_9^N$ where
\beq
I^N_6&=&2\sum_{i,j=1}^{d_1}\BE^Z\(\cN ^2\frc{\vartheta_1^2}{\|k\|^2}k\cdot n(x_1)(k-k\cdot n(x_1)n(x_1))_jn(x_1)_i\partial_i\partial_j f(z)\openone_{\gamma(x,k)\neq x_2}\openone_{\|k\|<\delta}\);\\
I^N_7&=&2\sum_{i,j=1}^{d_1}\BE^Z\(\cN ^2\frc{\vartheta_1^2}{\|k\|^2}k\cdot n(x_1)(k-k\cdot n(x_1)n(x_1))_jn(x_1)_i\partial_i\partial_j f(z)\openone_{\|k\|\ge\delta}\);\\
I^N_8&=&\sum_{i=1}^{d_1}\BE^Z\(\cN ^2[\gamma(x,k)-x-k]_i\(\partial_if(z^t)-\partial_if(z^t_1)\phantom{\sum_{i=1}^{d_1}}\right.\right.\cr
&&{}+\left.\left.\sum_{j=1}^{d_1}\left[\h(\gamma(x,k)-x+k)-\frc{\vartheta_1}{\|k\|}(k-k\cdot n(x_1)n(x_1))\right]_j\partial_i\partial_jf(z)\)\openone_{\gamma(x,k)=x_2}\openone_{\|k\|<\delta}\)
\qquad\\
&=&\sum_{i=1}^{d_1}\BE^Z\(\cN ^2[\gamma(x,k)-x-k]_i\(\partial_if(z^t)-\partial_if(z^t_1)
+\sum_{j=1}^{d_1}[
x_1-x
]_j\partial_i\partial_jf(z)\)\right.\cr&&\left.\kern9cm\phantom{\sum_{i=1}^{d_1}}{}\times\openone_{\gamma(x,k)=x_2}\openone_{\|k\|<\delta}\)
\qquad\\
\mbox{and}&&\cr
I^N_9&=&-2\sum_{i,j=1}^{d_1}\BE^Z\(\cN ^2\frc{\vartheta_1^2}{\|k\|^2}k\cdot n(x_1)(k-k\cdot n(x_1)n(x_1))_jn(x_1)_i\partial_i\partial_j f(z)\).
\eq
We observe that $|I^N_6|$ is bounded by a half of the RHS of \Ineq.bound4. and therefore it too tends to zero uniformly on $\FR^\delta\times[0,1]$.
We bound $I^N_7$ easily:
\be
|I^N_7|\le2\delta^{-1}\|\D\D f\|_\infty\BE^Z\(\cN ^2\|k\|^3\)\le\frc{8\|\D\D f\|_\infty}{\delta \cN }\(\sigma_\infty^3\BE^Z(\|Z\|^3)+\frc{b_\infty^3}{\cN ^3}\)\to0
\ee
uniformly on $\FR^\delta\times[0,1]$. Note that
\beq
\left|\partial_if(z^t)-\partial_if(z^t_1)
+\sum_{j=1}^{d_1}[x_1-x]_j\partial_i\partial_jf(z)\right|\kern-1pt&=&\kern-1pt\left|\sum_{j=1}^{d_1}[x-x_1]_j\int_0^1[\partial_i\partial_j f(sz^t+(1-s)z^t_1)-\partial_i\partial_jf(z)]\,\d s\right|\cr
&&\\
&\le&\|k\|\sum_{j=1}^{d_1}w(\partial_i\partial_jf,\|k\|+\cN ^{-2}H^\infty_{R(f)}R(f))
\eq
and
\beq
\sum_{j=1}^{d_1}w(\partial_i\partial_jf,\|k\|+\cN ^{-2}H^\infty_{R(f)}R(f))&\le&\sum_{j=1}^{d_1}w(\partial_i\partial_jf,\sigma_\infty\|Z\|\cN ^{-1}+b_\infty \cN ^{-2}+\cN ^{-2}H^\infty_{R(f)}R(f))\qquad
\eq
which is bounded and, by the uniform continuity of $\partial_i\partial_jf$, tends to zero uniformly on $\FR^\delta\times[0,1]$ at each $Z\in\BR^{d_1}$. Therefore
\beq
|I^N_8|&\le&2\BE^Z\(\cN ^2\|k\|^2\sum_{j=1}^{d_1}w(\partial_i\partial_jf,\sigma_\infty\|Z\|\cN ^{-1}+b_\infty \cN ^{-2}+\cN ^{-2}H^\infty_{R(f)}R(f))\)\\
&\le&4\BE^Z\([\sigma_\infty^2\|Z\|^2+b_\infty^2 \cN ^{-2}]\sum_{j=1}^{d_1}w(\partial_i\partial_jf,\sigma_\infty\|Z\|\cN ^{-1}+b_\infty \cN ^{-2}+\cN ^{-2}H^\infty_{R(f)}R(f))\)\qquad
\eq
which goes to zero uniformly on $\FR\times[0,1]$ by the dominated convergence theorem. Let
\be
I^N_{10}=-2\sum_{i,j=1}^{d_1}\BE^Z\(\cN ^2\frc{\vartheta_1^2}{\|k\|^2}k\cdot n(x)(k-k\cdot n(x)n(x))_jn(x)_i\partial_i\partial_j f(z)\).
\ee
Then proposition~\ref{prop:normal} implies that
\beq
|I_{10}^N-I^N_9|&\le& 8A_0\BE^Z(\cN ^2\|k\|^3)\|\D\D f\|_\infty\cr
&\le& 32A_0\cN ^{-1}[\sigma_\infty^3\BE^Z\(\|Z\|^3\)+b_\infty^3\cN ^{-3}]\|\D\D f\|_\infty\to0
\eq
uniformly on $\FR^\delta\times[0,1]$.
We need to show that $I_{10}^N\to0$ uniformly, to this end set $\bar Z$ to be the rotation of $Z$ by angle $\pi$ about the axis $n(x)$:
\be
\bar Z=2Z\cdot n(x)n(x)-Z,
\ee
and notice that for any integrable function $g(Z)$ we have $\BE^Z(g(Z))=\BE^Z(\h[g(Z)+g(\bar Z)])$. We also define
\be
\bar k=\frc{\sigma_{u^N_s}(z)\bar Z}{\cN }+\frc{b_{u^N_s}(z)}{\cN ^2},
\ee
so that $\bar k\cdot n(x)=k\cdot n(x)$ and
\be
\bar k-\bar k\cdot n(x)n(x)=-(k-k\cdot n(x)n(x))+\frc{2(b_{u^N_s}(z)-b_{u^N_s}(z)\cdot n(x)n(x))}{\cN ^2}.
\ee
Thus $I^N_{10}=I^N_{11}+I^N_{12}+I^N_{13}+I^N_{14}$ where
\beq
I^N_{11}&=&-2\sum_{i,j=1}^{d_1}\BE^Z\(k\cdot n(x)\frc{\bar\vartheta_1^2}{\|\bar k\|^2}(b_{u^N_s}(z)-b_{u^N_s}(z)\cdot n(x)n(x))_jn(x)_i\partial_i\partial_jf(z)\);\\
I^N_{12}&=&\sum_{i,j=1}^{d_1}\BE^Z\(\cN ^2k\cdot n(x)\(\frc{\bar\vartheta_1^2}{\|\bar k\|^2}-\frc{\vartheta_1^2}{\|k\|^2}\)(k-k\cdot n(x)n(x))_jn(x)_i\partial_i\partial_jf(z)
 \right.\cr
&&\kern8cm\left.\phantom{\frc\vartheta\|}{}\times
\openone_{
k\cdot n(x)<-B_0\{\|k\|\vee\|\bar k\|\}^2}
\);
\\
I^N_{13}&=&\sum_{i,j=1}^{d_1}\BE^Z\(\cN ^2k\cdot n(x)\(\frc{\bar\vartheta_1^2}{\|\bar k\|^2}-\frc{\vartheta_1^2}{\|k\|^2}\)(k-k\cdot n(x)n(x))_jn(x)_i\partial_i\partial_jf(z)
\right.\cr
&&\kern7.5cm\left.\phantom{\frc\vartheta\|}{}\times
\openone_{-B_0\{\|k\|\vee\|\bar k\|\}^2\le k\cdot n(x)\le 0}\);\\
&&\kern-1.9cm\mbox{and}\cr
I^N_{14}&=&\sum_{i,j=1}^{d_1}\BE^Z\(\cN ^2k\cdot n(x)\(\frc{\bar\vartheta_1^2}{\|\bar k\|^2}-\frc{\vartheta_1^2}{\|k\|^2}\)(k-k\cdot n(x)n(x))_jn(x)_i\partial_i\partial_jf(z)\openone_{k\cdot n(x)>0}\).\qquad
\eq
We continue to bound:
\be
|I^N_{11}|\le 2b_\infty\|\D\D f\|_\infty\BE^Z\(\|k\|\)\le\frc{2b_\infty\|\D\D f\|_\infty}{\cN }\(\sigma_\infty\BE^Z(\|Z\|)+\frc{b_\infty}{\cN }\)\to0
\ee
uniformly on $\FR^\delta\times[0,1]$ as $N\to\infty$. Furthermore if $k\cdot n(x)<-B_0\|k\|^2$ then for all $r\in[0,1]$ we have
that
\beq
\omega(x+rk)&=&\omega(x)+rk\cdot\D\omega+
\sum_{i,j=1}^{d_1}r^2\int_0^1(1-s)k_ik_j\partial_i\partial_j\omega(x+srk)\,\d s\\
&\le&\omega(x)-r\left\{\|\D\omega(x)\|B_0-\h r\|\D\D\omega\|_\infty\right\}\|k\|^2\le0
\eq
where we have recalled that $B_0=\h\|\D\D\omega\|_\infty$ and $\|\D\omega(x)\|\ge1$ as $z\in\FR^\delta$. It follows that $\vartheta_0=\|k\|$ and therefore that $\vartheta_1=0$.
This and the corresponding remark for $\bar\vartheta_1$ imply that $I^N_{12}=0$. Turning our attention to $I^N_{13}$ we derive the bound
\beq
|I^N_{13}|&\le&\|\D\D f\|_\infty\BE^Z\(\cN ^2\|k\|^2\openone_{|k\cdot n|<B_0\|k\|^2}\)+\|\D\D f\|_\infty\BE^Z\(\cN ^2\|k\|^2\openone_{|k\cdot n|<B_0\|\bar k\|^2}\)\\
&\le&\frc{2\|\D\D f\|_\infty F_2(B_0)}{\cN \sigma^{\rm min}_{s+1,R(f)}}+2\sigma_\infty^2\|\D\D f\|_\infty\BE^Z(\|Z\|^2\openone_{|Z\cdot n(x)|
<(2B_0b_\infty^2\cN ^{-2}+b_\infty+2B_0\sigma_\infty^2\|Z\|^2)
/\cN\sigma^{\rm min}_{s+1,R(f)}})
\cr&&\kern8cm{}+\frc{2\|\D\D f\|_\infty b_\infty^2}{\cN ^2}.
\eq
This tends to zero uniformly on $\FR^\delta\times[0,1]$ by using proposition~\ref{prop:Zestimate}.

Next we consider the expression
\be
\(\frc{\vartheta_1^2}{\|k\|^2}-\frc{\bar\vartheta_1^2}{\|\bar k\|^2}\)k\cdot n(x).
\ee
If $\vartheta_1,\bar\vartheta_1\neq0$ we deduce that
\beq
\left|\(\frc{\vartheta_1^2}{\|k\|^2}-\frc{\bar\vartheta_1^2}{\|\bar k\|^2}\)k\cdot n(x)\right|&=&\kern-1pt\left|\frc1{\|\D\omega(x)\|}\(\frc{\vartheta_1}{\|k\|}+\frc{\bar\vartheta_1}{\|\bar k\|}\)\sum_{i,j=1}^{d_1}
\left[\frc{\vartheta_0^2}{\|k\|^2}k_ik_j\int_0^1(1-s)\partial_i\partial_j\omega(x+s\vartheta_0k)\,\d s\right.\right.\cr
&&\left.\left.\kern1cm{}
-\frc{\bar\vartheta_0^2}{\|\bar k\|^2}\bar k_i\bar k_j\int_0^1(1-s)\partial_i\partial_j\omega(x+s\bar\vartheta_0\bar k)\,\d s\right]\right|
\\
&\le&\|\D\D\omega\|_\infty(\|k\|^2+\|\bar k\|^2)\le2(A_0+B_0)(\|k\|+\|\bar k\|)^2,
\eq
whereas if $\bar\vartheta_1=0$ and $\vartheta_1\neq0$ so that $\omega(x+\vartheta_0\hat k_0)=0$ and $\bar\vartheta_1=0$, we may set $x_1=x+\vartheta_0\hat k_0$ and $K=\bar k-\vartheta_0\hat k_0$ to derive
\beq
\frc{\vartheta_1}{\|k\|}k\cdot n(x)&=&\bar k\cdot n(x)-\frc{\vartheta_0}{\|k\|}k\cdot n(x)\\
&=&K\cdot n(x_1)+K\cdot (n(x)-n(x_1))\\
&\le&A_0\|K\|\,\|k\|+B_0\|K\|^2\le(A_0+B_0)(\|k\|+\|\bar k\|)^2
\eq
where we have used proposition~\ref{prop:normal} and the result that if $K\cdot n(x_1)>B_0\|K\|^2$ then
\be
\omega(x+\bar k)=\omega(x_1+K)\ge\omega(x_1)+K\cdot\D\omega(x_1)-B_0\|K\|^2>0.
\ee
This allows us to conclude that $\omega(x+\bar k)\notin \bar\Omega$. Accordingly if $k\cdot n(x)>0$ then
\be
\left|\(\frc{\vartheta_1^2}{\|k\|^2}-\frc{\bar\vartheta_1^2}{\|\bar k\|^2}\)k\cdot n(x)\right|
\le2(A_0+B_0)(\|k\|+\|\bar k\|)^2.
\ee
In consequence we are led to a bound for $I^N_{14}$:
\be
|I^N_{14}|\le\frc{8(A_0+B_0)\|\D\D f\|_\infty}{\cN }\BE^Z\(\left[\sigma_\infty\|Z\|+\frc{b_\infty}{\cN }\right]^3\)\to0
\ee
uniformly on $\FR^\delta\times[0,1]$.
\end{proof}

\begin{theorem}\LABEL{thm:uniform-internal}
If $f\in\FF_c\cup\{f_0\}$ with extension $\fext\in C^2(\FR^\#)$, $u^N_s\to s$ and $\cN \to\infty$ as $N\to\infty$ then
\beq
&&\lim_{N\to\infty}\  \BE^Z\(\cN ^2
[f\(\pim(z),\gamma\(x,k\),j_{\rm ID}^N(u^N_s,z)\)-f(z)]\)\cr
&&\kern4cm{}=\h a_s(z)\triangle f(z)+b_s(z)\cdot\D f(z)+H_s(z)\cdot\hat\D f(z).
\eq
uniformly on $\FR$ where $x=\pi_{\bar\Omega}(z)$, $Z\sim N(0,\Id1)$ and $k=\frc{\sigma_{u^N_s}(z)Z}{\cN }+\frc{b_{u^N_s}(z)}{\cN ^2}$.
\end{theorem}
\begin{proof} The case $f=f_0$ is trivial, so consider $f\in\FF_c$.
By remark~\ref{remark:GN}, we must prove that $\BE^Z\(\cN ^2G^N(u^N_s,z,Z)\)\to0$ uniformly on $\FR$ as $N\to\infty$. In terms of the modulus of continuity of $\hat \partial_i\fext$ we compute
\beq
&&\kern-1cm\left|\cN ^2[j_{\rm ID}^N(u,z)-\pi_\Gamma(z)]\cdot\int_0^1
\BE^Z
\(\hat\D \fext(\pim(z),(1-t)x+t\gamma(x,k),(1-t)\pi_\Gamma(z)+tj_{\rm ID}^N(u,z))-\hat\D f(z)\)
\,\d t\right|\cr
&&\kern2cm\le H^\infty_{R(f)}R(f) \sum_{\mu=1}^{d_2}\BE^Z\( w\(\hat\partial_\mu\fext, \frc{\sigma_\infty\|Z\|}{\cN }+\frc{b_\infty}{\cN ^2}
+\frc{H^\infty_{R(f)}R(f)}{\cN ^2}\)\)
\to0
\eq
uniformly on $\FR$ by noting that
\be
w\(\hat\partial_\mu\fext,
\frc{\sigma_\infty\|Z\|}{\cN }+\frc{b_\infty}{\cN ^2}+\frc{H^\infty_{R(f)}R(f) }{\cN ^2}\)\le2\:\hat\D \fext\:^\#
\ee
and using the dominated convergence theorem together with the fact that $\hat\partial_i \fext$ is uniformly continuous. Similarly,
\beq
&&\left|\cN ^2\BE^Z(k)\cdot\int_0^1\(\D f(\pim(z),x,(1-t)\pi_\Gamma(z)+tj_{\rm ID}^N(u,z))-\D f(z)\)\,\d t\right|\cr
&&\kern3cm=\left|b_u(z)\cdot\int_0^1\(\D f(\pim(z),x,\pi_\Gamma(z)(1-t)+j_{\rm ID}^N(u,z)t)-\D f(z)\)\,\d t\right|\qquad\qquad\\
&&\kern3cm\le b_\infty\sum_{i=1}^{d_1} w\(\partial_if,\frc{H^\infty_{R(f)}R(f)}{\cN ^2}\)\to0
\eq
uniformly on $\FR$ as $N\to\infty$,
\beq
&&\kern-0.5cm\left|\cN ^2\sum_{i,j=1}^{d_1}\int_0^1\int_0^s\BE^Z\((\gamma(x,k)-x)_i(\gamma(x,k)-x)_j
\right.\right.\cr
&&\kern2cm\left.\vphantom\int{}\left.\times\(\partial_i\partial_j\fext(\pim(z),(1-t)x+t\gamma(x,k),(1-s)\pi_\Gamma(z)+sj_{\rm ID}^N(u,z))-\partial_i\partial_jf(z)\)\)\,\d t\,\d s
\right|\cr
&&\le\left|
\BE^Z(\cN ^2\|k\|^2)
\sum_{i,j=1}^{d_1}w\(\partial_i\partial_j\fext,\frc{\sigma_\infty\|Z\|}{\cN }+\frc{b_\infty}{\cN ^2}
+\frc{H^\infty_{R(f)}R(f)}{\cN ^2}\)
\right|\to0
\eq
and
\be
\left|\cN ^2\sum_{i,j=1}^{d_1}\frc{b_ib_j\partial_i\partial_jf(z)}{2\cN ^4}\right|\le\frc{b_\infty^2\|\D\D f\|_\infty}{2\cN ^2}\to0
\ee
uniformly on $\FR$ as $N\to\infty$. 
Using proposition~\ref{prop:internallimit} it remains only to prove that
\beq
&&
\left|\cN ^2
\int_0^1
\BE^Z\([\gamma(x,k)-x-k]\cdot\D f(\pim(z),x,\pi_\Gamma(z)(1-t)+j_{\rm ID}^N(u,z)t)\)
\,\d t
\right.\cr&&\left.\kern3cm{}+\h \cN ^2\sum_{i,j=1}^{d_1} \BE^Z\(
(\gamma(x,k)-x-k)_i(\gamma(x,k)-x+k)_j\partial_{i}\partial_jf(z)\)
\right|\to0\qquad\quad
\eq
uniformly on $\FR$ as $N\to\infty$ which is a direct consequence of proposition~\ref{prop:uniform-result}.
\end{proof}


\ssection{Convergence of the Sequence of Stochastic Processes}

In this section we use a relative compactness argument to deduce the existence of a sequence of measures that converge weakly to a weak solution to \Eq.model. in the case of the Linear Growth Condition or a vague solution to \Eq.model. if the more general growth conditions are satisfied. The strategy is to consider spaces of discrete measures with bounded total mass and investigate the two component process $(P^N,u^N)$ where $u^N$ is an approximation to the time. By using $u^N$ instead of the true time we may utilize a stochastic generator that has no explicit time dependence and it is moderately easy to write down a set of rules that solve the associated martingale problem. We proceed to show that the sequence of stochastic processes is tight and Prohorov's theorem reveals the existence of an accumulation point. We verify that this accumulation point is a solution to \Eq.model. in the relevant sense. Once this is done we prove a result that allows us to dispense with the fictitious time variable.

It is advantageous to introduce an extra process associated with the fictitious time variable, the purpose of which is to bound the total rate of all processes away from zero. This {\em clock process} does not alter the current state of the measure describing the stochastic particle ensemble, it merely advances the fictitious time variable. In this way we can be sure that if all the other rates momentarily become zero the computed holding time in the current state remains finite.

Throughout this section we assume that $\aN $ and $\cN $ are sequences of positive real numbers that tend to infinity as $N\to\infty$.

\begin{definition}
\it Define spaces of measures by  
\be
\BS_C=\{P\in\CM_B(\FR): P(\pim)\le C\},\and\BS^N=\left\{P=\frc1N\sum_{i=1}^n\delta_{z_i}:z_i\in \FR,\ P(\pim)\le \aN \right\}
\ee
where $\delta_z$ is the Dirac measure concentrated at $z$. 
It is convenient to set $E^N=\BS^N\times[0,\infty)$.
\end{definition}

\begin{definition}
We put the metric $D(\cdot,\cdot)$ on $E=\CM_B(\FR)\times[0,\infty)$ and its subspaces in terms of \Eq.metric0. by defining
\be
D((P,u),(Q,v))=d_{\rm weak}(P,Q)+|u-v|
\ee
and let $\tilde D$ be the Skorohod metric on the space of c\`adl\`ag functions $D_E[0,\infty)$ derived from $D$ where 
if 
$\rho$ is a metric on a space $X$  then the Skorohod metric $\tilde\rho$ on $D_X[0,\infty)$
is defined by (cf.  equation~(3.5.2) in \cite{EK})
\be
\tilde \rho(x,y)=\inf_{\lambda\in\Lambda}\left[\gamma(\lambda)\vee\int_0^\infty\e^{-T}\sup_{t\ge0}\rho(x_{t\wedge T},y_{\lambda(t)\wedge T})\,\d T\right];\qquad
\gamma(\lambda)=\esssup_{0\le t<s}\left|\log\frc{\lambda(s)-\lambda(t)}{s-t}\right|
\LABEL{eq:skorohod}
\ee
where
\be
\Lambda=\{\lambda\in C([0,\infty),[0,\infty)):
\mbox{$\lambda$ is a strictly increasing Lipschitz function with $\gamma(\lambda)<\infty$}\}.
\ee
\end{definition}

\begin{remark}\em
We adopt the notational conventions for the mass-weighted source moments:
\be
\Lambda_t\equiv I_t(\FR)\and\Lambda^{(q)}=\sup\{I_t(\pim^{(q)}):0\le t<\infty\}.
\ee
\end{remark}

\begin{definition}\LABEL{def:selectionmeasure}
Let
$\nu^N_n:\CM_B(\FR)\times\FB(\FR)^n\to\BR$ be defined by $\nu^N_n(P,B_1\times\cdots\times B_n)=N^{-n}\nu_n(NP,B_1\times\cdots\times B_n)$ where
\be
\nu_n(P,B_1\times\cdots\times B_n)=\sum_{\sigma\in S_n}\zeta(\sigma)\kern-10pt\prod_{
\begin{array}{c}\scriptstyle
{\rm disjoint\ cycles}\\
\scriptstyle
{(i_1,\ldots,i_k){\rm\ of\ }\sigma}
\end{array}
}\kern-10pt
P(B_{i_1}\cap\ldots\cap B_{i_k})\LABEL{eq:nu_n}.
\ee
Here $S_n$ is the symmetric group on $\{1,\ldots,n\}$
and $\zeta$ is the signature homomorphism.
\end{definition}

\begin{remark}\em\LABEL{rem:distinct}
By proposition~\ref{prop:particleselectionmeasures}
\be
\nu_i^N\(\frc1N\sum_{k=1}^n\delta_{w_k},\{z_1\}\times\cdots\times\{z_i\}\)=\begin{cases}N^{-i}&{\vcenter{\hbox{
if there exists an injection $\alpha:\{1,\ldots,i\}\to\{1,\ldots,n\}$}
\hbox{
such that $z_r=w_{\alpha(r)}$ for $r=1,\ldots,i$;}}
}
\\
0&\mbox{otherwise.}
\end{cases}
\ee
\end{remark}

\begin{definition}
Define the sequence of generators $\CA^N:M(E^N)\to M(E^N)$ to be
\beq
\CA^N(F)(P,u)&=&{Nr}\{F(J_{\rm CLK}^N(P,u))-F(P,u)\}\cr
&&\kern0cm{}+N\int_\FR\{F(J_{I}^N((P,u),z))-F(P,u)\}I_u(\d z)\cr
&&\kern0cm{}+{N\cN ^2}\int_\FR\BE^Z\left\{F\(J_0^N\((P,u),z,\frc{\sigma_u(z)Z}{\cN }+\frc{b_u(z)}{\cN ^2}\)\)-F(P,u)\right\}P(\d z)\cr
&&\kern0cm{}+N\sum_{i=1}^\CI\sum_{j=1}^\CJ\frc 1{i!\,j!}\int_{\FR^{i+j}}K_{u,i,j}(z_1,\ldots,z_i,\d w_1,\ldots,\d w_{j})\cr
&&\kern1cm{}\times\{F(J^N_{i,j}((P,u),z_1,\ldots,z_i,w_1,\ldots, w_{j}))-F(P,u)\}\nu_i^N(P,\d z_1,\ldots,\d z_i)\qquad
\LABEL{eq:generatorA_N}
\eq
where $Z\sim N(0,\Id1)$, and the jumps: $J^N_{\rm CLK}:E^N\to E$, $J^N_I:E^N\times\FR\to E$,  $J_0^N:E^N\times\FR\times\BR^{d_1}\to E$, and $J^N_{i,j}:E^N\times\FR^i\times\FR^j\to E$ for $i=1,\ldots,\CI$ and $j=1,\ldots,\CJ$
are given by
\beq
J_{\rm CLK}^N(P,u)&=&\(P,u+\frc1{\rho^N(P,u)}\);\LABEL{eq:jumpCLK}\\
J_{I}^N((P,u),z)&=&\(j^N_I(P,z),u+\frc1{\rho^N(P,u)}\);\LABEL{eq:jumpI}\\
j^N_I(P,z)&=&\begin{cases}P+\frc1N\delta_z&\mbox{if
$P(\pim)+\frc{\pim(z)}N\le \aN $}
\\P&\mbox{otherwise;}\end{cases}\\
J_0^N((P,u),z,k)&=&\(P+\frc1N(\delta_{(\pim(z),\gamma(\pi_{\bar\Omega}(z),k),j_{\rm ID}^N(u,z))}-\delta_z),u+\frc1{\rho^N(P,u)}\);\qquad\LABEL{eq:jump0}\\
J_{i,j}^N((P,u),z_1,\ldots,z_i,w_1,\ldots, w_j)&=&\(P+\frc1N\(\sum_{m=1}^j\delta_{w_m}-\sum_{n=1}^i\delta_{z_n}\),u+\frc1{\rho^N(P,u)}\)\cr
&&\kern3cm\qquad\mbox{for $i=1,\ldots,\CI$, $j=1,\ldots,\CJ$};\LABEL{eq:jumpij} 
\eq
and
\be
\rho^N(P,u)=N\(r+\Lambda_u+
\cN ^2P(\FR)\vphantom{\sum_{j=1}^\CJ\int}
+\sum_{i=1}^\CI\sum_{j=1}^\CJ\frc1{i!\,j!}\int_{\FR^i}K_{u,i,j}(z_1,\ldots,z_i,\FR^j)\nu^N_i(P,\d z_1,\ldots,\d z_i)\right).
\LABEL{eq:rho^N}
\ee
Notice that the image of each of the jumps is a member of $E^N$ whenever $P(\{z\})\neq0$ in \Eq.jump0. and
$\nu^N_i(P,\{z_1\}\times\cdots\times\{z_i\})\neq0$ with $\sum_{m=1}^j\pim(w_m)=\sum_{n=1}^i\pim(z_n)$ in \Eq.jumpij.. This provides the motivation for definition~\ref{def:selectionmeasure}
\end{definition}

\begin{remark}\em
The generators $\CA^N$ define a standard Markov jump process:
\be
\CA^N(F)(P,u)=\rho^N(P,u)\int_{E^N}\{F(X)-F(P,u)\}\mu^N((P,u),\d X)
\ee
with transition function $\mu^N:E^N\times\FB(E^N)\to[0,1]$ defined by
\beq
\mu^N((P,u),\cdot)&=&\frc {rN}{\rho^N(P,u)}
\delta_{J_{\rm CLK}^N(P,u)}\cr
&&{}+\frc N{\rho^N(P,u)}
\int_\FR I_u(\d z)\delta_{J_I^N((P,u),z)}\cr
&&\kern0cm{}+\frc {N\cN ^2}{\rho^N(P,u)}
\int_\FR P(\d z)\BE^Z\(\delta_{J_0^N((P,u),z,\sigma_u(z)Z\cN ^{-1}+b_u(z)\cN ^{-2})}\)\cr
&&{}+\frc N{\rho^N(P,u)}\sum_{i=1}^\CI
\sum_{j=1}^\CJ\frc1{i!\,j!}\int_{\FR^{i+j}} K_{u,i,j}(z_1,\ldots,z_i,\d w_1,\ldots,\d w_j)\cr
&&\kern3cm\times\nu^N_i(P,\d z_1,\ldots,\d z_i)\delta_{J^N_{i,j}((P,u),z_1,\ldots,z_i,w_1,\ldots,w_j)}.
\eq
 Note that
\be
\rho^N(P,u)\le N(r+\Lambda^{(0)}+\cN ^2\aN +\CI \CJ K_\infty  \aN ^{\CI}+K^\infty_{\aN }).
\ee
In consequence the generators are bounded for each $N$ and according to section~4.2 \cite{EK} there exists a sequence of $E^N$-valued random variables $(P^N_t,u^N_t)$ such that
\be
M^N_t(F)=F(P^N_t,u^N_t)-F(P^N_0,u^N_0)-\int_0^t\CA^N(F)(P^N_s,u^N_s)\,\d s\LABEL{eq:martingales}
\ee
is a martingale for each $F\in M(E^N)$ for any random initial condition $(P_0^N,u^N_0)\in E^N$.
\end{remark}

\begin{assumption}
\it Let $u^N_0=0$ for all $N$, and suppose that there exists a sequence $P^N_0\in\BS^N$ with $P^N_0\weakto P_0$ as $N\to\infty$ such that $P^N_0$ is deterministic
with $P^N_0(\pim)$  bounded above. For notational convenience we define
\be
\Xi=\sup_NP_0^N(\pim).
\ee
\end{assumption}

\begin{remark}\em\LABEL{rem:martingales}
By considering the bounded measurable functions $\pi_u^L(P,u)=u\wedge L$ and $\Phi_f(P,u)=P(f)$ where $f\in\FF$
we find that
\be
M^N_t(\pi_u^L)=u_t^N\wedge L-\int_0^t\rho^N(P^N_s,u^N_s)\left\{\(u^N_s+\frc1{\rho^N(P^N_s,u^N_s)}\)\wedge L-u^N_s\wedge L\right\}\d s\LABEL{eq:mart1}
\ee
and
\beq
M^N_t(\Phi_f)&=&P^N_t(f)-P^N_0(f)-\int_0^t\int_\FR f(z)\openone_{P_s^N(\pim)+\pim(z)/N\le \aN }\,I_{u^N_s}(\d z)\,\d s\cr
&&\kern-1cm{}-\cN ^2\int_0^t\int_\FR\BE^Z\{f(\pim(z),\gamma(\pi_{\bar\Omega}(z),\sigma_u(z)Z\cN ^{-1}+b_u(z)\cN ^{-2}),j_{\rm ID}^N(u_s^N,z))-f(z)\}
\,P_s^N(\d z)\,\d s\cr
&&\kern-1cm{}-\sum_{i=1}^\CI\sum_{j=1}^\CJ\frc1{i!\,j!}\int_0^t\int_{\FR^{i+j}}K_{u^N_s,i,j}(z_1,\ldots,z_i,\d w_1,\ldots,\d w_j)\left\{\sum_{m=1}^jf(w_m)-\sum_{n=1}^if(z_n)\right\}\cr
&&\kern8cm{}\times\nu^N_i(P^N_s,\d z_1,\ldots,\d z_i)\,\d s
\LABEL{eq:mart2}
\eq
where $Z\sim N(0,\Id1)$, are martingales for each $N$.
\end{remark}

\begin{proposition}\LABEL{prop:pi_u}
 Let $\(\CF^N_t\)_{t\ge0}$ be a filtration adapted to $(P^N,u^N)$. Suppose that $\pi_u(P,u)=u$ and $\pi_u^{(2)}(P,u)=u^2$ then
\be
M^N_t(\pi_u)=u^N_t-t\and M^N_t(\pi_u^{(2)})
\ee
are $\CF^N_t$-martingales for each $N$.
\end{proposition}
\begin{proof}
We have
\beq
\BE(u^N_t\wedge L)&=&t+\int_0^t\BE(\CA^N(\pi_u^L)(P^N_\tau,u^N_\tau)-1)\,\d\tau\\
&=&t+\int_0^t\BE\(\rho^N(P^N_\tau,u^N_\tau)\left\{\(u^N_\tau+\frc1{\rho^N(P^N_\tau,u^N_\tau)}\)\wedge L\right.\right.\cr
&&\kern5cm\left.\left.{}-u^N_\tau\wedge L-\frc1{\rho^N(P^N_\tau,u^N_\tau)}\right\}\)\,\d\tau\le t.
\eq
Hence by the monotone convergence theorem $u^N_t$ is of class $L^1$ and $\BE(u^N_t)\le t$.
Furthermore if $t\ge s$,
\beq
\BE|\BE\(u^N_t\wedge L-t\mid\CF^N_s\)-(u^N_s\wedge L-s)|&=&\BE\left|\int_s^t\BE\(\rho^N(P^N_\tau,u^N_\tau)\left\{\(u^N_\tau+\frc1{\rho^N(P^N_\tau,u^N_\tau)}\)\wedge L\right.\right.\right.\cr
&&\kern1cm\left.\left.\left.{}-u^N_\tau\wedge L-\frc1{\rho^N(P^N_\tau,u^N_\tau)}\right\}\mid\CF^N_s\)\,\d\tau\right|\\
&\le&\BE\left\{(t-s)\BP\(u^N_t>L-\frc 1{rN}\mid\CF^N_s\)\right\}\\
&\le&\frc{(t-s)\BE(\BE(u^N_t|\CF^N_s))}{L-N^{-1}r^{-1}}\le\frc{t(t-s)}{L-N^{-1}r^{-1}}\to0
\eq
as $L\to\infty$, the conditional form of the monotone convergence theorem and Fatou's lemma then imply that $u_t^N-t$ is an $\CF^N_t$-martingale. Similarly
\beq
\BE((u^N_t)^{2}\wedge L^2)&=&t^2+\int_0^t\BE\(\rho^N(P^N_\tau,u^N_\tau)\left\{\(u^N_\tau+\frc1{\rho^N(P^N_\tau,u^N_\tau)}\)^2\wedge L^2\right.\right.\cr
&&\kern5cm\left.\left.{}-(u^N_\tau)^2\wedge L^2-\frc{2u^N_\tau}{\rho^N(P^N_\tau,u^N_\tau)}\right\}\)\,\d\tau\qquad\qquad\\
&\le&t^2+tr^{-1}N^{-1}.
\eq
Thus the monotone convergence theorem gives that $\BE((u^N_t)^2)<\infty$.
Moreover
\beq
&&\kern-1cm\BE\left|\BE\((u^N_t)^2\wedge L^2-\int_s^t\rho^N(P^N_\tau,u^N_\tau)\left\{\(u^N_\tau+\frc1{\rho^N(P^N_\tau,u^N_\tau)}\)^2-(u^N_\tau)^2
\right\}
\,\d\tau\mid\CF^N_s\)-(u^N_s)^2\wedge L^2\right|\cr
&&\kern0cm{}\le(t-s)\BE\left\{\BE\(\{2u^N_t+N^{-1}r^{-1}\}\openone_{u^N_t\le L-\frac1{rN}}\mid\CF^N_s\)\right\}
\\
&&\kern3cm{}=(t-s)\BE\(\{2u^N_t+N^{-1}r^{-1}\}\openone_{u^N_t\le L-\frac1{rN}}\)
\to0
\eq
as $L\to\infty$ by the dominated convergence theorem. Applying the conditional form of the monotone convergence theorem and Fatou's lemma on the LHS we deduce that
\be
\BE\((u^N_t)^2-\int_s^t\rho^N(P^N_\tau,u^N_\tau)\left\{\(u^N_\tau+\frc1{\rho^N(P^N_\tau,u^N_\tau)}\)^2-(u^N_\tau)^2
\right\}
\,\d\tau\mid\CF^N_s\)=(u^N_s)^2\qquad{\mbox{a.s.}}
\ee
and thus $M^N_t(\pi_u^{(2)})$ is an $\CF_t^N$-martingale.
\end{proof}

\begin{proposition}\LABEL{prop:Qpoly}
For all $N$ and $t\ge0$,
\be
\BE((P^N_t(\pim))^i)\le Q_i(t)\qquad i=0,\ldots,\CI,
\ee
where $Q_i$ is an increasing polynomial on $[0,\infty)$ of degree less than or equal to $i$.
\end{proposition}
\begin{proof}
We prove the result by induction. If $i=0$, then we can take $Q_0(t)=1$. For $0<i\le\CI$, we assume the result holds for all $r<i$, set $F((P,u))=(P(\pim))^i\le \aN ^i$ in \Eq.martingales. and take expectations to deduce that
\beq
\BE((P^N_t(\pim))^i)&=&(P^N_0(\pim))^i+N\int_0^t\BE\(\left[P^N_s(\pim)+\frc{\pim(z)}N\right]^i-[P^N_s(\pim)]^i\)\cr
&&\kern4cm{}\times\openone_{P_s^N(\pim
)+\pim(z)/N\le \aN }
\,I_{u^N_s}(\d z)
\,\d s\qquad\\
&\le&\Xi^i+\sum_{r=0}^{i-1}\frc{i!\,\Lambda^{(i-r)}}{(i-r)!\,r!}\int_0^tQ_r(s)\,\d s\equiv Q_i(t),
\eq
where we have made use of hypothesis~\ref{hyp:cont}.
\end{proof}

\begin{proposition}\LABEL{prop:two9}
For all $T>0$, $f\in\FF_c\cup\{f_0\}$ there exist
an increasing polynomial $\CQ$ on $[0,\infty)$ of degree no more than $\CI+1$ such that
\beq
\mbox{(a)}\kern2cm
\BE\sup_{0\le t\le T}|M^N_t(\Phi_f)|&\le&\(\frc
{\CQ(T)}N
\)^{\frac12}K^\infty_{R(f)}H^\infty_{R(f)}\:\fext\:^\#;\kern3.8cm\qquad\LABEL{eq:exp21}\\
\mbox{(b)}\kern1.8cm
\BE\sup_{0\le t\le T}|M^N_t(\Phi_{\pim})|&\le& 2\(\frc{T\Lambda^{(2)}}N\)^{\frac12}\LABEL{eq:sup2}
\eq
(c) for all $\epsilon>0$ there exists $R_T(\epsilon)>0$ such that
\be
\BP\(\sup_{0\le t\le T} |u^N_t-t|\ge R_T(\epsilon)\)\le\epsilon\LABEL{eq:exp20}.
\ee
(d) In particular
\be
\BP\(\sup_{0\le t\le T+1}|u^N_t-t|\ge N^{-\frac14}\)\le2(T+1)^{\frac12}r^{-\frac12}N^{-\frac14}\to0\as N->\infty.\LABEL{eq:uNt}
\ee
\end{proposition}
\begin{proof}
The previsible increasing processes associated with the martingales equations~(\ref{eq:mart1}-\ref{eq:mart2}) are given by
\beq
\<M^N(F)\>_t&=&\int_0^t\rho^N(P^N_s,u^N_s)\int_{E^N}\{F(X)-F(P^N_s,u^N_s)\}^2\mu^N((P^N_s,u^N_s),\d X)\,\d s
\\
&&\kern-1.5cm{}=rN\int_0^t\{F(J_{\rm CLK}^N(P^N_s,u^N_s))-F(P_s^N,u^N_s)\}^2\,\d s\cr
&&\kern-1cm{}+N\cN ^2\int_0^t\int_\FR\BE^Z\{F(J_0^N((P^N_s,u^N_s),z,\sigma_{u^N_s}(z)Z\cN ^{-1}+b_{u^N_s}(z)\cN ^{-2}))-F(P^N_s,u^N_s)\}^2 P^N_s(\d z)\,\d s\cr
&&\kern-1cm{}+N\int_0^t \int_{\FR}\{F(J_I^N((P^N_s,u^N_s),z))-F(P_s^N,u^N_s)\}^2\,I_{u^N_s}(\d z)\,\d s
\qquad\qquad\cr
&&\kern-1cm{}+N\sum_{i=1}^\CI\sum_{j=1}^\CJ\frc1{i!\,j!}\int_0^t\int_{\FR^{i+j}}K_{u^N_s,i,j}(z_1,\ldots,z_i,\d w_1,\ldots,\d w_j)\cr
&&\kern-0.6cm{}\times
\{F(J^N_{i,j}((P^N_s,u^N_s),z_1,\ldots,z_i,w_1,\ldots,w_j))-F(P^N_s,u^N_s)\}^2\nu^N_i(P^N_s,\d z_1,\ldots,\d z_i)\,\d s\cr
&&\LABEL{eq:increasing}
\eq
leading to
\beq
\<M^N(\Phi_f)\>_t&\le&\frc1N\int_0^t\int_\FR f(z)^2I_{u^N_s}(\d z)\,\d s\cr
&&{}+\frc{\cN ^2}N\int_0^t
\int_\FR\BE^Z\{f(\pim(z),\gamma(\pi_{\bar\Omega}(z),\sigma_{u^N_s}(z)Z\cN ^{-1}+b_{u^N_s}(z)\cN ^{-2}),j_{\rm ID}^N(u_s^N,z))\cr
&&\kern9.5cm{}-f(z)\}^2P^N_s(\d z)\,\d s \cr
&&+\frc1{N}\sum_{i=1}^\CI\sum_{j=1}^\CJ\frc1{i!\,j!}\int_0^t
\int_{\FR^{i+j}}K_{u^N_s,i,j}(z_1,\ldots,z_i,\d w_1,\ldots,\d w_j)\cr
&&\kern3cm{}\times\left\{\sum_{m=1}^jf(w_m)-\sum_{n=1}^if(z_n)\right\}^2
\nu^N_i(P^N_s,\d z_1,\ldots,\d z_i)\,\d s.
\LABEL{eq:mart3}
\eq
Write $k=\sigma_{u^N_s}(z)Z\cN ^{-1}+b_{u^N_s}(z)\cN ^{-2}$ and $x=\pi_{\bar\Omega}(z)$. The inequality
\beq
&&\cN ^2\BE^Z([f(\pim(z),\gamma(x,k),j_{\rm ID}^N(u_s^N,z))-f(z)]^2)\kern8cm\ \cr
&&\kern2cm\le2\cN ^2\BE^Z([f(\pim(z),\gamma(x,k),j_{\rm ID}^N(u_s^N,z))-f(\pim(z),x,j_{\rm ID}^N(u_s^N,z))]^2)\cr
&&\kern7cm{}+2\cN ^2[f(\pim(z),x,j_{\rm ID}^N(u_s^N,z))-f(z)]^2\\
&&\kern2cm\le2\cN ^2\BE^Z(\|k\|^2\|\D \fext\|_\infty^{\#\,2})+4\|f\|_\infty H^\infty_{R(f)}\|\hat\D f\|_\infty\pim(z)\\
&&\kern2cm{}\le4\left[\(\sigma_\infty^2d_1+b_\infty^2\sup_N\cN ^{-2}\)\|\D \fext\|_\infty^{\#\,2}+\|f\|_\infty H^\infty_{R(f)}\|\hat\D f\|_\infty\right]\pim(z)
\eq
then shows that
\be
\<M^N(\Phi_f)\>_t\le\frc{\CS(t)}{4N}
\ee
where
\beq
\CS(t)&=&\left\{4\Lambda^{(0)} t+
\left[16\sigma_\infty^2d_1+16b_\infty^2\sup_N\cN ^{-2}+ 16\right]\int_0^tP^N_s(\pim)\,\d s\right.\cr
&&\kern0cm\left.{}+4\int_0^tP^N_s(\pim)\,\d s+4\CI\CJ(\CI+\CJ)^2K_\infty\int_0^t\(P^N_s(\pim)\)^\CI\,\d s\right\}K^\infty_{R(f)}H^\infty_{R(f)}\:\fext\:^{\#\,2}.\qquad
\eq
It follows by proposition~\ref{prop:Qpoly} that if
\be
\CQ(t)=4\left[4\sigma_\infty^2d_1+4b_\infty^2\sup_N\cN ^{-2}+5\right]\int_0^tQ_1(s)\,\d s+4\Lambda^{(0)} t+4\CI\CJ(\CI+\CJ)^2K_\infty \int_0^tQ_\CI(s)\,\d s
\ee
then $\BE\<M^N(\Phi_f)\>_T\le \frcn14\CQ(T)N^{-1}K^\infty_{R(f)}H^\infty_{R(f)}\:\fext\:^{\#\,2}$.
Doob's $L^2$-Inequality
\be
\(\BE\sup_{0\le t\le T}|M^N_t(F)|\)^2\le4\BE\<M^N(F)\>_T\LABEL{eq:Doob}
\ee
gives part (a). Our definitions $R(f_0)=0$ and $H^\infty_0=K^\infty_0=1$ are essential to check the $n=0$ case.
Furthermore \Ineq.mart3. and the mass preserving property of the interaction kernels given by hypothesis~\ref{hyp:interactionker} give $\<M^N(\Phi_{\pim})\>_t\le t\Lambda^{(2)}/N$, and thus \Ineq.Doob. leads to part (b).

By proposition~\ref{prop:pi_u}, $M^N_t(\pi_u)$ and $M^N_t(\pi_u^{(2)})$ are martingales. In consequence the increasing process associated with  $M^N_t(\pi_u)$ is given by \Eq.increasing. with $F=\pi_u$, thus
\be
\<u^N-t\>_t\le\frc{t}{rN}\le\frc t{r}\LABEL{eq:Esup2}
\ee
and hence
\be
\BE\sup_{0\le s\le t}|u^N_s-s|\le2\(\frc t{rN}\)^{\frac12}.\LABEL{eq:expuNs}
\ee
Inequality~(\ref{eq:exp20}) follows from \Ineq.Doob. applied to $u_t^N-t$ and the Markov inequality, with $R_T(\epsilon)=2T^{\frac12}r^{-\frac12}N^{-\frac12}\epsilon^{-1}$ establishing part (c). If we set $\epsilon=2(T+1)^{\frac12}r^{-\frac12}N^{-\frac14}$ then part (d) follows.
\end{proof}

\begin{proposition}\LABEL{prop:Qhat}
There exist an increasing polynomial $\hat\CQ$ on $[0,\infty)$ of degree no more than $\CI$ such that (a)
\be
\BE|\CA^N(\Phi_f)(P_t^N,u_t^N)|\le \hat\CQ(t)K^\infty_{R(f)}H^\infty_{R(f)}\:\fext\:^\#
\LABEL{eq:expA}
\ee
for all $f\in\FF_c\cup\{f_0\}$ and $t\ge0$ and (b)
\be
|\CA^N(\Phi_{\pim})(P^N_t,u^N_t)|\le\Lambda^{(1)}.\LABEL{eq:exp2}
\ee
for all $t\ge0$.
\end{proposition}
\begin{proof}
(a) Setting $x=\pi_{\bar\Omega}(z)$, $k=\sigma_{u^N_t}(z)Z\cN ^{-1}+b_{u^N_t}(z)\cN ^{-2}$ and $Z\sim N(0,\Id1)$ we have
\beq
|\CA^N(\Phi_f)(P^N_t,u^N_t)|&\le& I_{u^N_t}(|f|)+\cN ^2\int_\FR|\BE^Z\(f(\pim(z),\gamma(x,k),j_{\rm ID}^N(u^N_t,z))-f(z)\)| P^N_t(\d z)\cr
&&{}+\sum_{i=1}^\CI\sum_{j=1}^\CJ\frc1{i!\,j!}\int_{\FR^{i+j}}K_{u^N_t,i,j}(z_1,\ldots,z_i,\d w_1,\ldots,\d w_j)
\cr
&&\kern2cm{}\times\left|\sum_{m=1}^jf(w_m)-\sum_{n=1}^if(z_n)\right|\nu^N_i(P^N_t,\d z_1,\ldots,\d z_i)
\LABEL{eq:exp3}
\\
&\le&\Lambda^{(0)}\|f\|_\infty+LH^\infty_{R(f)}\:\fext\:^\#P^N_t(\pim)\cr
&&\kern1cm{}+K^\infty_{R(f)}P_t^N(\pim)\|f\|_\infty+\CI\CJ(\CI+\CJ)K_\infty P^N_t(\pim)^\CI\|f\|_\infty
\eq
where we have used proposition~\ref{prop:boundedness} with $L$ given by \Eq.L_1.. Accordingly \Eq.expA. follows with
\be
\hat\CQ(t)=\Lambda^{(0)}+(L+1)Q_1(t)+\CI\CJ(\CI+\CJ)K_\infty Q_\CI(t)
\ee
where we have made use of proposition~\ref{prop:Qpoly}. The $n=0$ case is dealt with by using the definitions $H^\infty_0=K^\infty_0=1$.

\ni(b) Inequality~(\ref{eq:exp2}) follows directly from \Ineq.exp3. together with the mass conserving property of the interaction kernel (hypothesis~\ref{hyp:interactionker}).
\end{proof}

\begin{proposition}\LABEL{prop:two10}\it
\it For all $T>0$ and $\epsilon>0$ there exists $L_T(\epsilon)$ such that for all $N\ge1$ we have that
\be
\BP\(\sup_{t\le T}P^N_t(\pim)>L_T(\epsilon)\) <\epsilon.\LABEL{eq:LT}
\ee
\end{proposition}
\begin{proof}
If $t\le T$ we have from the martingale definition
\be
P^N_t(\pim)\le P_0^N(\pim)+\int_0^T|\CA^N(\Phi_{\pim})(P^N_t,u^N_t)|\,\d t+|M^N_t(\Phi_{\pim})|\le
\Xi +T\Lambda^{(1)}+\sup_{t\le T}|M^N_t(\Phi_{\pim})|.
\ee
Take the supremum and expectation and apply proposition~\ref{prop:two9}(b) to find
\be
\BE\(\sup_{t\le T}P^N_t(\pim)\)\le K_T
\ee
where
\be
K_T=\Xi +T\Lambda^{(1)}+2(T\Lambda^{(2)})^{\frac12}
\ee
and thus, if $L_T(\epsilon)=K_T/\epsilon$ then Markov's inequality gives \Eq.LT..
\end{proof}

\begin{corollary}
For all $T>0$ and $0<\eta\le1$,
\be
\inf_N\BP\((P^N_t,u^N_t)\in\BS_{L_T(\eta/2)}\times[0,T+R_T(\eta/2)],\,\ 0\le t\le T\)\ge 1-\eta.
\ee
\end{corollary}
\begin{proof}
This is immediate from propositions~\ref{prop:two9} and~\ref{prop:two10}.
\end{proof}

\begin{proposition}
For all $C>0$, the space $\BS_C$ is a compact subset of $\CM_B(\FR)$.
\end{proposition}
\begin{proof}  
First we show that $\BS_C$ is closed. Suppose $p_n\in\BS_C$ with $p_n\weakto p$, then for each $L\in\BN$ set $\pim^L(z)=L\wedge\pim(z)$, we have that $C\ge p_n(\pim^L)\to p(\pim^L)$. Thus
$p(\pim^L)\le C$ and by taking $L\to\infty$, the monotone convergence theorem implies that $p(\pim)\le C$, i.e., $p\in \BS_C$ since $p(\FR)\le p(\pim)\le C$ and so $p\in\CM_B(\FR)$.

Next we show sequential compactness. Suppose $p_n\in\BS_C$, write $p_n=\lambda_nq_n$ with $q_n(\FR)=1$ and $\lambda_n=p_n(\FR)$. Then
$\lambda_n\in[0,C]$ and
therefore there exists a convergent subsequence $\lambda_{n_k}\to\lambda$. If $\lambda=0$, then for any bounded continuous function $f$ on $\FR$ we have
$p_{n_k}(f)\le\lambda_{n_k}\|f\|_\infty\to0$ as $k\to\infty$, i.e., $p_{n_k}\weakto0$. If $\lambda\neq0$ then we may assume, by passing to a further subsequence if necessary, that $\lambda_{n_k}>\h\lambda$ for all $k$ and we find that for any $\eta>0$ we have
\be
q_{n_k}\(\pim^{-1}\{1,\ldots, \lceil 2C/\lambda\eta\rceil\}\)=1-q_{n_k}\(\pim^{-1}\{\lceil 2C/\lambda\eta\rceil+1,\ldots\}\)\ge1-\frc{\lambda_{n_k}}\eta Cq_{n_k}(\pim)\ge1-\eta
\ee
where $\lceil x\rceil=\min\{n\in\BZ:n\ge x\}$.
Note that since $\pim$ is continuous $\pim^{-1}\{1,\ldots, \lceil 2C/\lambda\eta\rceil\}$ is  closed and bounded, and thus compact. Accordingly the $q_{n_k}$ form a tight sequence
of probability measures.
Prohorov's theorem (theorem~3.2.2, \cite{EK}) implies that there exists a further subsequence $q_{n_{k(r)}}$ which converges weakly to a probability measure $q$. Hence
\be
d_{\rm weak}(p_{n_{k(r)}},\lambda q)=\sum_{\ell=0}^\infty\frc1{4^\ell}(|\lambda_{n_{k(r)}}q_{n_{k(r)}}(f_\ell)-\lambda q(f_\ell)|\wedge1)\to0\as r->\infty
\ee
with limiting measure $\lambda q\in\BS_C$ by closure. It follows that $\BS_C$ is sequentially compact.
\end{proof}

\begin{corollary}\LABEL{cor:compact}
The sequence $(P^N_t,u^N_t)$ obeys the compact support condition (remark~3.7.3, \cite{EK}).
\end{corollary}

\begin{definition} The modulus of continuity of an element $x\in D_E[0,\infty)$ is defined to be
\be
w'(x,\delta,T)=\inf_{t_i}\ \max_i\sup_{s,t\in[t_{i-1},t_i)}D(x_s,x_t)
\ee
where $0=t_0\le t_1\le\cdots\le t_{n-1}\le T\le t_n$ and $\min\{t_{i}-t_{i-1}\}>\delta$.
\end{definition}

\begin{proposition}\LABEL{prop:modulus}\it
For any $T>0$ and $\eta>0$,
\be
\limsup_N\BP(w'((P^N,u^N),\delta,T)\ge\eta)\le\eta
\ee
for sufficiently small $\delta$.
\end{proposition}
\begin{proof}
We investigate the partition given by $t_j=2j\delta$, then for $x\in D_E[0,\infty)$,
\be
w'(x,\delta,T)\le2\max_j\sup_{2j\delta\le s<2(j+1)\delta}D(x_{2j\delta},x_s).
\ee
It follows that
\beq
\BP(w'((P^N,u^N),\delta,T)\ge\eta)&\le&\BP\(\max_j\sup_{2j\delta\le s<2(j+1)\delta}D((P^N_{2j\delta},u^N_{2j\delta}),(P^N_s,u^N_s))\ge\h\eta\)\qquad\\
&\le&\frc2\eta\BE\(\max_j\sup_{2j\delta\le s<2(j+1)\delta}D((P^N_{2j\delta},u^N_{2j\delta}),(P^N_s,u^N_s))\).
\eq
If $2j\delta\le T$ and $2j\delta\le s\le2(j+1)\delta$ and we write $\Phi_{f_n}(P,u)=P(f_n)$ we have
\beq
D((P^N_{2j\delta},u^N_{2j\delta}),(P^N_s,u^N_s))&\le&
\sum_{n=0}^\infty\frc1{4^n}(|P^N_s(f_n)-P^N_{2j\delta}(f_n)|\wedge1)+|u^N_s-u^N_{2j\delta}|\\
&\le&\sum_{n=0}^\infty\frc1{4^n}\(\left|M^N_s(\Phi_{f_n})-M^N_{2j\delta}(\Phi_{f_n})+\int_{2j\delta}^s\CA^N(\Phi_{f_n})(P^N_\tau,u^N_\tau)\,\d \tau\right|\wedge1\)\cr
&&\kern5cm{}+|u^N_s-u^N_{2j\delta}|\\
&\le&\sum_{n=0}^\infty\frc1{4^n}\(\left|2\sup_{0\le\tau\le T+2\delta}|M^N_\tau(\Phi_{f_n})|+\int_{2j\delta}^s|\CA^N(\Phi_{f_n})(P^N_\tau,u^N_\tau)|\,\d\tau\right|\wedge1\)\cr
&&\kern3cm{}+2\sup_{0\le\tau\le T+2\delta}|u^N_\tau-\tau|+{2\delta}.
\eq
Thus propositions~\ref{prop:two9} and~\ref{prop:Qhat} give
\beq
\BE\(\max_j\sup_{2j\delta\le
s<2(j+1)\delta}D((P^N_{2j\delta},u^N_{2j\delta}),(P^N_s,u^N_s))\)&&\cr
&&\kern-7cm\le
\sum_{n=0}^\infty\frc1{4^n}\(\left|2\BE\sup_{0\le\tau\le T+2\delta}|M^N_\tau(\Phi_{f_n})|+2\delta
\hat\CQ(T+2\delta)K^\infty_{R(f_n)}H^\infty_{R(f_n)}\:\fext_n\:^\#\right|\wedge1\)\cr
\qquad{}+2\BE\sup_{0\le\tau\le T+2\delta}|u^N_\tau-\tau|+{2\delta}&&\\
&&\kern-7cm\le
\sum_{n=0}^\infty\frc1{4^n}\(\left| 2K^\infty_{R(f_n)}H^\infty_{R(f_n)}\:\fext_n\:^\#\left[\(\frc{
\CQ(T+2\delta)}N\)^{\frac12}\right.\right.\right.\cr
&&\kern-2cm\left.\left.\left.\phantom{\frc\|\|}{}+\delta
\hat\CQ(T+2\delta)\right]\right|\wedge1\)+4\(\frc{T+2\delta}{Nr}\)^{\frac12}+{2\delta}
\eq
where we have made use of inequalities~(\ref{eq:exp21}) and~(\ref{eq:Esup2}) 
\new
and checked the $n=0$ case separately. 
\endnew
Consequently
\be
\limsup_N\BP(w'((P^N,u^N),\delta,T)\ge\eta)\le\frc{4\delta(1+\hat\CQ(T+2)\sum_{n=0}^\infty2^{-\frac12n})}\eta\le\eta
\ee
whenever $\delta\le\frcn14\eta^2[1+\sqrt2\hat\CQ(T+2)/(\sqrt2-1)]^{-1}\wedge1$.
\end{proof}

\begin{corollary}\LABEL{cor:laws}\it The laws of the random variables $(P^N,u^N)$ are relatively compact in the space of probability measures defined on $D_E[0,\infty)$ equipped with the weak topology.
\end{corollary}
\begin{proof}
Note that since $\FR$ is polish so too is $E$. 
The result follows by applying proposition~3.7.4 of
\cite{EK},
together with corollary~\ref{cor:compact} and proposition~\ref{prop:modulus}.
\end{proof}

\begin{corollary}\LABEL{cor:cts}\it
There exists a subsequence $(P^{N_n},u^{N_n})$ converging in distribution (which we denote by $\imp$) to $(P,u)\in C_E[0,\infty)$ almost surely.
\end{corollary}
\begin{proof} Corollary~\ref{cor:laws} implies the existence of a subsequence converging to an element of $D_E[0,\infty)$. By theorem~3.10.2 \cite{EK} to prove the corollary it will suffice to
check that the distance between neighbouring states vanishes uniformly as $N\to\infty$. We have for any $(p,v)\in E^N$, $z,z_1,\ldots,z_\CI,w_1,\ldots,w_\CJ\in\FR$ and $k\in\BR^{d_1}$ that by using $\|f_n\|_\infty\le\:f_n\:\le\:\fext_n\|^\#\le2^n$ for all $n\in\BN_0$ that
\beq
D((p,v),J_{\rm CLK}^N(p,v))&=&\frc1{\rho^N(p,v)}\le\frc1{rN};\\
D((p,v),J^N_I((p,v),z))&\le&\sum_{n=0}^\infty\frc1{4^n}\(\left|\frc{f_n(z)}N\right|\wedge1\)+\frc1{\rho^N(p,v)}
\le\frc2N+\frc1{rN};\\
D((p,v),J^N_0((p,v),z,k))&=&\sum_{n=1}^\infty4^{-n}\(\left|\frc{f_n(\pim(z),\gamma(\pi_{\bar\Omega}(z),k),j_{\rm ID}^N(v,z))-f_n(z)}N\right|\wedge1\)
+\frc1{\rho^N(p,v)}
\cr&&\\
&\le&
\frc2N+\frc1{rN};\\
D((p,v),J_{i,j}^N((p,v),z_1,\ldots,z_i,w_1,\ldots, w_j))\kern-2.7cm&&\cr
&&\kern-3cm=\sum_{n=0}^\infty4^{-n}\(\left|\frc{\sum_{m=1}^jf_n(w_m)-\sum_{\ell=1}^if_n(z_\ell)}N\right|\wedge1\)+\frc1{\rho^N(p,v)}
\le\frc{2(\CI+\CJ)}N+\frc1{rN}.
\eq
\end{proof}

\begin{proposition}\LABEL{prop:bounded-mass}
For all $t\ge0$,
\be
\sup_{0\le s\le t}P_s(\pim)<\infty\qquad\mbox{a.s.}
\ee
\end{proposition}
\begin{proof}
Observe that $P^{N_n}_s(\pim\wedge L)\imp P_s(\pim\wedge L)$ by weak convergence, moreover, since by corollary~\ref{cor:cts} the limit is almost surely continuous in time the convergence is almost surely uniform on compact sets (theorem 3.10.1 in \cite{EK}). Thus $\sup_{0\le s\le t}P^{N_n}_s(\pim\wedge L)\imp\sup_{0\le s\le t}P_s(\pim\wedge L)$.
Now Fatou's lemma implies that
\be
\BE\(\sup_{0\le s\le t}P_s(\pim\wedge L)\)\le\liminf_n\BE\(\sup_{0\le s\le t}P^{N_n}_s(\pim\wedge L)\)\le\Xi +t\Lambda^{(1)}+2(t\Lambda^{(2)})^{\frac12}.
\ee
Taking $L\to\infty$ and using the monotone convergence theorem yields
\be
\BE\(\sup_{0\le s\le t}P_s(\pim)\)\le\Xi +t\Lambda^{(1)}+2(t\Lambda^{(2)})^{\frac12}.
\ee
Accordingly $\BP\(\sup_{0\le s\le t}P_s(\pim)<\infty\)=1$.
\end{proof}

\begin{definition}
Let  $\CV_t=\{(p,v)\in D_E[0,\infty):\sup_{0\le s\le t}p_s(\pim)<\infty\}$, then proposition~\ref{prop:bounded-mass} implies that $\BP(P\in\CV_t)=1$ for all $t\ge0$.
\end{definition}

\begin{remark}\em
Since $E$ is separable, so too is $D_E[0,\infty)$ by theorem~3.5.6 of \cite{EK}. Now the Skorohod representation theorem (theorem~3.1.8 {\em Ibid.\/}) implies that there exists a common probability space $(\b\Omega,\CF,\BP)$ on which $(P^{N_n},u^{N_n})$ are defined and $(P^{N_n},u^{N_n})(\varpi)\stackrel {\tilde D}\to(P,u)(\varpi)$ as $n\to\infty$ for all $\varpi\in\b\Omega$.
\end{remark}

\begin{proposition}\LABEL{prop:inception-limit}
If for all $t\ge0$, $(p^{n},v^{n})\in \CV_t\cap D_{E^{N_n}}[0,\infty)$ and $(p,v)\in C_E[0,\infty)\cap\CV_t$ with $p^{n}_s\weakto p_s$ and $\sup_{s\le t}|v^{n}_s-s|\to0$ as $n\to\infty$ then for every   $f\in\FF_c$,
(a)
\be
\int_0^t\int_\FR f(z)\openone_{p_s^{n}(\pim)+\pim(z)/N_n\le \aNn }I_{v^{n}_s}(\d z)\,\d s\to\int_0^tI_{s}(f)\,\d s
\ee
as $n\to\infty$ and (b)
\beq
&&\int_0^t\int_{\FR^{i+j}} K_{v^n_s,i,j}(z_1,\ldots,z_i,\d w_1,\ldots,\d
w_j)\left\{\sum_{m=1}^jf(w_m)-\sum_{n=1}^if(z_n)\right\}\nu^{N_n}_i(p^n_s,\d
z_1,\ldots,\d z_i)\,\d s\cr
&&\kern4cm\to\int_0^t\int_{\FR^{i+j}} K_{s,i,j}(z_1,\ldots,z_i,\d
w_1,\ldots,\d w_j)\cr
&&\kern6cm{}\times\left\{\sum_{m=1}^jf(w_m)-\sum_{n=1}^if(z_n)\right\}p_s(\d
z_1)\ldots p_s(\d z_i)\,\d s\LABEL{eq:coaglim}\qquad\quad
\eq
as $n\to\infty$.
\end{proposition}
\begin{proof}
(a) For all compactly supported functions $f$,
\be
f(z)\openone_{p^{n}_s(\pim)+\pim(z)/N_n>\aNn }=0\LABEL{eq:zero}
\ee
for large enough $n$,
as  by Fatou's lemma we have $\limsup_n p^n_s(\pim\wedge L)\le p_s(\pim\wedge L)\le p_s(\pim)$
 and thus the monotone convergence theorem implies that $\limsup_np_s^{n}(\pim)\le p_s(\pim)$. Equation~(\ref{eq:zero}) follows by the condition $\aNn \to\infty$.
Continuity of the map $t\mapsto I_t$ (hypothesis~\ref{hyp:cont}) implies that
$I_{v^n_s}(f)\to I_{s}(f)$
and thus
\beq
\left|\int_\FR f(z)\openone_{p^{n}_s(\pim)+\pim(z)/N_n\le \aNn } 
I_{v^{n}_s}(\d z)-I_{s}(f)\right|&\le&
\left|\int_\FR f(z)\openone_{p^{n}_s(\pim)+\pim(z)/N_n>\aNn } I_{v^{n}_s}(\d z)\right|\cr
&&\kern2cm{}+|I_{v^n_s}(f)-I_{s}(f)|\to0
\eq
as $n\to\infty$.
We know that $\left|\INT_\FR f(z)\openone_{p^{n}_s(\pim)+\pim(z)/N_n\le \aNn }
I_{v^n_s}(\d z)\right|\le\Lambda^{(0)}\|f\|_\infty$ and hence the dominated convergence theorem implies that
\be
\int_0^t\left[\int_\FR f(z) \openone_{p^{n}_s(\pim)+\pim(z)/N_n\le \aNn }I_{v^n_s}(\d z)-I_{s}(f)\right]\,\d s\to0\as n->\infty.
\ee
(b) We make separate arguments depending on whether we are considering the self-interaction kernel corresponding to $(i,j)=(1,1)$ or not. If $(i,j)\neq(1,1)$ we have by hypothesis~\ref{hyp:interactionker}
that the function
\be
\CK(s,z_1,\ldots,z_i)=\int_{\FR^{j}} K_{s,i,j}(z_1,\ldots,z_i,\d w_1,\ldots,\d
w_j)\left\{\sum_{m=1}^jf(w_m)-\sum_{n=1}^if(z_n)\right\}
\ee
is bounded in modulus by $(i+j)i!\,j!\,K_\infty\|f\|_\infty\pim(z_1)\cdots\pim(z_i)$. 
We assume that $n>n_0$ where
$\sup_{n>n_0}\sup_{0\le s\le t}v^n_s\le t+1$. Moreover since $\FR$ is
separable, $\CM_B(\FR^r)\cong\bigotimes_{i=1}^r\CM_B(\FR)$ and
$p^n_s\weakto p_s$ implies that $(p^n_s){}^{\otimes r}\weakto
p_s{}^{\otimes r}$. Then
if $f\in\FF_c$ and $r<i$ and $\alpha:\{1,\ldots,i\}\to\{1,\ldots,r\}$
we find
\beq
\frc1{N_n^{i-r}}\left|\int_{\FR^r}\CK(v^n_s,z_{\alpha(1)},\ldots,z_{\alpha(i)})p^n_s(\d z_1)\cdots p^n_s(\d z_r)\right|
&&\cr
&&\kern-5cm{}\le\frc{(i+j)i!\,j!}{N_n^{i-r}}\|f\|_\infty K_\infty R(f)^{i-r}\(\sup_{0\le s\le t}p_s^n(\pim)\)^{r}
\to0\qquad
\eq
as $n\to\infty$.
Recalling that
\beq
\nu_i^{N_n}(p^n_s,B_1\times\cdots\times B_i)&=&N_n^{-i}\sum_{\sigma\in S_{n}}\zeta(\sigma)\prod_{
\begin{array}{c}\scriptstyle
{\rm disjoint\ cycles}\\
\scriptstyle
{(i_1,\ldots,i_k){\rm\ of\ }\sigma}\end{array}
}
N_np^n_s(B_{i_1}\cap\cdots\cap B_{i_k})
\eq
it follows that
\be
\left|\int_{\FR^i}\CK(v^n_s,z_1,\ldots,z_i)\nu_i^{N_n}(p^n_s,\d z_1,\ldots,\d z_i)-\int_{\FR^i}\CK(v^n_s,z_1,\ldots,z_i)p^n_s(\d z_1)\cdots p^n_s(\d z_i)\right|\to0.
\ee
Furthermore we have
\beq
\left|\int_{\FR^i}\CK({v^n_s},z_1\ldots,z_i)\,p^n_s(\d z_1)\cdots p^n_s(\d z_i)-\int_{\FR^i}\CK(s,z_1,\ldots,z_i)\,p_s(\d z_1)\cdots p_s(\d z_i)
\right|&&
\cr
&&\kern-10cm{}\le\int_{\FR^i}
w(\CK,|v^n_s-s|)\,
p^n_s(\d z_1)\cdots p^n_s(\d z_i)\cr&&\kern-8cm{}+\left|\int_{\FR^i}\CK({s},z_1,\ldots,z_i)\,p^n_s(\d z_1)\cdots p^n_s(\d z_i)\right.\cr
&&\kern-5cm\left.{}-\int_{\FR^i}\CK({s},z_1,\ldots,z_i)\,p_s(\d z_1)\cdots p_s(\d z_i)\right|.\qquad
\eq
The first term on the LHS vanishes as $n\to\infty$ by the uniform continuity of $\CK$ (which follows by the uniform continuity given by hypothesis~\ref{hyp:interactionker}), \Eq.ctsmod., combined with the fact that $p^n_s(\FR)\to p_s(\FR)<\infty$, the second term vanishes by the aforementioned weak convergence of $(p^n_s){}^{\otimes i}$. Accordingly for each $s$,
\be
\int_{\FR^i}\CK({v^n_s},z_1,\ldots,z_i)\nu_i^{N_n}(p^n_s,\d z_1,\ldots,\d z_i)\to \int_{\FR^i}\CK({s},z_1,\ldots,z_i)\,p_s(\d z_1)\cdots p_s(\d z_i)
\ee
as $n\to\infty$. The LHS is c\`adl\`ag and therefore integrable in time, for large enough $n$  the LHS has modulus bounded above by
\be
\left|\int_{\FR^i}\CK(v^n_s,z_1,\ldots,z_i)\,p^n_s(\d z_1)\cdots p^n_s(\d z_i)\right|\le
(i+j)i!\,j!\,\|f\|_\infty K_\infty\(1+\sup_{0\le s\le t}p_s(\pim)^i\)
<\infty
\ee
then the dominated convergence theorem yields \Eq.coaglim..

The proof for the self-interaction kernel is similar, though in this case we remark that $\CK$ is supported on $[0,\infty)\times\supp f$ and therefore the uniform continuity of $\CK$ comes from the mass preserving property of the self-interaction kernel.
\end{proof}

\begin{definition}
Let $(P,u)$ be as in corollary~\ref{cor:cts} then for all $f\in\FF_c$ (or for all $f\in\FF$ if the linear growth condition is satisfied),
define
\beq
M_t(f)&=&P_t(f)-P_0(f)-\int_0^t I_s(f)\,\d s-\int_0^tP_s\(\h a_s\triangle f+b_s\cdot\D f\)\,\d s-\int_0^t P_s\(H_s\cdot\hat\D f\)\,\d s\cr
&&\kern1cm{}-\sum_{i=1}^\CI\sum_{j=1}^\CJ\frc1{i!\,j!}\int_0^t\int_{\FR^{i+j}}K_{s,i,j}(z_1,\ldots,z_i,\d w_1,\ldots,\d w_j)\cr
&&\kern4cm{}\times\(\sum_{m=1}^jf(w_m)-\sum_{n=1}^if(z_n)\)P_s(\d z_1)\ldots P_s(\d z_i)\,\d s.\LABEL{eq:M_t(f)}
\eq
\end{definition}

\begin{proposition}
\LABEL{prop:Mconv}
For all $f\in \FF_c$
and all $t\ge0$,
\be
M^{N_n}_t(\Phi_f)\imp M_t(f)\as n->\infty.
\ee
\end{proposition}
\begin{proof}
Since $(P^{N_n},u^{N_n})$ converges in distribution to $(P,u)$ with $\BP((P,u)\in C_E[0,\infty))=1$
we may prove this result by considering those sequences of sample paths with limit in $C_E[0,\infty)$.
 
By corollary~\ref{cor:cts} and the convergence in distribution of $u^{N_n}\imp u$, an almost surely continuous limit in time, theorem~3.10.1 in \cite{EK} implies that the convergence is almost surely
uniform on the compact set $[0,t]$, i.e., for almost all of the sample paths
\be
\sup_{0\le s\le t}|v^n_s-s|\to0\as n->\infty
\ee
where for convenience we have written $v^{n}=u^{N_n}(\varpi)$.

Let $p^{n}_s=P^{N_n}_s(\varpi)$, with $p^n_s\weakto p_s$. If $k=\frc{\sigma_{v^n_s}(z)Z}{\cNn }+\frc{b_{v^n_s}(z)}{\cNn ^2}$ with $Z\sim N(0,\Id1)$ then
\beq
&&\left|\cNn ^2
\int_\FR\BE^Z(f(\pim(z),\gamma(\pi_{\bar\Omega}(z),k),j^{N_n}_{\rm ID}(v^n_s,z))-f(z))\,p^n_s(\d z)
\right.\cr
&&\kern4cm\left.{}-\int_\FR(\h a_s(z)\triangle f(z)+b_s(z)\cdot\D f(z)+H_s(z)\cdot\hat\D f(z))\,p_s(\d z)
\right|\qquad\\
&&\le\left\|\cNn ^2\BE^Z(f(\pim(z),\gamma(\pi_{\bar\Omega}(z),k),j^{N_n}_{\rm ID}(v^n_s,z))-f(z))\right.\cr&&\left.\kern3cm{}-(\h a_s(z)\triangle f(z)+b_s(z)\cdot\D f(z)+H_s(z)\cdot\hat\D f(z))\right\|_\infty\sup_np^n_s(\FR)\qquad\quad\cr
&&\kern1cm{}+\left|\int_\FR(\h a_s(z)\triangle f(z)+b_s(z)\cdot\D f(z)+H_s(z)\cdot\hat\D f(z))\,p^n_s(\d z)\right.\cr
&&\left.\kern3cm{}-\int_\FR(\h a_s(z)\triangle f(z)+b_s(z)\cdot\D f(z)+H_s(z)\cdot\hat\D f(z))\,p_s(\d z)\right|.\LABEL{eq:unif-weak}
\eq
The first term on the RHS of \Ineq.unif-weak. converges to zero by theorem~\ref{thm:uniform-internal} and the weak convergence implying $p^n_s(\FR)\to p_s(\FR)$.
The second term tends to zero by the weak convergence $p^n_s\weakto p_s$ since the integrand is a bounded continuous function on $\FR$. The c\`adl\`ag property implies that
\be
\cNn ^2\int_\FR\BE^Z(f(\pim(z),\gamma(\pi_{\bar\Omega}(z),k),j^{N_n}_{\rm ID}(v^n_s,z))-f(z))\,p^n_s(\d z)
\ee
is integrable in time, and proposition~\ref{prop:boundedness} implies it is bounded in modulus by $LH^\infty_{R(f)}\:\fext\:^\#(1+p_s(\pim))$ for large enough $n$, and thus the dominated convergence theorem gives
\beq
\cNn ^2\int_0^t\int_\FR\BE^Z(f(\pim(z),\gamma(\pi_{\bar\Omega}(z),k),j^{N_n}_{\rm ID}(v^n_s,z))-f(z))\,p^n_s(\d z)\,\d s&&\cr
&&\kern-7cm\to\int_0^t\int_\FR(\h a_s(z)\triangle f(z)+b_s(z)\cdot\D f(z)+H_s(z)\cdot\hat\D f(z))\,p_s(\d z)\,\d s.\qquad
\eq
Proposition~\ref{prop:inception-limit} and the weak convergence results $p^n_t(f)\to p_t(f)$, $p^n_0(f)\to P_0(f)$ then complete the proof.
\end{proof}

\begin{proposition}
\LABEL{prop:bp}
Suppose $f\in \FF_c$ (or $f\in\FF$ if the linear growth condition holds),
$f=\bplim f_{n_k}$,
$\D f=\bplim\D f_{n_k}$, $\partial_i\partial_jf=\bplim\partial_i\partial_jf_{n_k}$,
 $\hat\D f=\bplim \hat\D f_{n_k}$ and $R(f_{n_k})\le R(f)$ (with $n_k\ge1$) whenever $f\in\FF_c$, then for all $t\ge0$,
\be
M_t(f_{n_k})\imp M_t(f)\as k->\infty.
\ee
\end{proposition}
\begin{proof}
Set $F=\sup_k\:f_{n_k}\:$. Then for every $(p,v)\in C_E[0,\infty)$
we have $|f_{n_k}(z)|\le F$, $p_t(F)<\infty$ and hence by the dominated convergence theorem $p_t(f_{n_k})\to p_t(f)$ as $k\to\infty$. Furthermore $I_s(F)\le F\Lambda^{(0)}$, so
that $I_s(f_{n_k})\to I_s(f)$ and therefore $\left|\int_0^t\,I_s(F)\,\d s\right|\le tF\Lambda^{(0)}$. A final application of the dominated convergence theorem
yields $\int_0^t I_s(f_{n_k})\,\d s\to \int_0^t I_s(f)\,\d s$, where we have made use of the continuity condition, hypothesis~\ref{hyp:cont}, to justify the integrability of the LHS with respect to Lebesgue measure.
 
Additionally by proposition~\ref{prop:bounded-mass} almost surely $(p,v)\in C_E[0,\infty)\cap\CV_t$, therefore we have for $s\in[0,t]$
\be
|\h a_s(z)\triangle f_{n_k}(z)+b_s(z)\cdot\D f_{n_k}(z)+H_s(z)\cdot\hat\D f_{n_k}(z)|\le (\h\sigma_\infty^2+b_\infty+1)FH^\infty_{R(f)}\pim(z)
\ee
when $f\in\FF_c$. We use the bound
\be
|\h a_s(z)\triangle f_{n_k}(z)+b_s(z)\cdot\D f_{n_k}(z)+H_s(z)\cdot\hat\D f_{n_k}(z)|\le (\h\sigma_\infty^2+b_\infty+1)FH_\infty\pim(z)
\ee
when the linear growth condition holds (recall remark~\ref{rem:lineargrowth}). To deal with both conditions at once set $H_\infty=H^\infty_{R(f)}$ when $f\in\FF_c$. Therefore
$p_s((\h\sigma_\infty^2+b_\infty+1)FH_\infty\pim)<\infty$ thus the
dominated convergence theorem implies that
\beq
&&p_s(\h a_s\triangle f_{n_k}+b_s\cdot\D f_{n_k}+H_s\cdot\hat\D f_{n_k})\cr
&&\kern4cm{}\to p_s(\h a_s\triangle f+b_s\cdot\D f+H_s\cdot\hat\D f)\\
&&\kern6cm{}\le (\h\sigma_\infty^2+b_\infty+1)F H_\infty\sup_{0\le s\le t}p_s(\pim)\qquad\quad
\eq
which is integrable on $[0,t]$ when $f\in\FF_c$.
and hence
\beq
\int_0^tp_s(\h a_s\triangle f_{n_k}+b_s\cdot\D f_{n_k}+H_s\cdot\hat\D f_{n_k})\,\d s
\to\int_0^t p_s(
\h a_s\triangle f+b_s\cdot\D f+H_s\cdot\hat\D f)\,\d s.
\eq
The integrability of the LHS comes from the continuity of
$p_s(\h a_s\triangle f_{n_k}+b_s\cdot\D f_{n_k}+H_s\cdot\hat\D f_{n_k})$.

For the interaction kernel terms with $(i,j)\neq(1,1)$ the dominated convergence theorem implies that
\beq
\int_{\FR^j}K_{s,i,j}(z_1,\ldots,z_i,\d w_1,\ldots,\d w_j)\left\{\sum_{m=1}^if_{n_k}(w_m)-\sum_{n=1}^jf_{n_k}(z_n)\right\}&&\cr
&&\kern-7cm{}\to\int_{\FR^j}K_{s,i,j}(z_1,\ldots,z_i,\d w_1,\ldots,\d w_j)\left\{\sum_{m=1}^if(w_m)-\sum_{n=1}^jf(z_n)\right\}\qquad
\eq
the LHS being bounded by $(i+j)K_\infty F\pim(z_1)\cdots\pim(z_i)$. In consequence a further application of the dominated convergence theorem implies that
\beq
\int_{\FR^{i+j}}K_{s,i,j}(z_1,\ldots,z_i,\d w_1,\ldots,\d w_j)\left\{\sum_{m=1}^if_{n_k}(w_m)-\sum_{n=1}^jf_{n_k}(z_n)\right\}\,p_s(\d z_1)\cdots p_s(\d z_i)&&\cr
&&\kern-13cm{}\to\int_{\FR^{i+j}}K_{s,i,j}(z_1,\ldots,z_i,\d w_1,\ldots,\d w_j)\left\{\sum_{m=1}^if(w_m)-\sum_{n=1}^jf(z_n)\right\}\,p_s(\d z_1)\cdots p_s(\d z_i).\qquad
\eq
We deduce that
\beq
\int_0^t\int_{\FR^{i+j}}K_{s,i,j}(z_1,\ldots,z_i,\d w_1,\ldots,\d w_j)\left\{\sum_{m=1}^if_{n_k}(w_m)-\sum_{n=1}^jf_{n_k}(z_n)\right\}\,p_s(\d z_1)\cdots p_s(\d z_i)\,\d s&&\cr
&&\kern-14.5cm{}\to\int_0^t\int_{\FR^{i+j}}K_{s,i,j}(z_1,\ldots,z_i,\d w_1,\ldots,\d w_j)\left\{\sum_{m=1}^if(w_m)-\sum_{n=1}^jf(z_n)\right\}\,p_s(\d z_1)\cdots p_s(\d z_i)\,\d s\qquad\quad
\eq
since the integrand on the LHS is continuous by hypothesis~\ref{hyp:interactionker} (and assuming that $p_s$ is continuous which is almost surely true) and is bounded in modulus by
$K_\infty F(i+j)\sup_{0\le s\le t}p_s(\pim)^i$
which in turn is integrable on $[0,t]$.

For the self-interaction kernel term $i=j=1$, we have for all $z\in\FR$ by the dominated convergence theorem,
\beq
\int_\FR K_{s,1,1}(z,\d w)\{f_{n_k}(w)-f_{n_k}(z)\}\to\int_\FR K_{s,1,1}(z,\d w)\{f(w)-f(z)\}
\eq
which is bounded by $2FK^\infty_{R(f)}$ if $f\in\FF_c$ and bounded by $2FK_\infty\pim(z)$ when the linear growth condition holds. It follows by the dominated convergence theorem that
\be
\int_{\FR^2} K_{s,1,1}(z,\d w)\{f_{n_k}(w)-f_{n_k}(z)\}p_s(\d z)\to\int_{\FR^2} K_{s,1,1}(z,\d w)\{f(w)-f(z)\}p_s(\d z).
\ee
Under either the linear growth condition or $f\in\FF_c$ the RHS is bounded and integrable in time, so that
\be
\int_0^t\int_{\FR^2} K_{s,1,1}(z,\d w)\{f_{n_k}(w)-f_{n_k}(z)\}p_s(\d z)\,\d s\to\int_0^t\int_{\FR^2} K_{s,1,1}(z,\d w)\{f(w)-f(z)\}p_s(\d z)\,\d s.
\ee
Putting these results together proves the proposition.
\end{proof}

\begin{proposition}\LABEL{prop:M_t}
Almost surely, for all $f\in\FF_c$ (or for all $f\in\FF$ if the linear growth condition holds) and for all $t\ge0$,  $M_t(f)=0$.
\end{proposition}
\begin{proof}
For each $t\in\BQ^+$ and $n\in\BN$ we have by proposition~\ref{prop:Mconv}, proposition~\ref{prop:two9} and Fatou's lemma
\be
\BE|M_t(f_n)|
\le\liminf_k\BE|M^{N_k}_t(\Phi_{f_n})|=0,
\ee
so that $M_t(f_n)=0$ almost surely. Observe that
\be
\BP(\{M_t(f_n)\neq0\mbox{ for some $n$ and $t\in\BQ^+$}\})=0
\ee
and since every $f\in\FF_c$ (or $f\in\FF$ if the the linear growth condition holds) is the limit in the sense of proposition~\ref{prop:bp} of a subsequence $f_{n_k}$ (with $n_k\neq0$) we have that $M_t(f)=0$ for all $t\in\BQ^+$. Recall from corollary~\ref{cor:cts} that almost all the sample paths of $(P,u)$ are continuous in time, and hence $M_t(f)=0$ for all $f$ and for all $t\ge0$ except on a set of zero probability.
\end{proof}

\begin{definition}\it Define the random variables
\be
v^N_t=\sup \ (u^{N})^{-1}[0,t]\and Q^{N}_t=P^{N}_{v^N_t}.
\ee
\end{definition}

\begin{theorem}\LABEL{thm:main} The sequence of random variables $Q^{N_k}\in D_{\CM_B(\FR)}[0,\infty)$ converges in distribution to $P$ as $k\to\infty$, where $D_{\CM_B(\FR)}[0,\infty)$ is equipped with the Skorohod metric  $\tilde d_{\rm weak}$ derived from $d_{\rm weak}$ using \Eq.metric0..
\end{theorem}
\begin{proof}
Observe that if $\sup_{0\le t\le T+1}|u_t^{N}-t|\le N^{-\frac14}$, then $\sup_{0\le t\le T}|v_t^{N}-t|\le2N^{-\frac14}$.
Now \Eq.skorohod. implies that
\be
\tilde d_{\rm weak}(Q^{N},P^{N})\le\int_0^\infty\e^{-T}\sup_{0\le t\le T}d_{\rm weak}(Q^{N}_t,P^{N}_t)\,\d T.
\ee
Since $\tilde d_{\rm weak}(Q^{N},P^{N})$ is a non-negative random variable bounded above by $\frcn43$ we have that
\beq
\BE\sup_{0\le t\le T}d_{\rm weak}(Q^{N}_t,P^{N}_t)&=&\BE\sum_{n=0}^\infty4^{-n}(\sup_{0\le t\le T}|P^{N}_{v^{N}_t}(f_n)-P^{N}_t(f_n)|\wedge1)\\
&&\kern-2cm{}\le\sum_{n=0}^\infty4^{-n}\(1\wedge\BE\sup_{0\le t\le T}\left|M^{N}_{v^{N}_t}(\Phi_{f_n})-M^{N}_t(\Phi_{f_n})+\int_t^{v_t^{N}}\CA^N(\Phi_{f_n})(P^N_s,u^N_s)\,\d s\right|
\)\cr&&\\
&&\kern-2cm
\le\sum_{n=0}^\infty4^{-n}\(2\BE\sup_{0\le t\le T+2}|M^{N}_t(\Phi_{f_n})|+2 N^{-\frac14}\hat\CQ(T+2)
K^\infty_{R(f_n)}H^\infty_{R(f_n)}\:\fext_n\:^\#\)\cr
&&\kern3cm{}+\frcn43\BP\(\sup_{0\le t\le T+1}|u^{N}_t-t|\ge N^{-\frac14}\)\\
&&\kern-2cm{}\le\frc{2\sqrt2}{\sqrt2-1}[(\CQ(T+2))^{\frac12}N^{-\frac12}
+N^{-\frac14}\hat\CQ(T+2)]+\frcn83(T+1)^{\frac12}r^{-\frac12}N^{-\frac14}
\eq
on using proposition~\ref{prop:two9}(d). Thus
\be
\BE(\tilde d_{\rm weak}(Q^{N},P^{N}))\to0\as N->\infty,
\ee
it follows that $\tilde d_{\rm weak}(Q^{N_k},P^{N_k})\to0$ in probability. Now apply corollary~3.3.3 and remark~3.3.4 of \cite{EK} to deduce that $Q^{N_k}\imp P$.
\end{proof}

\ssection{The Discrete Approximation Scheme}\label{sect:DAS}

Theorem~\ref{thm:main} combined with the construction of sample paths to the martingale problem for jump processes
(e.g., section~4.2 \cite{EK}) provide us with a technique for generating the sample paths of the processes $Q^N$. We start with an initial weak approximation $p^N_0\in\BS^N$ to $P_0$  at time $t=0$. Given the sample path generated up to and including time $t$, the procedure is then to compute $\tau^{-1}=\rho^N(p^N_t,t)$ from \Eq.rho^N. and set $p^N_s=p^N_t$ for $s\in(t,t+\tau)$. With probabilities
\be
\frc{\Lambda_t}{\rho^N(p_t^N,t)},\qquad
\frc{\cN ^2p^N_t(\FR)}{\rho^N(p_t^N,t)},\qquad
\frc{r}{\rho^N(p_t^N,t)},\qquad
\INT_{\FR^i}\frc{K_{t,i,j}(z_1,\ldots,z_i,\FR^j)\nu_i^N(p^N_t,\d z_1,\ldots,\d z_i)}{i!\,j!\,\rho^N(p_t^N,t)}
\ee
for $i=1,\ldots,\CI$ and $j=1,\ldots,\CJ$, determine $p^N_{t+\tau}$ by rules $(1)$, $(2)$, $(3)$ and $(i,j)$ respectively where rule (1) is to select an element $z\in\FR$ according to the law $I_t/\Lambda_t$ and let $p^N_{t+\tau}=p^N_t+N^{-1}\delta_z$ if this is in $\BS^N$, otherwise set $p^N_{t+\tau}=p^N_t$. Rule (2) is to choose an element of $\FR$ according to  $p^N_t/p^N_t(\FR)$ and compute its evolution under the diffusion, external and internal drifts: i.e., choose $Z$ according to the $N(0,\Id1)$ distribution and let $z'=(\pim(z),\gamma(\pi_{\bar\Omega}(z),\sigma_t(z)Z\cN ^{-1}+b_t(z)\cN ^{-2}),j_{\rm ID}^N(z))$ then let $p^N_{t+\tau}=p^N_t+N^{-1}\({\delta_{z'}-\delta_z}\)$. Rule (3) is to do nothing: $p^N_{t+\tau}=p^N_t$. Rule $(i,j)$ is to select $z_1,\ldots,z_i$ according to the law
\be
\frc{K_{t,i,j}(z_1,\ldots,z_i,\FR^j)\nu^N_i(p^N_t,\d z_1,\ldots,\d z_i)}
{\INT_{\FR^i} K_{t,i,j}(z_1,\ldots,z_i,\FR^j)\nu^N_i(p^N_t,\d z_1,\ldots,\d z_i)}
\ee
then choose $w_1,\ldots,w_j$ in accordance with the law
\be
\frc{K_{t,i,j}(z_1,\ldots,z_i,\d w_1,\ldots,\d w_j)}
{K_{t,i,j}(z_1,\ldots,z_i,\FR^j)}
\ee
and set $p^N_{t+\tau}=p^N_t+N^{-1}\(\Sum_{m=1}^j\delta_{w_m}-\Sum_{n=1}^i\delta_{z_n}\)$. In terms of computer implementation, the scheme may be augmented by the use of majorants and their associated fictitious jumps. The computation of efficient majorants for these processes will depend on the precise form of the interactions.

\ssection{Applications to Gas Phase Dynamics}

The approximation scheme presented in section~6 may be readily applied to the problem of the simulation of the formation of particle clusters in a reactor. In this section we point out some of the limitations of the method and to what extent they limit practical applications. We also point out how the rather general interaction kernels may be constructed to describe a selection of common phenomonena.

One of the major limitations of the current work is that hypothesis~\ref{hyp:interactionker} precludes a purely local interaction, i.e., where $\pi_{\bar\Omega}(z_1)=\ldots=\pi_{\bar\Omega}(z_i)=\pi_{\bar\Omega}(w_1)=\ldots=\pi_{\bar\Omega}(w_j)$, $K_{t,i,j}(z_1,\ldots,z_i,\d w_1,\ldots,\d w_j)$-a.e. Further restrictions on the diffusion, drift and interaction kernels need to be made to yield the corresponding existence results. In the literature, e.g., \cite{Guiasthesis,Gui01,Deaconu,Wrzosek}, the diffusion coefficients are usually assumed to be spatially homogeneous. This rules out some important physical applications, in particular we are not able to deal with those situations where the temperature of the gas may vary according to its position in the reactor vessel.

Included within the model defined by \Eq.model. is that discussed in~\cite{Turing}. In this paper Turing showed by analysing the Fourier series of the concentration that a diffusion-reaction system can exhibit catastrophic instability from a homogeneous equilibrium when non-homogeneous perturbations are introduced. The catastrophic instability refers to the divergence of the local particle concentration in finite time and arises when the magnitude of one of the Fourier coefficients becomes infinite. It follows from Parseval's theorem that the particle density is not square integrable for all time and a localized version of \Eq.model. with quadratic or higher order terms will not in general have a solution in the weak/vague  sense. For this reason it seems unlikely that solutions to a localized version of the model can be found except in a number of special cases. Let us mention a number of approaches to this problem. In \cite{Grosskinsky}, the existence and uniqueness of a suitably regular solution to the Smoluchowski coagulation equation was assumed, the stochastic particle approximation could then be compared to the solution and convergence established. However Norris~\cite{James} has shown that uniqueness need not hold for general coagulation kernels even in the homogeneous setting. Deaconu~\& Fournier~\cite{Deaconu} use a spatially inhomogeneous mass-flow technique to generate solutions to a mollified version of the pure coagulation equation. In this case the spatial independence of the diffusion coefficients (and absence of drift term) in the domain $\BR^{d_1}$ allows them to exploit the properties of the heat kernel to make an {\em a priori\/} estimate of the local particle density. This estimate is uniform in the quantity $\epsilon$ determining the mollification, and allows the $\epsilon\to0$ limit to be taken. The theorem is proved under conditions of uniform ellipticity on the diffusion terms. In the discrete case the coagulation kernel obeys a bound of the form  $K(m,m')\le C(m+m')$. Lauren\c cot \&~Mischler~\cite{LM01} have used weak compactness and a diagonal argument together with the results of \cite{Wrz97} to deduce an existence theorem under more general conditions, in particular fragmentation is taken into account and the uniform ellipticity condition is dropped.
 
In those cases where no existence result is known for the localized version of \Eq.model. we can make the following heuristic arguments from order of magnitude considerations. Supposing  a well-behaved limit exists it might be expected that our model approximates the localized version well on condition that  both the effective range of the particle interaction is large compared to the hydrodynamic scale ($\sim N^{-1/d_1}$) and tends to zero, i.e., the interaction localizes, as $N\to\infty$. In this way the interaction terms correctly detect the limiting particle density as $N\to\infty$.
From the point of view of the simulation of molecular dynamics 
it would be relevant to mollify the interactions 
on any macroscopic scale over which the limiting particle concentration has a suitably small variation, indeed the pointlike particle limit assumed in the derivation of \Eq.model. is itself unphysical, and would need to be modified in a way that prevents infinite particle densities from arising if the molecular dynamics were to be correctly described.

Now let us turn to the possible applications of the model. The individual particles can carry additional information that evolve with time and during interactions with other particles; these might include discretely changing quantities such as electric charge or the number of active sites, such as in the theory of soot formation arising from combustion theory, where free radicals on the surface of the soot particles lead to sites where binding to other particles can occur. The number of these sites typically decay with time. 
To incorporate this decay process into the model we include an internal coordinate taking values on a closed finite interval which contains all non-negative integers less than the largest possible number of active sites for a particle of the given size. The coagulation kernel will depend on the number of active sites of each of the potential coagulation partners and the self-interaction kernel will describe, possibly amongst other things, the rate at which the active sites decay.

Often it useful to have additional information about the morphology of the particles under investigation.  For instance the distribution of surface area is a valuable piece of information for certain chemical engineering applications. To account for this we take the particle's surface area as one of the internal coordinates and allow this to decay in time according to a size-dependent internal drift term. In this way sintering can be quantitatively described in the model.
 
The generality of the diffusion and drift terms in \Eq.model. allows us to incorporate the effect of a thermophoretic force, this is a force arising from the differential change in pressure acting on the particles due to a temperature gradient in the reactor. The thermophoretic force  manifests itself as a mass-dependent drift term in the modelling \Eq.model.. 

Chemical reactions may be modelled by labelling each of a finite number of chemical species by an integer, then, provided the chemical reactions proceeds by the law of mass action (if not it may be possible to write the reaction as a number of more elementary reactions that do obey the law of mass action) a suitable coagulation or fragmentation kernel can be used to implement the 
reactions.

\ssection{Conclusions}
In this paper we have constructed a purely stochastic scheme to solve a complicated partial differential integral equation. After imposing a number of hypotheses on the growth of the terms in the equation we were able to show the convergence, in a suitable sense, of this scheme to a solution of the equation. Our approach differs from many results in the mathematical literature in that we incorporate many possible physical processes rather than concentrating on just one or two. In addition we chose our hypotheses in such a way so as to maximize the physical usefulness of the resulting scheme; for this reason we are satisfied with using mollified interaction kernels if it means, for example, being able to model non-trivial diffusion, convection and source terms and realistic reactor geometries. Very few of the existing mathematical results allow for temporal inhomegeneity in the parameters describing the system. In contrast our result allows for such inhomogeneity, which may typically result from time dependent thermal variations in the reactor vessel.

The aim of this work has been to develop an approximation scheme and associated convergence result with wide applicability for the numerical study of molecular dynamics. In the scientific modelling literature it is relatively easy to write down a reasonable Monte Carlo scheme depending on many variables, though in those cases it is highly unlikely that any existing convergence theorem will apply directly. By proving the convergence for a sizable class of such models we have been able to put some of these schemes on a firm mathematical foundation.

One particular objective of this current work was to avoid the necessity of using a hybrid approach of Monte Carlo and finite element methods. A strength of the Monte Carlo approach is its superior scaling properties with the dimension of the configuration space compared to finite element methods. It is useful therefore to develop a purely stochastic approach to solving these problems. The author believes that there is considerable merit in having a number of different approaches to tackling such problems and the development of new schemes of all types is welcome, and all the more if the convergence of these schemes can be demonstrated.

\ssection{Acknowledgements}

I would like to thank James Norris, Alastair Smith and Neal Morgan 
for useful discussions.  

\aappendix{Separability of $\FF_c$}


\begin{proposition}\LABEL{prop:separable}
\it The space $(\FF_c,\:\cdot\:)$ defined in definition~\ref{def:testfunctions} is separable.
\end{proposition}
\begin{proof}
Let $\{p_n\}_{n=1}^\infty$ be the set of all polynomials on $\FR$ with rational coefficients and suppose that $G\in C^\infty(\BR)$ is a fixed smooth map satisfying $\openone_{(-\infty,0]}\le G\le\openone_{(-\infty,1]}$. Set
\be
\{g_k\}_{k=1}^\infty=\{z\mapsto G(\|z\|-r)p_n(z)\,:\  n,r\in\BN\},
\ee
and notice that  by using proposition~A.7.1 \cite{EK} and the $\:\cdot\:$-density of $(p_n)_{n=1}^\infty$ in the space of all polynomials on each compact set of $\FR$ that $\{g_k\}_{k=1}^\infty$ is a countable set of compactly supported functions that is dense in $(C_c^2(\FR),\:\cdot\:)$. For each $n,m\in\BN$, let $C_{n,m}=\{f\in\FF_c:\:f-g_m\:<n^{-1}\}$, and for every pair of integers $(n,m)$ such that $C_{n,m}\neq\emptyset$ choose $f_{n,m}\in C_{n,m}$, it follows that
\be
\{f_{n,m}:C_{n,m}\neq\emptyset,\ n,m\in\BN\}
\ee
 is a countable collection of functions that are dense in $(\FF_c,\:\cdot\:)$.
\end{proof}


\aappendix{Whitney Extension Result}

\begin{proposition}\LABEL{prop:line-segments}
 Suppose $x\in\partial\Omega$ and $y,z\in\bar\Omega$ such that $\|y-x\|,\|z-x\|<\frcn1{15}\|\D\D\omega\|^{-1}_\infty\wedge\delta$ where $\delta $ is such that for all  $w$ with $\|x-w\|<\delta$
we have that $\|\D\omega(w)\|>1$.  Let 
\be
p=y-2\|z-y\|\D\omega(y)/\|\D\omega(y)\|
\ee
then the line straight line segments from $y$ to $p$ and from $p$ to $z$ are contained within $\bar\Omega$.
\end{proposition}
\begin{proof}
For $\theta\in[0,1]$ we use Taylor's theorem to compute
\beq
\omega(y(1-\theta)+p\theta)&\le&\omega(y)-2\theta\|z-y\|\,\|\D\omega(y)\|+2\theta^2\|y-z\|^2\,\|\D\D\omega\|_\infty\\
&\le&-2\theta\|z-y\|(1-\|y-z\|\,\|\D\D\omega\|_\infty)\le0.
\eq
and
\beq
\omega(z(1-\theta)+p\theta)&\le&\omega(z)+\theta[(z-y)\cdot\D\omega(y)-2\|z-y\|\,\|\D\omega(y)\|]+3\theta\|y-z\|^2\|\D\D\omega\|_\infty\cr
&&\kern6cm{}+\frcn92\theta^2\|y-z\|^2\|\D\D\omega\|_\infty\\
&\le&-\theta\|y-z\|(1-\frcn{15}2\|y-z\|\,\|\D\D\omega\|_\infty)\le0
\eq
since $\|y-z\|\le\|y-x\|+\|x-z\|<\frcn2{15}\|\D\D\omega\|_\infty^{-1}$. This proves the proposition.
\end{proof}

\begin{proposition}\LABEL{prop:Whitney}
(a) If $f\in\FF_c\cupp\{f_0\}$ then there exists a function $\fext\in C^2(\FR^\#)$ such that the restriction satisfies $\fext|_{\FR}=f$ and  $\:\fext\:^\#$ is finite.

\ni(b) If $f\in\FF_c$ then $\fext$ can be chosen such that
$\max\pim(\supp(\fext))=R(f)$.
\end{proposition}
\begin{proof}
If $f=f_0=1$ then let $\fext_0=1$ and the proposition follows. Assume that $f\in\FF_c$.
The result will follow from the Whitney extension theorem (\cite{Whitney}, \cite{Abraham} appendix A, \cite{EK} proposition A.6.1).
For each $m\in\BN$ let $\FR_m=\bar\Omega\times\Gamma_m$, $\FR_m^\#=\co(\bar\Omega)\times\Gamma_m$ and let $F_m\in C^2(\FR_m)$ be defined by $F_m(\cdot)=f(m,\cdot)$. If $F_m=0$ then let the extension
$F_m^\#\in C^2(\FR^\#_m)$ be identically zero. In this way part (b) is satisfied. We note that $\:\fext\:^\#=\sup\{\|\fext(z)\|,\|\D\fext(z)\|,\|\hat\D\fext(z)\|,\|\D\D\fext(z)\|:z\in\FR^\#,\ \pim(z)\le R(f)\}$ is finite by the compact support condition and continuity.  Suppose that $F_m\neq0$, then for $u,v\in\FR_m$ construct
\beq
R^0(u,v)&=&F_m(u)-F_m(v);\\
R^1(u,v)&=&
\begin{cases}
\frc{F_m(u)-F_m(v)-F_m'(v)(u-v)}{\|u-v\|}&\mbox{if $u\neq v$};\\
0&\mbox{ if $u=v$;}\\
\end{cases}\\
R^2(u,v)&=&
\begin{cases}
\frc{F_m(u)-F_m(v)-F_m'(v)(u-v)-F_m''(v)(u-v)(u-v)}{\|u-v\|^2}&\mbox{if $u\neq v$};\\
0&\mbox{if $u=v$.}\\
\end{cases}
\eq
To ascertain the existence of an extension $F_m^\#\in C^2(\FR_m^\#)$ it suffices to show that each of these functions is (uniformly) continuous on the compact set $\FR_m\times\FR_m$.
The function $R^0$ is clearly continuous, as are $R^1$ and $R^2$ if $u\neq v$. Examine continuity at $(z,z)$. Write $\pi_{\bar\Omega}:\FR_m\to\bar\Omega$ for the projection onto $\bar\Omega$.
If $x=\pi_{\bar\Omega}(z)\not\in\partial\Omega$ then there exists a ball, $B$ of radius $r>0$ centred at $x$ contained within $\bar\Omega$.  Let $\CC=B\times\Gamma_m$ then $\CC$ is a convex
open set (in the topology on $\FR_m$) and it follows that if $u,v\in\CC$ with $u\neq v$ then
\be
|R^1(u,v)|=\left|\int_0^1\{F_m'((1-t)v+tu)-F_m'(v)\}\frc{(v-u)}{\|u-v\|}\,\d t\right|\le w(F_m',\|u-v\|)\to0
\ee
as $(u,v)\to(z,z)$ and where the modulus of continuity is defined by \Eq.modulus. using the operator norm. Similarly
\be
|R^2(u,v)|=\left|\int_0^1\{F_m''((1-t)v+tu)-F_m''(v)\}\frc{(v-u)}{\|u-v\|}\frc{(u-v)}{\|u-v\|}\,(1-t)\,\d t\right|\le \h w(F_m'',\|u-v\|)\to0.
\ee
If $x\in\partial\Omega$ then there may not exist a convex neighbourhood of $z$ contained in $\FR_m$ in which case we make use of the path given by proposition~\ref{prop:line-segments}.
Write $u=(y,Y)$ and $v=(w,W)$, with $y,w\in\bar\Omega$ and $Y,W\in \Gamma_m$, set $p=y-2\|y-w\|\D \omega(y)/\|\D \omega(y)\|$ and $\xi=(p,Y)$ then
we have
for $u\neq v$,
\beq
|R^1(u,v)|&=&\left|\int_0^1\left\{[F_m'((1-t)\xi+tu)-F_m'(v)]\frc{(u-\xi)}{\|u-v\|}+[F_m'((1-t)v+t\xi)-F_m'(v)]\frc{(\xi-v)}{\|u-v\|}\right\}\,\d t\right|\cr
&\le&5w(F_m',3\|u-v\|)\to0\as(u,v)->(z,z).
\eq
since $\|u-\xi\|\le2\|u-v\|$ and $\|\xi-v\|\le3\|u-v\|$. Similarly
\beq
|R^2(u,v)|&=&\left|\int_0^1\left\{(1-t)[F_m''((1-t)\xi+tu)-F_m''(v)]\frc{(u-\xi)}{\|u-v\|}\frc{(u-\xi)}{\|u-v\|}
\right.\right.\cr
&&\left.\left.\kern1cm{}+(1-t)[F_m''((1-t)v+t\xi)-F_m''(v)]\frc{(\xi-v)}{\|u-v\|}\frc{(\xi-v)}{\|u-v\|}
\right.\right.\cr
&&\left.\left.\kern3cm{}+[F_m''((1-t)v+t\xi)-F_m''(v)]\frc{(\xi-v)}{\|u-v\|}\frc{(u-\xi)}{\|u-v\|}\right\}\d t\right|\\
&\le&\frcn{25}2w(F_m'',3\|u-v\|)\to0\as(u,v)->(z,z).
\eq
It follows by the Whitney extension theorem that there exists $F_m^\#\in C^2(\FR_m^\#)$ such that $F_m^\#|_\FR=F_m$. Define $f^\#(m,z)=F_m^\#(z)$ to complete the proof.
\end{proof}

\gobble
\begin{proposition}\LABEL{prop:hash-bounded}
If $f\in\FF$ then $\fext$ may be chosen such that there exists a constant $C>0$ such that $\:\fext\:^\#\le C\:f\:$ (in particular $\:\fext\:^\#$ is finite).
\end{proposition}
\begin{proof} Since $\|\D\omega\|>1$ on $\partial\Omega$ which is compact there exist $\delta>0$ such that for all $x\in\bar\Omega$ such that the distance from $x$ to the boundary,
$d(x,\partial\Omega)<\delta$ implies that $\|\D\omega(x)\|>1$, we may take $\delta<\frcn1{15}\|\D\D\omega\|^{-1}_\infty$. By the compactness of $\partial\Omega
$ there exists a finite collection of open balls of radius $\delta$ in $\BR^{d_1
}$ with centres on $\partial\Omega
$ that cover $\partial\Omega
$, call these $B_1,\ldots,B_n$ and set $U_\mu=B_\mu\times\BR^{d_2}
$, additionally define 
$V=\(\BR^{d_1}\setminus\bar{\cupp_{\mu=1}^n B_\mu\cupp\bar\Omega}\)\times\BR^{d_2}$.

We let $(\psi_\mu)_{\mu=0}^{n}$ be a smooth partition of unity subordinate to $\{V,U_1,\ldots,U_n\}$. Then for all $y\in\FR^\#_m
$ we have either (a) $y\in\FR_m$, (b) $y\in V$  or (c) $y\in U_\mu$ for some $\mu=1,\ldots,n$. We define an extension of $F_m$ on each of $\FR_m$ (where $F_m^\#=F_m$),
 $V$ (with extension $F_{m,0}^\#=0$) and $U_\mu$ (with extension $F_{m,\mu}^\#$ defined below) then let $F_m^\#(y)=\sum_{\mu=0}^n\psi_\mu(y)F_{m,\mu}^\#(y)$ for
 $y\in\FR^\#_m\setminus\FR_m$. It suffices to show that on each $U_\mu$ that there exists some $C_\mu$ such that $\:F_{m,\mu}\:\le C_\mu\:f\:$ since if
 $y\in\FR_m$ then trivially $\:F_m^\#(y)\:=\:F_m(y)\:\le\:f\:$ and if $y\in V$ then $F^\#_{m,0}(y)=0$.

If $y\in U_\mu$ for some $\mu=1,\ldots,n$,
we make use of the proof of the Whitney extension theorem given in \cite{Abraham}. In particular there exists a sequence of smooth maps $\varphi_j$
such that $0\le\varphi_j(y)\le1$ and $\sum_j\varphi_j(y)=1$ for all $y\in U_\mu\setminus(\bar\Omega\times\BR^{d_2})$, and for each $\mu$, there exist constants $\alpha$, $N$, $M_1$ and $M_2$,
and a sequence of points $x_j\in\bar\Omega\times\BR^{d_2}\cap U_\mu$ such that $\|x_j-y\|\le\alpha d(y)$ whenever $\varphi_j(y)\neq0$ and
$\|\varphi_j'\|\le M_1d(y)^{-1}$, $\|\varphi_j''\|\le M_2d(y)^{-2}$
with $d(y)=\inf\{\|x-y\|:x\in\FR_m\}$.
Furthermore every point of $U_\mu\setminus(\bar\Omega\times\BR^{d_2})$ has a neighbourhood on which all but at most $N$ of the $\varphi_j$ are identically zero.

We define
\be
F_{m,\mu}^\#(y)=\sum_{j=1}^\infty\varphi_j(y)[F_m(x_j)+F_m'(x_j)(y-x_j)+\h F_m''(x_j)(y-x_j)(y-x_j)]
\ee
and note it obeys $|F_{m,\mu}(y)|\le(1+2\delta+2\delta^2)\:f\:$ whenever $y\in U_\mu\setminus(\bar\Omega\times\BR^{d_2})$. Observe that 
\be
\sum_{j=1}^\infty\sum_{k=1}^\infty F_m(x_k)\varphi_k(y)\varphi'_j(y)=0
\ee
 so that
\beq
F_{m,\mu}^\#{}'(y)&=&\sum_{j=1}^\infty\varphi_j(y)[F_m'(x_j)+F_m''(x_j)(y-x_j)]\cr
&&\kern-1cm{}+\sum_{j,k=1}^\infty[F_m(x_j)-F_m(x_k)+F_m'(x_j)(y-x_j)
+\h F_m''(x_j)(y-x_j)(y-x_j)
]\varphi_k(y)\varphi_j'(y).\qquad
\eq
Let $x_j=(q_j,Q_j)$ and set $\xi_{jk}=(p_{jk},Q_j)$ where $p_{jk}=q_j-2\|q_j-q_k\|\D\omega(q_j)/\|\D\omega(q_j)\|$, then using proposition~\ref{prop:line-segments} we have
\be
F_m(x_j)-F_m(x_k)=\int_0^1\{F_m'((1-t)\xi_{jk}+tx_j)(x_j-\xi_{jk})+F_m'((1-t)x_k+t\xi_{jk})(\xi_{jk}-x_k)\}\,\d t.
\ee
Now $|F_m(x_j)-F_m(x_k)|\le5\|x_j-x_k\|\,\:f\:\le10\alpha d(y)\:f\:$ if both $\varphi_k(y)$ and $\varphi_j'(y)$ are non-zero. Thus
\be
\|F_{m,\mu}^\#{}'(y)\|\le[1+2\delta+11\alpha NM_1+\alpha NM_1\delta]\:f\:.
\ee
Similarly we have $\|F_m'(x_j)-F_m'(x_k)\|\le10\alpha d(y)\:f\:$ and
\beq
F_m(x_j)-F_m(x_k)+F_m'(x_j)(y-x_j)-F_m'(x_k)(y-x_k)\kern-8.5cm&&\cr
&=&\int_0^1\left\{(1-t)[F_m''((1-t)x_k+t\xi_{jk})(x_k-\xi_{jk})(x_k-\xi_{jk})
-F_m''((1-t)\xi_{jk}-tx_j)(\xi_{jk}-x_j)(\xi_{jk}-x_j)
\right.\cr
&&\left.\kern0cm{}+[F_m''((1-t)\xi_{jk}+tx_j)(x_j-\xi_{jk})-F_m''((1-t)\xi_{jk}+tx_k)(x_k-\xi_{jk})](y-\xi_{jk})\right\}\d t,
\eq
leading to
\be
\|F_m(x_j)-F_m(x_k)+F_m'(x_j)(y-x_j)-F_m'(x_k)(y-x_k)\|\le38\alpha^2d(y)^2\:f\:
\ee
if neither $\varphi_j''(y)$ nor $\varphi_k(y)$ are zero. Thus we have
\beq
F_{m,\mu}^\#{}''(y)(H,K)&=&\sum_{j=1}^\infty\varphi_j(y)F_m''(x_j)(H,K)\cr
&&\kern-0.5cm{}+\sum_{j,k=1}^\infty\varphi_j'(y)(H)\varphi_k(y)[F_m'(x_j)-F_m'(x_k)
+F_m''(x_j)(y-x_j)](K)\cr
&&\kern-0.5cm{}+\sum_{j,k=1}^\infty\varphi_j'(y)(K)\varphi_k(y)[F_m'(x_j)-F_m'(x_k)+F_m''(x_j)(y-x_j)](H)\cr
&&\kern-0.5cm{}+\sum_{j,k=1}^\infty\varphi_j''(y)(H,K)\varphi_k(y) [F_m(x_j)-F_m(x_k)+F_m'(x_j)(y-x_j)-F_m'(x_k)(y-x_k)\cr
&&\kern5cm{}+\h F_m''(x_j)(y-x_j)(y-x_j)]
\eq
so that
\be
\|F_{m,\mu}^\#{}''(y)\|\le[1+22\alpha NM_1+\frcn{77}2\alpha^2NM_2]\:f\:
\ee
and the result follows from the boundedness of $\psi_\mu$ and its first and second derivatives.
\end{proof}

\endgobble

\aappendix{Particle Selection Measures}

\begin{proposition}\LABEL{prop:particleselectionmeasures}
If $P=\sum_{i=1}^N\delta_{z_i}$ then  $\nu_n(P,B_1\times\cdots\times B_n)$ is the number of injections $\alpha:\{1,\dots,n\}\to\{1,\ldots,N\}$ such that $(z_{\alpha(1)},\ldots,z_{\alpha(n)})\in B_1\times\cdots\times B_n$.
\end{proposition}
\begin{proof}
A map $\beta:\{1,\dots,n+1\}\to\{1,\dots,N\}$ with $z_{\beta(r)}\in B_r$ for $r=1,\ldots,n+1$ and $\alpha=\beta|_{\{1,\ldots,n\}}$ an injection, is either an injection itself, or there exists a unique $j\in\{1,\ldots,n\}$ such that $\alpha(j)=\beta(n+1)$ with $z_{\alpha(j)}\in B_j\cap B_{n+1}$. Accordingly
\be
\nu_{n+1}(P,B_1\times\cdots\times B_{n+1})=\nu_n(P,B_1\times\cdots\times B_n)P(B_{n+1})-\sum_{j=1}^n\nu_n(P,B_1\times\cdots\times(B_j\cap B_{n+1})\times\cdots\times B_n)\qquad
\ee
with $\nu_1(P,B_1)=P(B_1)$.
We show that \Eq.nu_n. satisfies this equation. Clearly equality holds for $\nu_1$. We
regard  $S_n$ as the stabilizer of $n+1$ in $S_{n+1}$. Note that for any $\sigma\in S_{n+1}$ either $\sigma\in S_n$ or $\sigma=(n+1,j)\sigma'$ with $\sigma'\in S_n$ and $j=\sigma(n+1)\in\{1,\ldots, n\}$.
Thus
\beq
\sum_{\sigma\in S_{n+1}}\zeta(\sigma)\kern-10pt\prod_{
\begin{array}{c}\scriptstyle
{\rm disjoint\ cycles}\\
\scriptstyle
{(i_1,\ldots,i_k){\rm\ of\ }\sigma}\\\scriptstyle{\rm\ in\ }S_{n+1}
\end{array}
}\kern-10pt
P(B_{i_1}\cap\ldots\cap B_{i_k})
&=&\sum_{\sigma\in S_{n}}\zeta(\sigma) P(B_{n+1})\kern-10pt\prod_{
\begin{array}{c}\scriptstyle
{\rm disjoint\ cycles}\\
\scriptstyle
{n+1\notin(i_1,\ldots,i_k){\rm\ of\ }\sigma}\\\scriptstyle{\rm\ in\ }S_{n}
\end{array}
}\kern-10pt
P(B_{i_1}\cap\ldots\cap B_{i_k})\cr
&&\kern-6cm{}+\sum_{j=1}^n\kern-0.5cm\sum_{\begin{array}{c}\scriptstyle\sigma'\in S_n
\\
\scriptstyle
\sigma=(n+1,j)\sigma'
\end{array}}\kern-0.5cm\zeta(\sigma)
P((B_{n+1}\cap B_j)\cap B_{\sigma'(j)}\cap\ldots\cap B_{\sigma'^n(j)})
\kern-10pt\prod_{
\begin{array}{c}\scriptstyle
{\rm disjoint\ cycles}\\
\scriptstyle
{j\notin(i_1,\ldots,i_k){\rm\ of\ }\sigma'}\\\scriptstyle{\rm\ in\ }S_{n}\end{array}
}\kern-10pt
P(B_{i_1}\cap\ldots\cap B_{i_k})\\
&&\kern-6cm{}=
\nu_n(P,B_1\times\cdots\times B_n)P(B_{n+1})
-\sum_{j=1}^n
\nu_n(P,B_1\times\cdots\times(B_j\cap B_{n+1})\times\cdots\times B_n)\\
&&\kern-6cm{}=\nu_{n+1}(P,B_1\times\cdots\times B_{n+1}).
\eq
\end{proof}

\gdef\refname{{\normalsize\bf References}}

\end{document}